\documentclass{book}
\usepackage[margin=2.7cm]{geometry}
\usepackage{afterpage}
\usepackage{times}
\usepackage[T1]{fontenc}

\usepackage{graphicx}
\usepackage{makeidx}
\usepackage{multicol}
\usepackage{subeqn}
\usepackage{url}
\usepackage{listings}
\usepackage{bm}
\usepackage{amsfonts}
\usepackage{amssymb}
\usepackage{booktabs}
\usepackage{xcolor}
\usepackage{color,colortbl}
\usepackage[pageref]{backref}
\usepackage{float}
\usepackage{pgfplots}
\usepackage{pgfplotstable}
\usepackage[multiple]{footmisc}
\usepgfplotslibrary{groupplots}
\usepackage{tikz}
\usetikzlibrary{calc,trees,positioning,arrows,chains,shapes.geometric,%
    decorations.pathreplacing,decorations.pathmorphing,shapes,%
    matrix,shapes.symbols, decorations.markings, patterns,fit,quotes,%
    arrows.meta}

\renewcommand*{\backrefalt}[4]{%
\ifcase #1 %
\color{white}{ }
\or
\color{white}{ }
\else
\color{white}{ }
\fi
}

\makeatletter
\let\@fnsymbol\@arabic
\makeatother

\definecolor{codegreen}{rgb}{0,0.6,0}
\definecolor{codegray}{rgb}{0.5,0.5,0.5}
\definecolor{codepurple}{rgb}{0.58,0,0.82}
\definecolor{backcolour}{rgb}{0.95,0.95,0.92}

\lstdefinestyle{mystyle}{
    backgroundcolor=\color{backcolour},
    commentstyle=\color{codegreen},
    keywordstyle=\color{magenta},
    numberstyle=\tiny\color{codegray},
    stringstyle=\color{codepurple},
    basicstyle=\ttfamily\footnotesize,
    breakatwhitespace=false,
    breaklines=true,
    captionpos=b,
    keepspaces=true,
    numbers=left,
    numbersep=5pt,
    showspaces=false,
    showstringspaces=false,
    showtabs=false,
    tabsize=2
}

\ifdefined\norm\else
\newcommand{\norm}[1]{\Vert #1 \Vert}
\fi

\usetikzlibrary{chains,positioning,shapes.symbols, backgrounds}
\usetikzlibrary{shapes.geometric, arrows}
\usetikzlibrary{shapes.callouts}

\providecommand{\eqref}[1]{(\ref{#1})}

\ifdefined 
  \ifpdf
  \resizebox{1\textwidth}{!}{
    \input{\else.pdftex_t}
  }
  \else
  \resizebox{1\textwidth}{!}{
    \input{\else.pstex_t}
  }
  \fi

\fi
\usepackage{algorithmic}

\DeclareSymbolFont{mathx}{U}{mathx}{m}{n}
\DeclareMathSymbol{\bigtimes}{1}{mathx}{"91}

\usepackage{siunitx}
\usepackage{tikzsymbols}

\usetikzlibrary{calc,trees,positioning,arrows,chains,shapes.geometric,%
    decorations.pathreplacing,decorations.pathmorphing,shapes,%
    matrix,shapes.symbols, decorations.markings, patterns,fit}

\tikzset{
  punktchain/.style={
    rectangle,
    rounded corners,
    draw=black, very thick,
    text width=25em,
    minimum height=3em,
    text centered,
    on chain},
  line/.style={draw, thick, <-},
  element/.style={
    tape,
    top color=white,
    bottom color=blue!50!black!60!,
    minimum width=8em,
    draw=blue!40!black!90, very thick,
    text width=10em,
    minimum height=3.5em,
    text centered,
    on chain},
  every join/.style={->, thick,shorten >=1pt},
  tuborg/.style={decorate},
  tubnode/.style={midway, right=2pt},
}

\DeclareRobustCommand\sampleline[1]{%
  \tikz\draw[#1] (0,0) (0,\the\dimexpr\fontdimen22\textfont2\relax)
  -- (1.5em,\the\dimexpr\fontdimen22\textfont2\relax);%
}

\usepackage{pgfplots}
\pgfplotsset{compat=1.17}

\ifdefined 
  \ifpdf
  \resizebox{1\textwidth}{!}{
    \input{./images/\else.pdftex_t}
  }
  \else
  \resizebox{1\textwidth}{!}{
    \input{./images/\else.pstex_t}
  }
  \fi

\fi

\usepackage{comment}

\usepackage{array,multirow}
\newcommand{\ldb}{\mathopen{\lbrack\!\lbrack}}
\newcommand{\rdb}{\mathclose{\rbrack\!\rbrack}}

\usepackage{tikz-cd}

\ifdefined 
  \ifpdf
  \resizebox{1\textwidth}{!}{
    \input{\else.pdftex_t}
  }
  \else
  \resizebox{1\textwidth}{!}{
    \input{\else.pstex_t}
  }
  \fi

\fi

\usepackage{booktabs}

\newcommand{\SISSA}{\footnote[2]{mathLab, Mathematics Area, SISSA, Scuola Internazionale Superiore di Studi Avanzati, Trieste, Italy; \\{\color{white}{XXX}}email: gstabile@sissa.it, grozza@sissa.it}}
\newcommand{\SISSAtwo}{$^2$}
\newcommand{\unicatt}{\footnote[3]{Department of Mathematics and Physics, Catholic University of the Sacred Heart, Brescia, Italy; \\{\color{white}{XXX}}email: francesco.ballarin@unicatt.it}}

\newcommand{\uniat}{\footnote[1]{School of Mathematics, Aristotle University of Thessaloniki, Thessaloniki, Greece; \\{\color{white}{XXX}}email: ekaratza@math.auth.gr}}

\graphicspath{{images/}}

\begin{document}

\chapter{Reduced Basis, Embedded Methods and Parametrized Levelset  {Geometry}}
\pagestyle{plain}
%
\begin{author}
{  }Efthymios N. Karatzas\uniat, Giovanni Stabile\SISSA, Francesco Ballarin\unicatt, Gianluigi Rozza\SISSAtwo
\end{author}

\section{Introduction and Overview}
In this chapter we examine reduced order techniques for geometrical parametrized heat exchange systems, Poisson, and flows based on Stokes, steady and unsteady incompressible Navier-Stokes and Cahn-Hilliard problems. The full order finite element methods, employed in an embedded and/or immersed geometry framework, 
are  the Shifted Boundary (SBM) and the Cut elements (CutFEM) methodologies, with applications mainly focused in fluids. We start by introducing the Nitsche's method, for both SBM/CutFEM and parametrized physical problems as well as the high fidelity approximation. We continue with the full order parameterized Nitsche shifted boundary variational weak formulation, and the reduced order modeling ideas based on a Proper Orthogonal Decomposition Galerkin method and geometrical parametrization, quoting the main differences and advantages with respect to a reference domain approach used for classical finite element methods, while stability issues may overcome employing supremizer enrichment methodologies.
Numerical experiments verify the efficiency of the introduced ``hello world'' problems considering reduced order results in several cases for one, two, three and four dimensional geometrical kind of parametrization. We investigate execution times, and we illustrate transport methods and improvements. A list of important references related to unfitted methods and reduced order modeling are  \cite{KaBaRO18,karatzas_stabile_atallah,KaratzasStabileNouveauScovazziRozza2018,KaratzasStabileNouveauScovazziRozzaNS2019,KaratzasRozzaCH20,KaNoBaRo20,KKF22}.
 \subsection{ {The heat exchange model problem}}
  {Let us} consider the simplest problem, Poisson's equation:
\begin{eqnarray*}
-\Delta T = g \quad\text{in }{\mathcal D}, \text{ with }
 T = g_D \quad\text{on }\Gamma. 
\end{eqnarray*}
We are looking for a variational formulation that it is satisfied by the weak solution $T\in H^1({\mathcal{D}})$, i.e., consistent, that is symmetric, and has a unique solution--the bilinear form is coercive.
We start by taking the strong form of the equation, multiplying by a test function $v\in H^1({\mathcal{D}})$ and integrating by parts. Starting with the right-hand side, we add the productive zero $0=T-g_D$ on the boundary and we obtain,
$$
(g ,v)_{\mathcal D} = (-\Delta T,v)_{\mathcal D}=(\nabla T,\nabla v)_{\mathcal D} - \int_{\Gamma} {\bf n}\cdot{\nabla} T v\,ds   
{\,- \int_{\Gamma
} (T-g_D){\bf n}\cdot{\nabla}v\,ds.}
$$
Separating linear and bilinear forms, $\text{for all }v\in H^1({\mathcal D})$, the variational equation --symmetric bilinear formulation-- becomes,  
$$
(\nabla T 
,\nabla v)_{\mathcal D} - \int_{\Gamma 
} {\bf n}\cdot{\nabla} 
T  
v\,ds - {\int_{\Gamma 
} T 
 {\bf n}\cdot{\nabla} 
v\,ds }  
= -{\int_{\Gamma 
} g_D{\bf n}\cdot{\nabla} 
v\,ds }
+ \int_{\mathcal D} g v\,dx
.
$$
 {We note that in the above formulation, the symmetric bilinear form is  not coercive due to that it cannot be bounded} from below for $T
=v$ by $c\|v\|^2 _{H^1({\mathcal D})}$. So, we  add a symmetric term that vanishes for the true solution: $\eta\int_{\Gamma} ( {T-g_D})v\,ds$ for some $\eta>0$ large enough. This leads to the symmetric, consistent, coercive weak formulation: 
Find $T\in H^1({\mathcal D})$, $\forall v\in H^1({\mathcal D})$ such that
\begin{eqnarray}
(\nabla T
,\nabla v)_{\mathcal D} - \int_{\Gamma
} {\bf n}\cdot{\nabla} 
T 
v\,ds - {\int_{\Gamma
} T
{\bf n}\cdot{\nabla} 
v\,ds } +{\eta\int_{\Gamma
} T
v\,ds }= -{\int_{\Gamma
} g{\bf n}\cdot{\nabla}
 v\,ds} 
 \nonumber
\\
  + {\eta\int_{\Gamma
} g_D v\,ds} + \int_{\mathcal D} g v\,dx.\quad
\label{eq:Nitsche form}
\end{eqnarray}
For finite element discrete approximations $T
_h,v_h\in V_h\subset H^1(\Omega)
$ 
and for Nitsche penalty chosen as 
$\eta= ch^{-1}$ with $c>0$ sufficiently large, one can manage a stable discrete problem with respect to a suitable mesh-dependent norm.
\subsection{Shifted Nitsche boundary weak formulation}
\begin{figure} \centering
(i) CutFEM 
\\
 \includegraphics[width=0.99\textwidth]{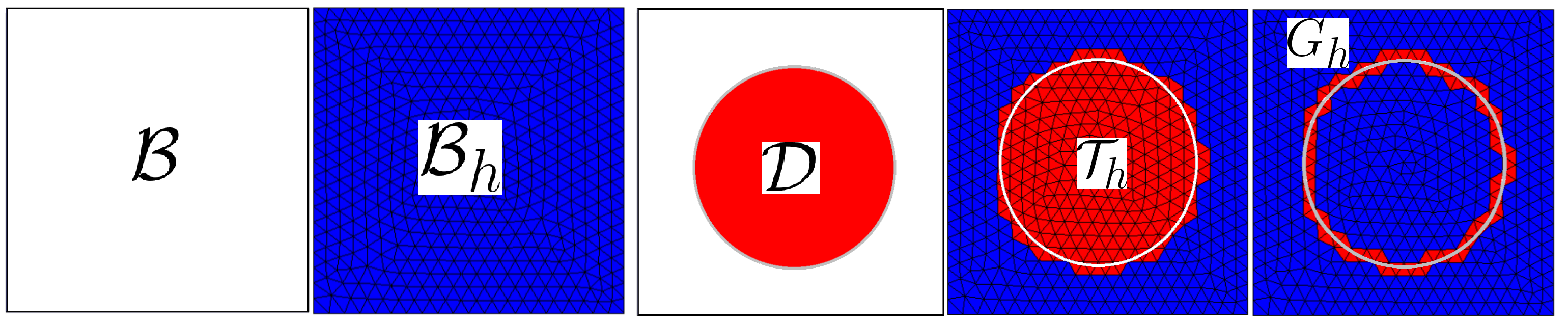} 
  \\
(ii) SBM
\\
  \includegraphics[width=0.99\textwidth]{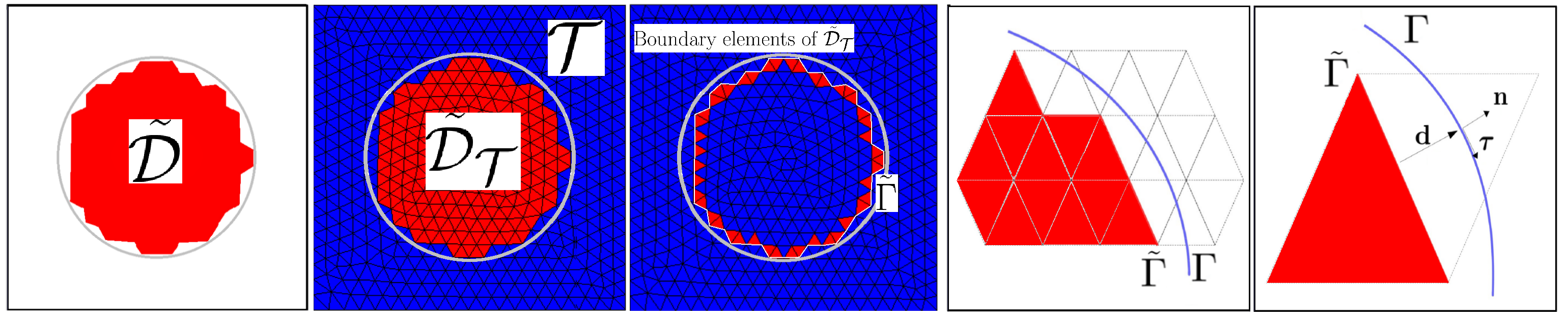}
\caption{Embedded FEM geometrical tools for (i) CutFEM and (ii) SBM. 
}
\label{SurrogateMesh}
\end{figure}

\index{Nitsche boundary weak formulation}
Starting from the aforementioned Nitsche form (\ref{eq:Nitsche form}) and based on a closest-point projection, see e.g. \cite{karatzas_stabile_atallah} and references therein, we derive a segmented/faceted nature of the employed surrogate boundary. In particular a smooth mapping  {$\bf M$ from points in a surrogate boundary ${\tilde{\Gamma}}$ to points in the true boundary} $ \Gamma $ is introduced:
$
   {\bf{M}}: {\tilde{\bf x}}|_{ {\tilde{\Gamma}
}}\to {\bf x}|_{\Gamma
}$, 
while the mapping ${\bf{M}}$ is defined though a distance vector function
${\bf d \equiv {d_M (x) = x}} {- {\tilde{\bf x}} = [{\bf M - I}]( {\tilde{\bf x}})}$, 
 see Figure \ref{SurrogateMesh}(ii).
Finally, we conclude to the Shifted Boundary weak formulation:
\begin{eqnarray*}
(\nabla T
,\nabla v)_{{{\text{$\tilde {{\mathcal D}}$}}}} 
- \hskip-2pt\int_{{{\text{$\tilde\Gamma
$}}}}{\tilde{\bf{n}}}\cdot{\nabla}
T
 v\,ds
 - 
 {\int_{{{\text{$\tilde {{\Gamma
}}$}}}} T
{\tilde{\bf{n}}}\cdot{\nabla} 
v\,ds 
} 
+{\eta\int_{{{\text{$\tilde\Gamma
$}}}} T
 v\,ds 
}= 
 - {\int_{{{\text{$\tilde\Gamma
$}}}} {{{\text{${\tilde{g}}_D
$}}}}{{\tilde{\bf{n}}}}\cdot{\nabla}
 v\,ds
 } 
\\ 
+ {\eta\int_{{{\text{$\tilde\Gamma
$}}}} {{{\text{${\tilde{g}}_D
$}}}} v\,ds
 } + \int_{{{\text{$\tilde{\mathcal D}$}}}} g v\,dx
,
\end{eqnarray*}
or  {equivalently}
\begin{eqnarray*}
(\nabla T
,\nabla v)_{{{\text{$\tilde {{\mathcal D}}$}}}} 
- \hskip-2pt\int_{{{\text{$\tilde\Gamma
$}}}}{\tilde{\bf{n}}}\cdot{\nabla} 
T
 v\,ds
 - 
 {\int_{{{\text{$\tilde {\Gamma
}$}}}} T
{\tilde{\bf{n}}}\cdot{\nabla} 
v\,ds 
} 
+{\eta\int_{{{\text{$\tilde\Gamma
$}}}} T
 v\,ds 
}
\qquad \qquad\qquad \qquad\\ 
= 
 - {\int_{{{\text{$\tilde\Gamma
$}}}} {{{\text{$(T
+\nabla T
 \cdot {\bf{d}})$}}}} {{\tilde{\bf{n}}}}\cdot{\nabla}
 v\,ds
 }
 + {\eta\int_{{{\text{$\tilde\Gamma
$}}}} {{{\text{$(T
+\nabla T
 \cdot {\bf{d}})$}}}} v\,ds
 } + \int_{{{\text{$\tilde{\mathcal D}$}}}} g v\,dx,
\end{eqnarray*}
where 
${{\tilde{g}}_D}$ denotes an extension of the Dirichlet boundary condition $ g_D$ to the boundary  $\tilde{\Gamma}$ of the surrogate domain 
based on a second order accurate Taylor expansion, {$T + \nabla T\cdot {\bf{d}} \approx \tilde{g}_D$}  {or if we let ${\bf  x}={\bf M}(\tilde {\bf  x})$:
 $T(\tilde{\bf  x}) +\nabla T(\tilde {\bf  x})\cdot{\bf{d}}
\approx g_D({\bf  x})$}.
\subsection{The parametrized thermal-heat exchange model}\index{thermal-heat exchange model, parametrized}
We consider parameter depended geometry with a parameter which lives in a $k$-dimensional parameter space $\mathcal P$, 
 and parameter vector ${\mu}\in \mathcal P \subset \mathbb{R} ^k
$.  
We denote by {${\mathcal{D}}(\mu)$}  
a bounded parametrized domain depending on $\mu$, with boundary $\Gamma(\mu)$. 
So, the model problem in $\mathcal D(\mu)$ becomes: find the temperature 
$T(\mu)$
\begin{eqnarray} \label{P:Poisson}
-\Delta  T(\mu)=g(\mu)
\quad\text{in }
\mathcal{D}(\mu), \text{ with }
T(\mu)=g_D(\mu) \quad   \text{on }  \Gamma
(\mu),
\end{eqnarray}
that in SBM weak formulation 
can be transformed in a system of linear equations (rewritten in matrix form):\begin{equation}\label{eq:system_linear}
{\bf
{A}}(\mu) {\bf{T}}(\mu)   = {\bf{F_g}}(\mu) \mbox{,}
\end{equation}
for 
 ${\bf
{A}}(\mu)\in {\mathbb R}^{N_h \times N_h}$ and $ {\bf{F_g}}(\mu)\in {\mathbb{R}}^{N_h \times 1}$ to correspond to the bilinear and linear form respectively. 
\subsection{Model Reduction Methodology}
{\textit{Offline stage.}}
 In this stage {, which is called ``training'',} one performs a certain number of full order solves for various parameters in order to use the solutions for the construction of a low dimensional reduced basis.
   This basis can approximate any member of the solution set to a prescribed accuracy, 
  while it is possible to perform a 
  {Galerkin projection} of the full order differential operators, describing the governing equations, onto the reduced basis space. 
  This procedure involves the solution of a possibly large number of high dimensional problems and the manipulation of high-dimensional structures.
  The required computational cost is high and therefore this operation is usually performed on a high performance system such as a computer cluster.
\\
{\textit{Online stage.}}
 Into this stage and often on a system with a reduced computational power and storage capacity, the reduced system 
 can be solved for any new value of the input parameters with predicted accuracy and reduced computational and time cost. For the interested reader we also refer to \cite{rozza2015book}.

\subsection{POD}\index{Proper Orthogonal Decomposition}\label{POD}
 Subsequently, with a Proper Orthogonal Decomposition one can generate the reduced basis space after 
the full-order model has been solved for each $\mu \in \mathcal{K}=\{ \mu^1, \dots, \mu^{
{N_s}}\} \subset {\mathbb{R}^k}
$, 
where $\mathcal{K}$ is 
a finite dimensional training set of parameters  
inside the parameter space ${\mathcal P}
$,
 $N_s$ is the  number of snapshots,
 $N_h$ denotes the number of degrees of freedom for the discrete full order solution,
 and the snapshots matrix  {${\mathcal{S}}$} is given by $N_s$ full-order snapshots   {properly extended to a fixed background mesh and defined on $\mu^i$ parameter depended domains}:
\begin{equation}
\bm{\mathcal{S}} = [{T}(\mu^1),\dots,{T}(\mu^{N_s})] \in \mathbb{R}^{N_h\times N_s}.
\end{equation}
Given a general scalar function ${T}:{\mathcal D} \to \mathbb{R}^d$, with a certain number of realizations ${T}_1,\dots, {T}_{N_s}$, the POD problem consists in finding, for each value of the dimension of POD space $N_{POD} = 1,\dots,N_s$, the scalar coefficients $a_1^1,$ $\dots,$ $a_1^{N_s},$ $\dots,$ $a_{N_s}^1,$ $\dots,$ $a_{N_s}^{N_s}$ and functions ${\varphi}_1,\dots,{\varphi}_{N_s}$ that minimize the quantity:
\begin{eqnarray*}\label{eq:pod_energy}
E_{N_{POD}} = \sum_{i=1}^{N_s}||{{T}}_i-\sum_{k=1}^{N_{POD}}{ {a_i^k {{\varphi}}_k}}||^2_{L^2({\mathcal D })}, \,
\forall
N_{POD} = 1,\dots,N
\end{eqnarray*}
{with} 
$
({{\varphi}_i,{\varphi}_j}) _{L^2({\mathcal D })} = \delta_{ij}, 
 i,j = 1,\dots,N_s
 $.
This is equivalent of solving the following eigenvalue problem:
$
{\bm{C}}\bm{Q} = \bm{Q}\bm{\lambda}$, for ${C}_{ij} = ({{T}_i,{T}_j}) _{L^2{(\mathcal D )}} \mbox{,\, } i,j = 1,\dots,N_s 
$ 
 ${\bm{C}}$ is the correlation matrix obtained starting from the snapshots $\bm{\mathcal{S}}$, 
 $\bm{Q}$ is a square matrix of eigenvectors,  
 $\bm{\lambda}$ is a diagonal matrix of eigenvalues. 
The basis functions can then be obtained with the formula: 
$
{\varphi_i} = \frac{1}{N_s\lambda^{1/2}_{ii}}\sum_{j=1}^{N_s} {T}_j Q_{ij}.
$
and the POD space:
${\bf{L}} = [{{\varphi}}_1, \dots , {{\varphi}}_{N^r}] \in \mathbb{R}^{N_h \times N^r},$
for $N^r < N_s$ chosen according to the eigenvalue decay of $\bm{\lambda}$, \cite{rozza2008reduced}. 
\newline
\subsection{The projection stage and the generation of the ROM}
The reduced solution, then, can be approximated with: 
\begin{equation}\label{eq:aprox_fields}
{T^r} \approx \sum_{i=1}^{N^r} a_i(\mu) {{\varphi}_i(\bm{x})} = {\bf{L}} \bm{{\bf{a}}}(\mu),
\end{equation} 
for which the reduced solution vectors ${\bf{a}} \in \mathbb{R}^{N^r\times1}$ depend only on the parameter values. The basis functions ${\varphi}_i$ depend only on the physical space, and the unknown vector of coefficients ${\bf{a}}$  can be obtained through 
a Galerkin projection of the full order system of the equations onto the POD reduced basis space, resulting to the consequent reduced algebraic system:
\begin{equation}\label{eq:system_linear_proj}
 {\bf{L}}^T   {\bf{A}}(\mu)   {\bf{L}}   {{\bf{a
}}(\mu)} = {\bf{L}}^T {\bf{F}}(\mu) \mbox{,}
\end{equation}
equivalent to the following algebraic reduced system: 
\begin{equation}\label{eq:system_linear_reduced}
 {\bf{A}}^r(\mu) {\bf{a}}(\mu) = {\bf{F}}^r (\mu)\mbox{,}
\end{equation}
where ${\bf{A}}^r(\mu) \in {\mathbb{R}}^{N^r \times N^r}$,  and ${\bf{F}}^r(\mu) \in {\mathbb{R}}^{N^r \times 1}$ are the reduced discretized operators and reduced forcing vector respectively. 
We highlight that the dimension of the latter reduced system, 
is much smaller with respect to the dimension of the full order system of equations and
 much cheaper to solve.
\subsection{Numerical experiments (Heat exchange/SBM)}
\begin{figure} 
\centering
(i) 
\includegraphics[width=0.58\textwidth]{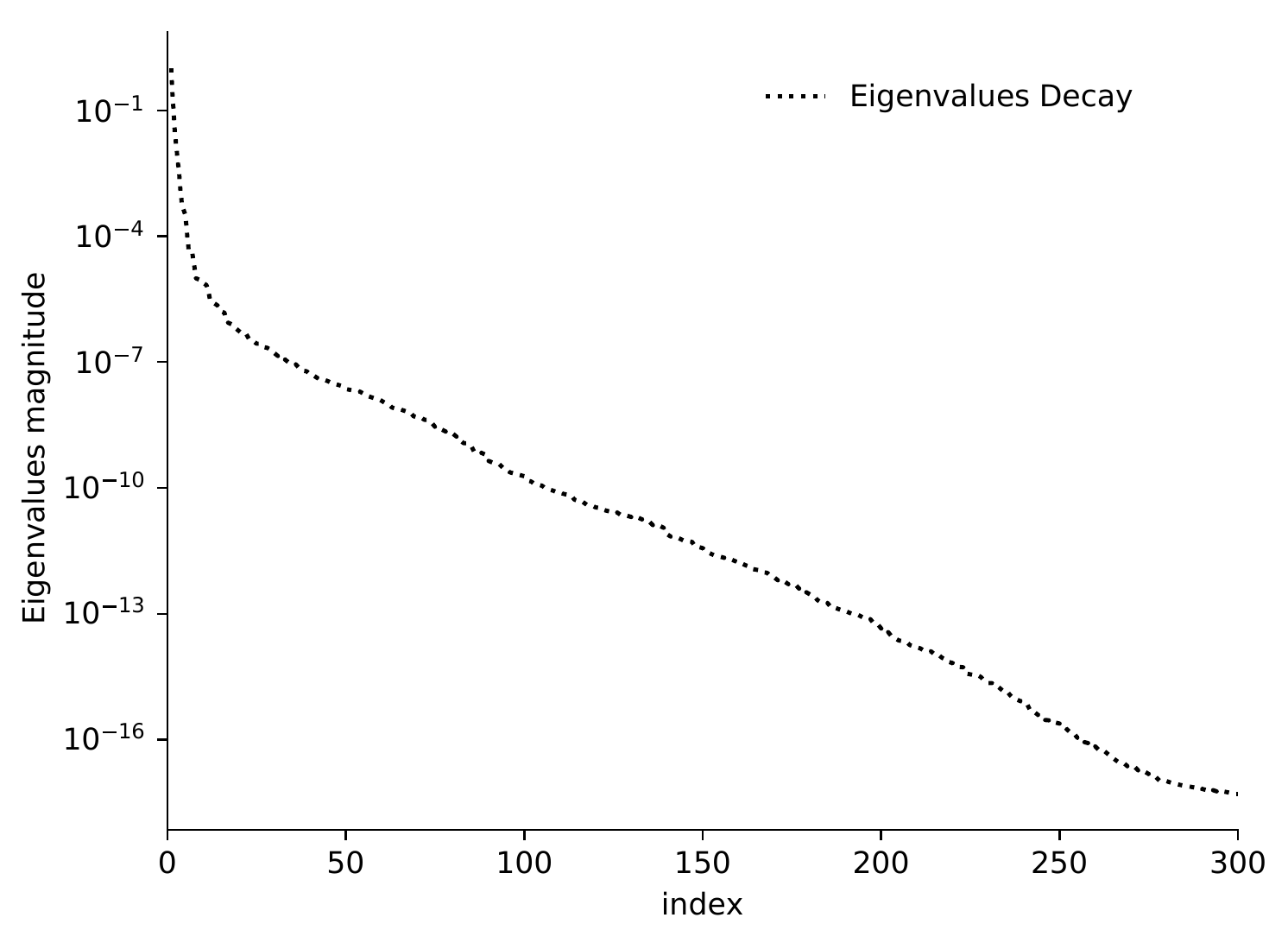}\\
(ii)
\includegraphics[width=0.57\textwidth]{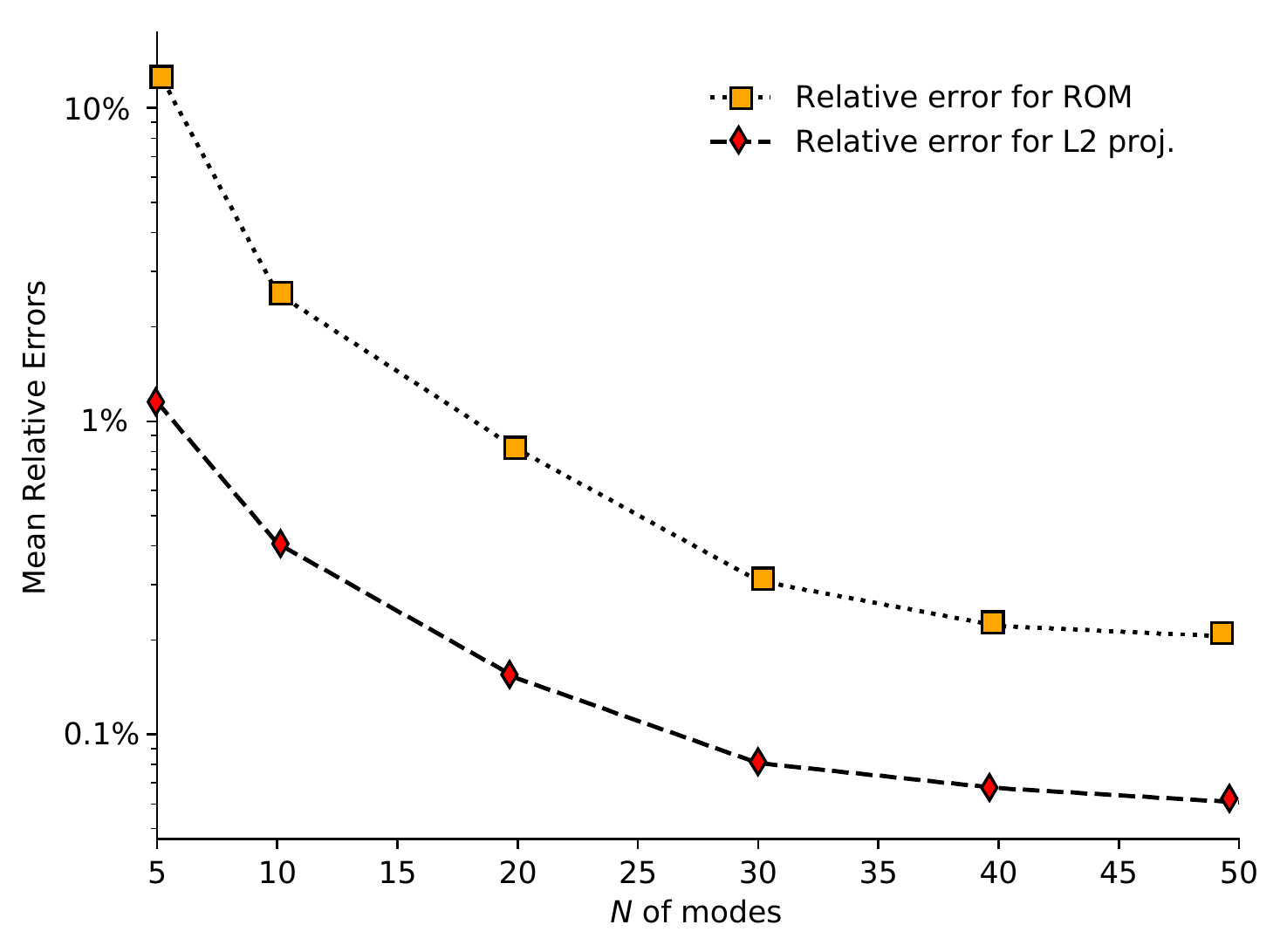}
\caption{Heat exchange problem:  
(i) eigenvalues decay and (ii) mean relative errors.
}
\label{fig:poisson_results}
\end{figure}

We assume that the embedded domain consists of a rectangle of size $0.8\times 0.7$, and  its position inside the domain is parametrized with a geometrical parameter $\mu$. 
 The position $(0,\mu)$ of the rectangle embedded domain depends on its parametrized $y$-center, while the horizontal coordinate of the center of the box is not parametrized and is located in the $x$-center of the domain. 
All this configuration is immersed in a background domain of size  $\mathcal{D}=[-2,2]\times[-1,1]$.
 The ROM has been trained with  $400$ samples {, which is the dimension of the offline  FOM,} for $\mu \in [-0.5,0.5]$ chosen randomly inside the parameter space,
while the ROM results have been compared with FOM for $50$ additional random samples.
\begin{figure}
 (i)
 \\
 \includegraphics[width=0.3\textwidth]{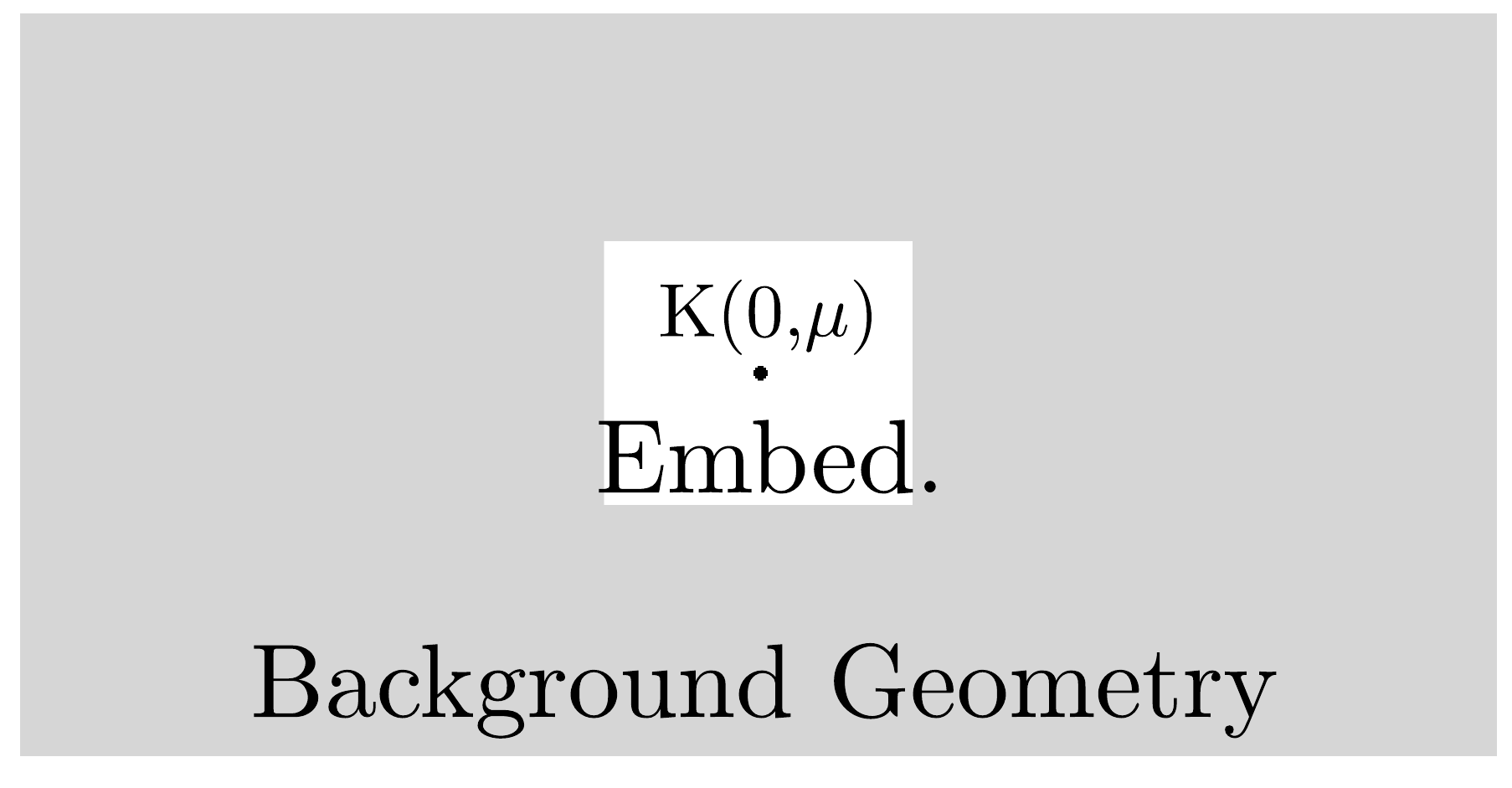}
 \hfil \hfil \hfil  \hfil\hfil \hfil  \hfil\hfil \hfil  \hfil\hfil \hfil  \hfil\hfil \hfil  \hfil\hfil \hfil  \hfil\hfil \hfil  \hfil\hfil \hfil  \hfil\hfil \hfil  \hfil\hfil \hfil  \hfil\hfil \hfil  \hfil\hfil \hfil  \hfil\hfil \hfil  \hfil\hfil \hfil  \hfil\hfil \hfil  \hfil\hfil \hfil  \hfil\hfil \hfil  \hfil\hfil \hfil  \hfil
\\(ii)
\\
\includegraphics[width=0.295\textwidth]{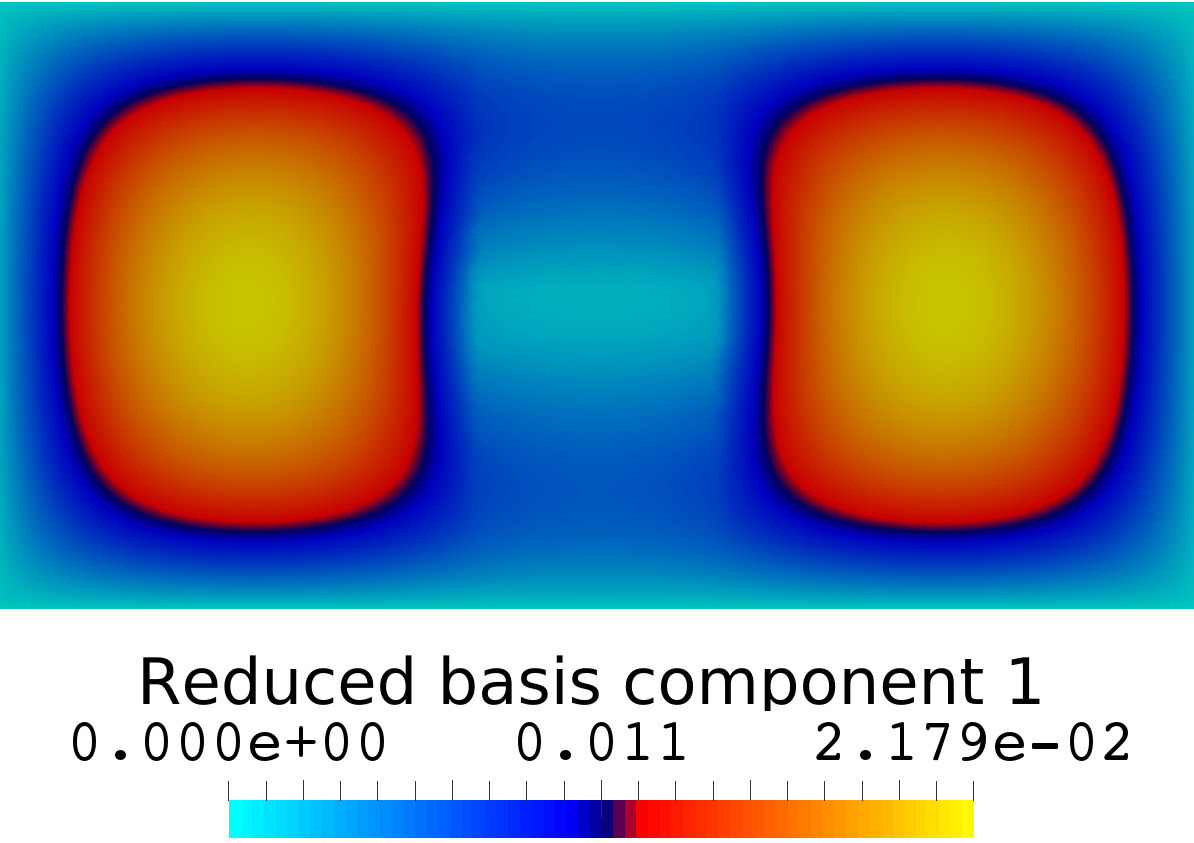} 
\includegraphics[width=0.295\textwidth]{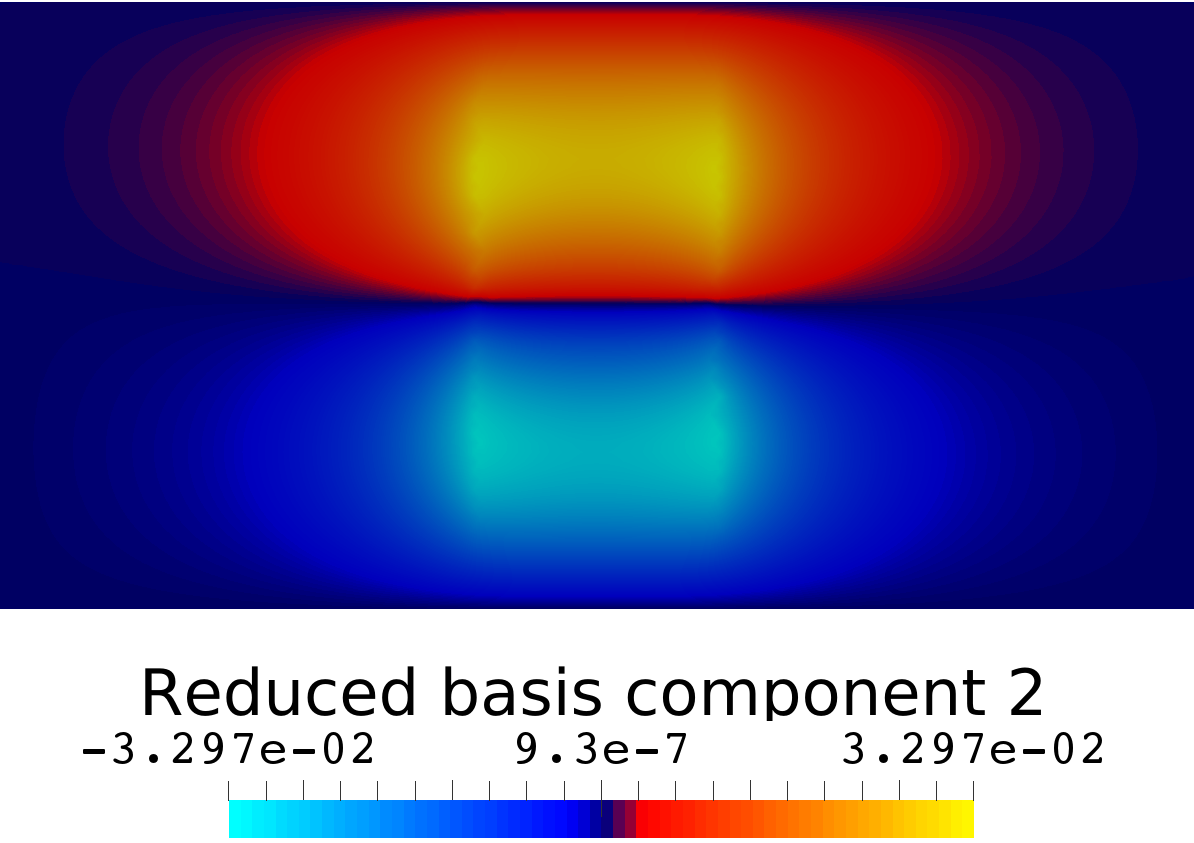}
\\
\includegraphics[width=0.295\textwidth]{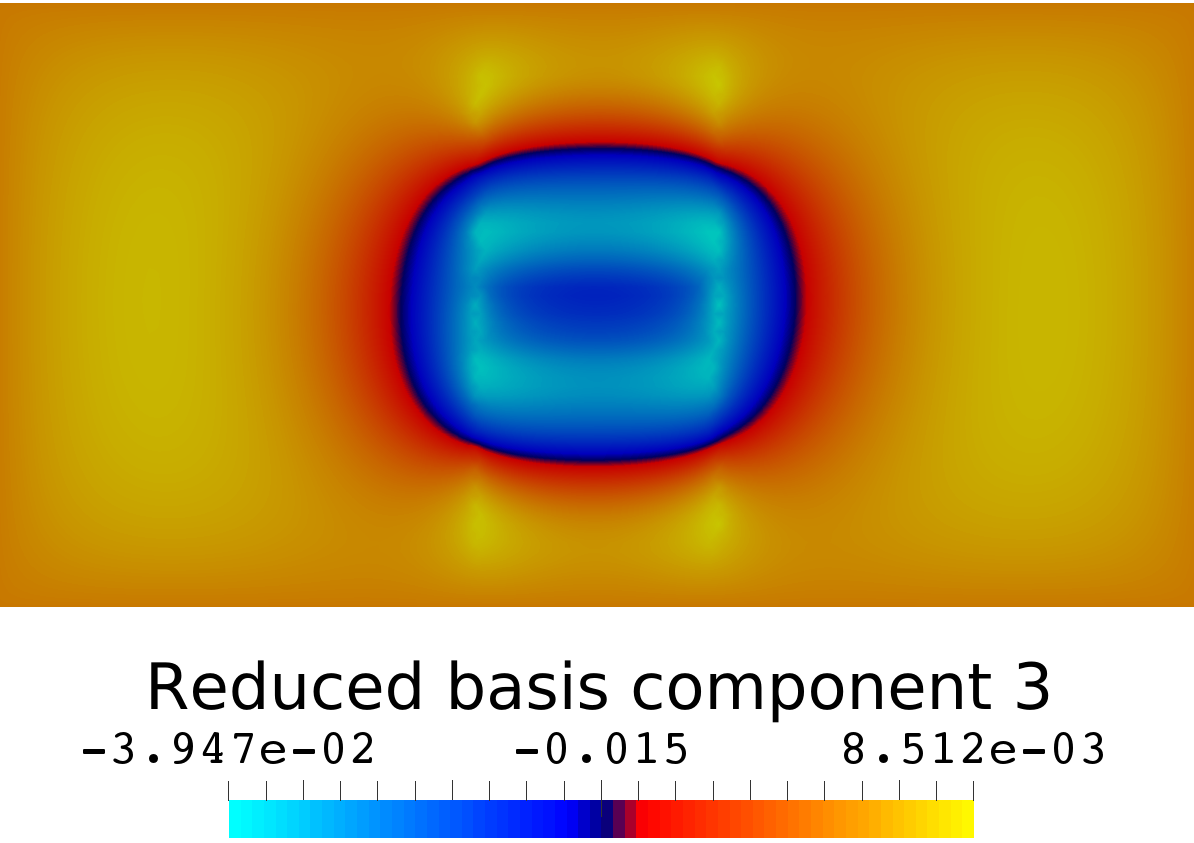}
\includegraphics[width=0.295\textwidth]{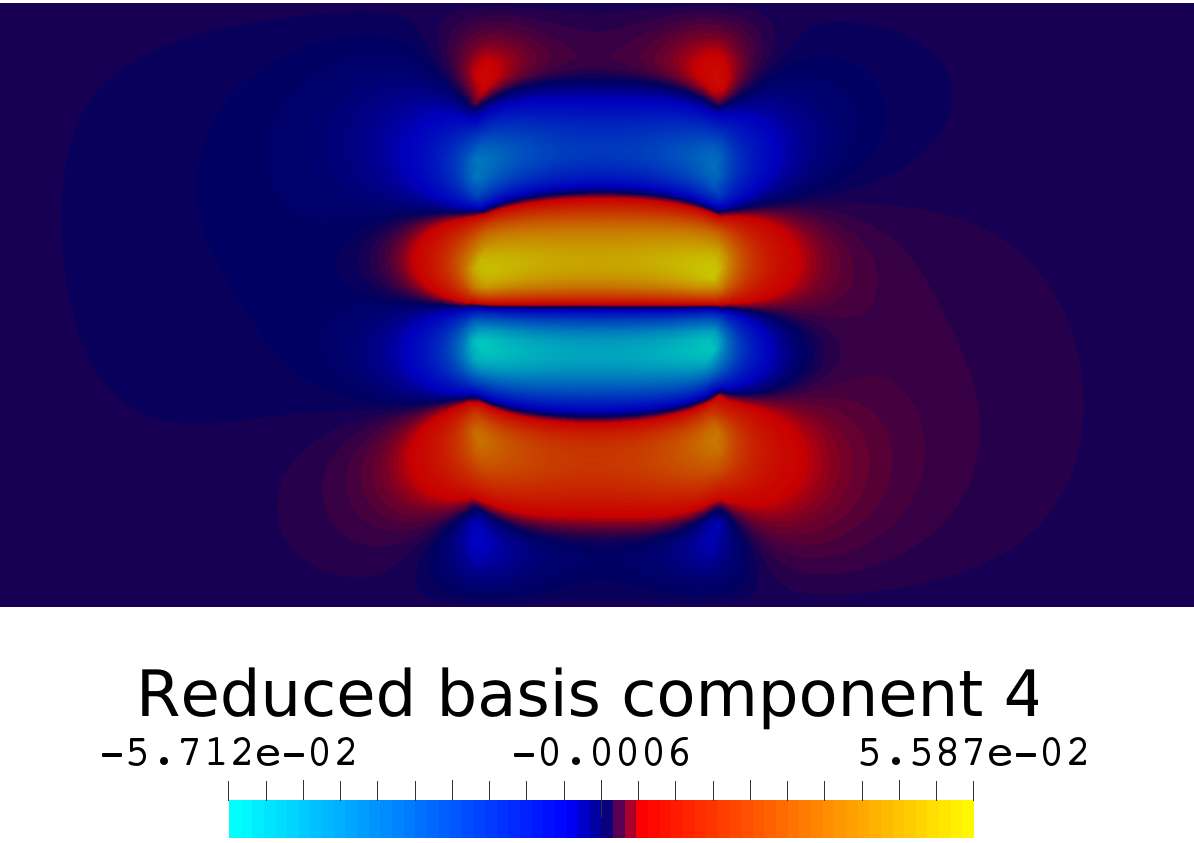}
\\
(iii)
\\
\includegraphics[width=0.3\textwidth]{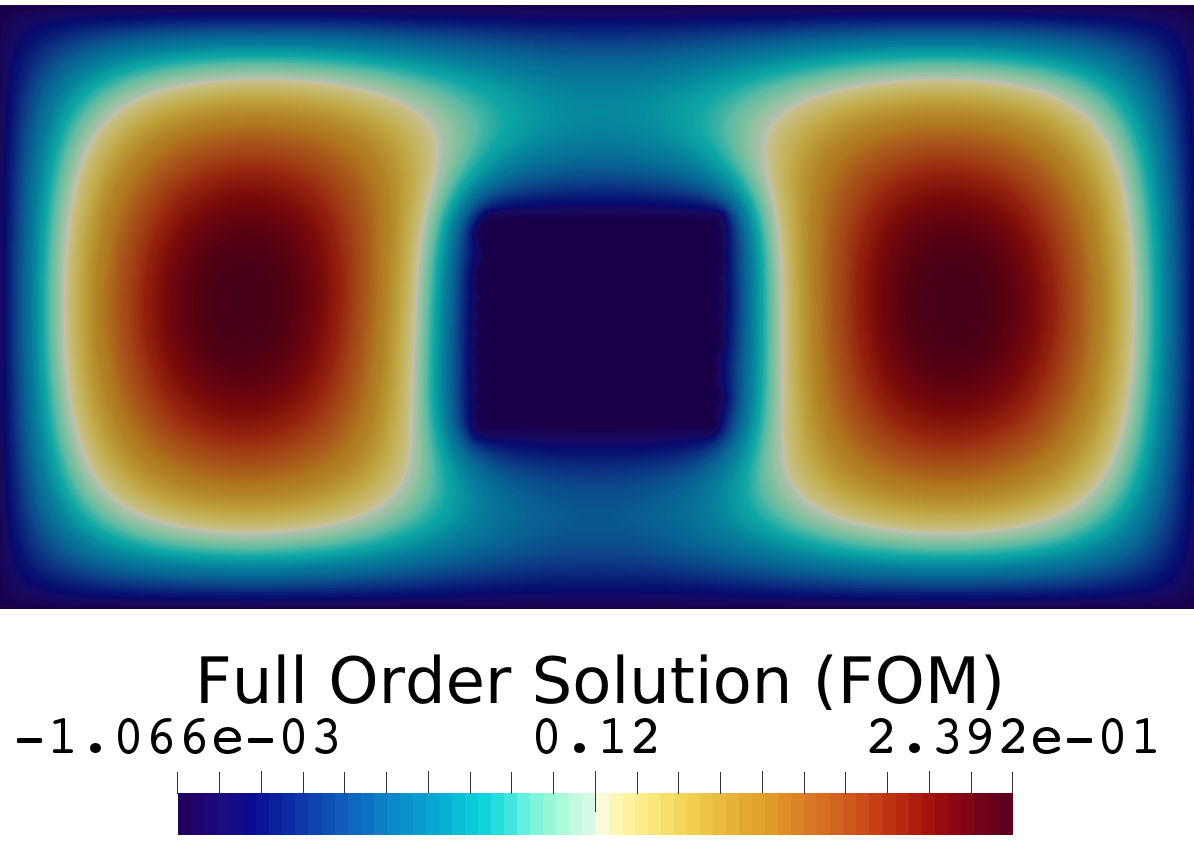}
\includegraphics[width=0.3\textwidth]{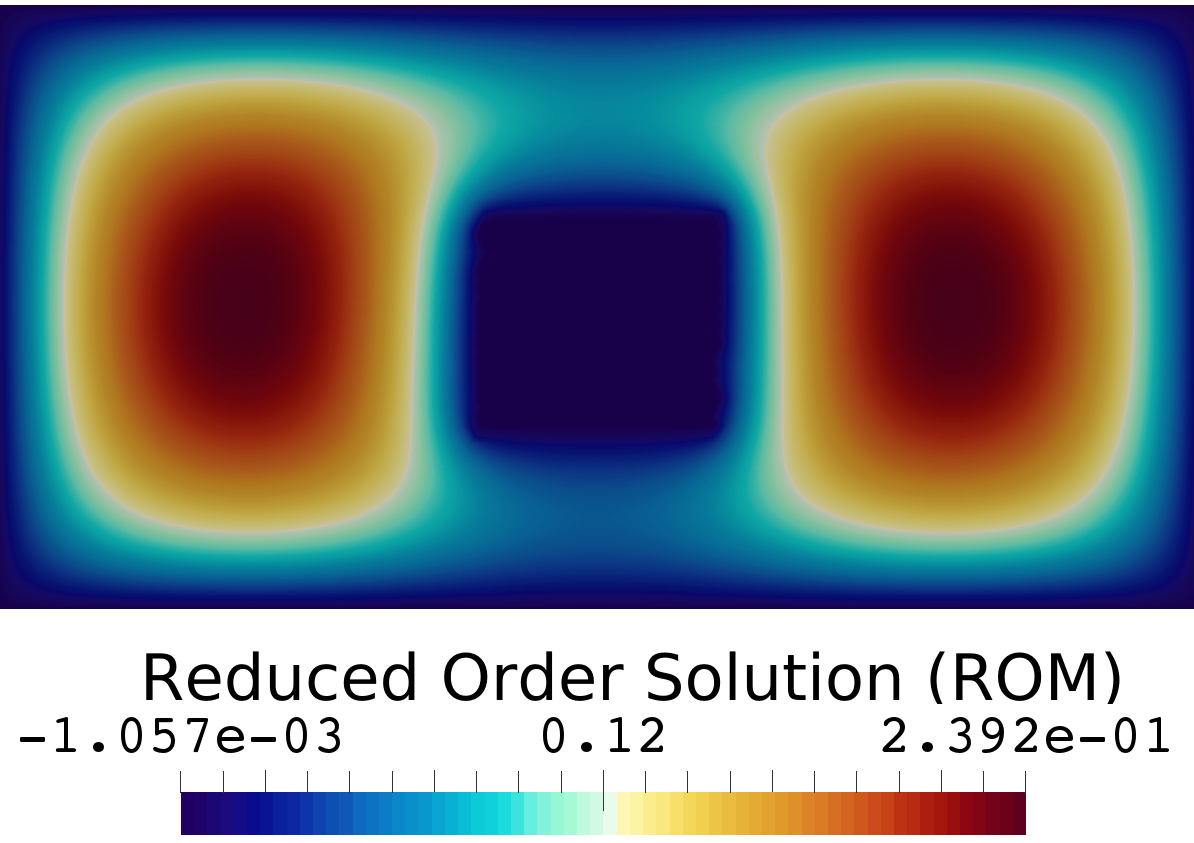}
\includegraphics[width=0.3\textwidth]{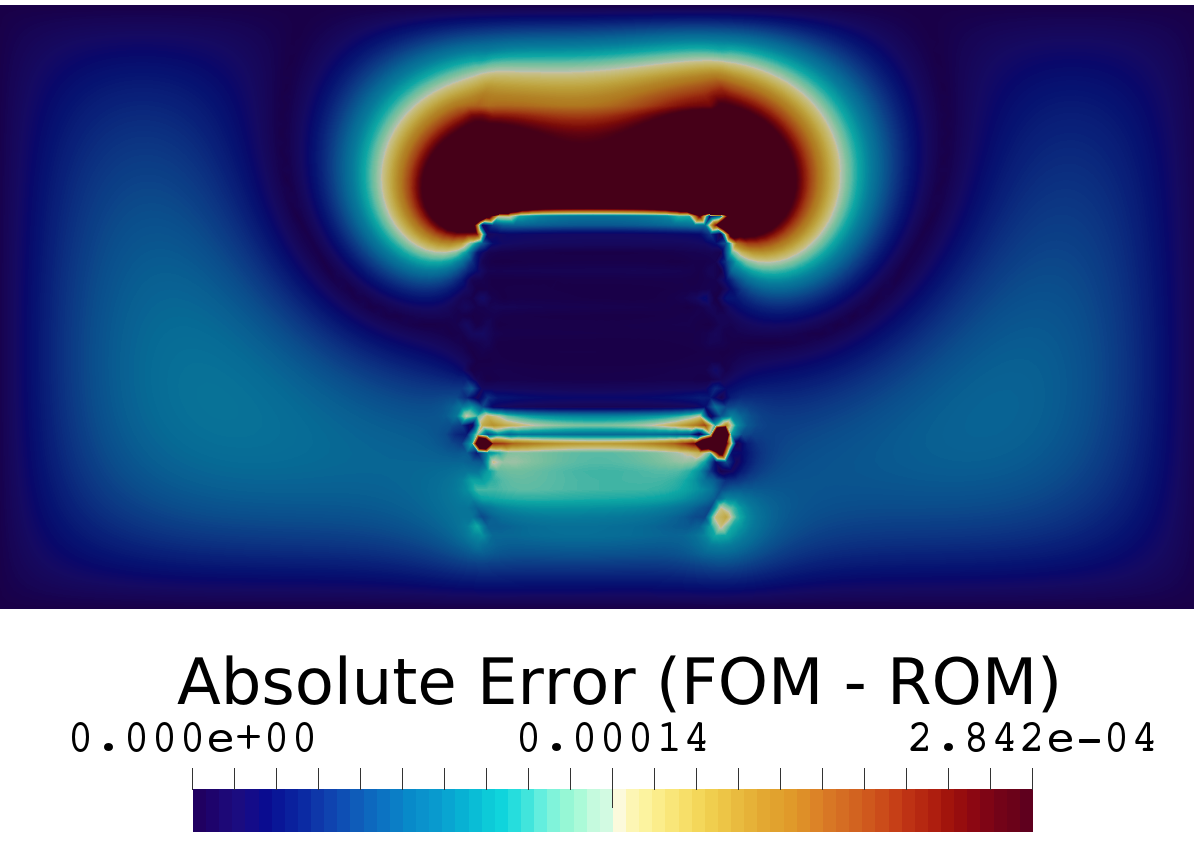}
\caption{(i) Parametrized geometry (ii)  {Some reduced} basis components for  $\mu \in [-0.5, 0.5]$ parametrized geometry. 
(iii) The full/reduced order solution and the absolute error ($\mu = -0.015$).}
\label{fig:poisson_field}
\end{figure}
The mean relative errors, 
the eigenvalue decay and the behavior of the results for one parameter sample 
are visualized in Figures \ref{fig:poisson_field} and \ref{fig:poisson_results} and Table \ref{table:poisson_rel_er}. The time savings are reported in Table \ref{tab:Heat_Times}  {where}
\begin{table}
\centering
\caption{ {Heat exchange problem: relative error results
.}}
\begin{tabular}{|c|c|c|} 
\hline
  Modes         & $L^2$ projection 
  &  \multicolumn{1}{c|}{Galerkin projection
  } 
  \\ \hline 
2   & 6.45035392e-02 &  7.10916700e-01  
\\\hline
10  & 4.83332393e-03 &  2.64459969e-02    
\\ \hline
20  & 2.19454585e-03 &  5.61736266e-03    
\\ \hline
30  & 1.27046941e-03 &  3.30372025e-03    
\\ \hline
40  & 7.72326410e-04 &  2.50189079e-03   
\\ \hline
50  & 5.39532759e-04 &  1.69903034e-03    
\\ \hline
100 & 6.79464703e-05 &  3.36531580e-04    
\\  
\hline
\end{tabular}
\label{table:poisson_rel_er}
\end{table}
\begin{table}
\centering
\caption{ {Heat exchange problem: execution time, savings and speed up
.
}}
\begin{tabular}{|c|c|c|c|}
\hline
\multirow{2}{*}{Modes} & \multirow{2}{*}{Execution time(s)
} & Savings & Speedup\\ 
 &    & $\frac{t_{{\text{\tiny FOM}}}-t_{{\text{\tiny RB}}}}{{t_{{\text{\tiny RB}}}}}
 $  &  $\frac{{t_{\text{\tiny FOM}}}}{{t_{\text{\tiny RB}}}}$ \\ \hline
2 &  $4.119470\times 10^{-2}$ & 96.399\% & 27.770 \\\hline
10 & $4.168334\times 10^{-2}$ & 96.356\% & 27.445   \\\hline
20 & $4.243647\times 10^{-2}$ & 96.290\% & 26.957   \\\hline
30 & $4.353909\times 10^{-2}$ & 96.194\% & 26.275   \\\hline
40 & $4.449359\times 10^{-2}$ & 96.110\% & 25.711   \\\hline
50 & $4.494564\times 10^{-2}$ & 96.071\% & 25.452   \\\hline
100& $4.992923\times 10^{-2}$ & 95.635\% & 22.912 \\
\hline 
\text{FOM}& $1.14540\times 10^{0}$ & -- & -- \\ \hline
\end{tabular}
\label{tab:Heat_Times}
\end{table}
the computation time includes the assembling of the full order matrices, their projection and the resolution of reduced problem, 
while we have avoided a) {the usage of reference domain}, b) {remeshing}, c) {adaptive refinement}.
\section{Parametrized Steady Stokes equations}\index{Stokes equations, Parametrized, Steady } 
In this section, we examine a Newtonian, incompressible viscous fluid flow, in a domain when 
the convective forces are negligible with respect to viscous forces, 
namely the Stokes system
\begin{eqnarray*}
-\nabla\cdot (2\nu{{\bm{\epsilon}} (\bm u(\mu))} - p(\mu){ \bm I}) &=& \bm g(\mu),\, \,\,\,\text{ in } {\mathcal{D}(\mu)}, 
\\
\nabla \cdot { \bm u(\mu)} &=& 0, \qquad \,\text{ in }{\mathcal{D}(\mu)}, 
\\ \bm u(\mu) &=& \bm g_D(\mu), \text{ on } \Gamma
(\mu),
\\
(2\nu { \bm \epsilon ( \bm u(\mu))} - p(\mu){ \bm I}) \cdot { \bm n} &=& { \bm g_N(\mu)} , \text{ on } \Gamma _N(\mu),
\end{eqnarray*}
where 
 ${ \bm  \epsilon ( \bm u)} = 1/2( \nabla   {  \bm u}+ \nabla {  \bm u}^T )$ is the velocity strain tensor (i.e., the symmetric gradient of the velocity), 
 $p$ is the pressure, $ \bm g$ a body force, 
 $ \bm  g_{D}$, the values of the velocity on the Dirichlet and 
 $  \bm g_{N}$ is the normal stress on the Neumman boundary. 

Based on \cite{KaratzasStabileNouveauScovazziRozza2018}, next we define the Shifted Boundary weak formulation 
{\it Find ${\bm u} \in  {\bm V}_h $ and $p \in Q_h $ such that, $\forall {\bm w} \in {\bm V}_h $ and $\forall q \in Q_h,$}
\begin{eqnarray}
\nonumber\label{eq:SBM}
( {\bm \epsilon ( \bm w)}, 2 {\nu} {\bm \epsilon (\bm u)})_{\tilde {\mathcal{D}}(\mu)} 
- \left\langle {\bm w} \otimes {\bm {\tilde n}}, 2 {\nu} {\bm \epsilon (\bm u)}\right \rangle _{\tilde\Gamma 
(\mu)}
- \left\langle 2 {\nu} {\bm \epsilon (\bm w)} , ({\bm u} + (\nabla {\bm u})\cdot{\bm d} ) \otimes {\bm{ \tilde n}}\right\rangle  _{\tilde\Gamma 
(\mu)}
\nonumber\\
+ \alpha\left \langle 2 {\nu}/h({\bm w} + ( \nabla{\bm w})\cdot{\bm d}), {\bm u} + (\nabla  \bm u)\cdot{\bm d} \right \rangle  _{\tilde\Gamma 
(\mu)}
+ \beta \left\langle 2 {\nu} h \nabla _{\bar \tau _i} {\bm w}, \nabla _{\bar \tau _i} {\bm u}  \right\rangle_{\tilde\Gamma 
(\mu)} 
 - (\nabla \cdot {\bm w}, p)_{\tilde {\mathcal{D}}(\mu)}
 \nonumber
\\
+\left\langle  {\bm w} \cdot  {\bm{\tilde n}}, p\right \rangle _{\tilde\Gamma 
(\mu)}=-({\bm w}, {\bm g})_{\tilde {\mathcal{D}}(\mu)}
+ \left\langle {\bm w}, {\bm g_N }  \right\rangle_{\tilde\Gamma _N(\mu)}
- \left\langle 2 {\nu} {\bm \epsilon (\bm w)},   {\bm {\bar g_D}} \otimes {\bm {\tilde  n}}\right\rangle  _{\tilde\Gamma 
(\mu)}
\nonumber
\\
+ \alpha\left \langle 2 {\nu}/h({\bm w} + ( \nabla{\bm w})\cdot{\bm d}),  {\bm {\bar g_D}}\right \rangle  _{\tilde\Gamma 
(\mu)}+ \beta \left\langle 2 {\nu} h \nabla _{\bar \tau _i} {\bm w},  \nabla _{\bar \tau _i} {\bm{ \bar 
g_D}} \right\rangle_{\tilde\Gamma 
(\mu)},
\nonumber\\
\text{and}\nonumber
\\
- (\nabla \cdot {\bm u},q)_{\tilde {\mathcal{D}}(\mu)} 
+\left\langle  {\bm u} \cdot  {\bm{\tilde n}}, q\right \rangle _{\tilde\Gamma 
(\mu)}
+\left\langle q {\bm d} \otimes {\bm{\tilde n}}, \nabla {\bm u}\right \rangle _{\tilde\Gamma 
(\mu)}
=\left\langle {\bm{\bar  g_D}}\cdot  \bm {\tilde{n}} ,q
\right \rangle _{\tilde\Gamma 
(\mu)},
\nonumber
\end{eqnarray}
or in a more abstract notation
\begin{eqnarray*}
a({\bm u},{\bm w};\mu)+b(p,{\bm w};\mu)=\ell _g(\bm w;\mu),
\\
b(q,{\bm u};\mu)+\hat b(q,{\bm u};\mu)=\ell _q(q;\mu),
\end{eqnarray*} 
 {where  $\bm n$ and  ${\bm \tau} _i$ are  the unit normal vector and unit tangential vectors to the boundary $\Gamma$  and can be extended to
the boundary $\tilde\Gamma$, namely: ${\bm{\tilde n}}(\bm{\tilde x}) \equiv {\bm{n}}(\bf M(\bm{\tilde x}))$ and $\bm{\bar \tau}_i(\bm{\tilde x}) \equiv \bm{\tau}_i({\bf M}((\bm{\tilde x}))$, 
${\bm{\bar{g}_ D}}(\bm{\tilde x}) = {\bm {g_D}}({\bf M}(\bm{\tilde x})) $.
We clarify that $\mu$ is defined similarly to subsection \ref{POD} and $h$ is a characteristic length of the elements.}
\subsection{POD adapted to flows}\label{subsec_POD_theory}\index{Proper Orthogonal Decomposition!adapted to flows}
Again, we collect
 $N_s $, number of snapshots, 
with $N_u^h$,  $N_p^h$, number of dof 
for the discrete full order solution  for the velocity and pressure, and we construct seperate basis for both velocity and pressure

\begin{equation}
\nonumber
 {L_u = [{\bm{\varphi}}_1, \dots , {\bm{\varphi}}_{N_u^r}] \in \mathbb{R}^{N_{u}^h \times N_u^r},}\quad
L_p = [{\chi_1}, \dots , {\chi_{N_p^r}}] \in \mathbb{R}^{N_{p}^h \times N_p^r}.
\end{equation}
where $N_u^r$, $N_p^r < N_s$ are chosen according to the eigenvalue decay of the vectors of eigenvalues ${\lambda}^u$ and ${\lambda}^p$. 
Furthermore for best approximation we employ the supremizer enrichment as in the work of \cite{RoVe07} 
and we manage a solvable and stable problem that
satisfies a reduced and also parametric, version of the  ``inf-sup'' condition. 
Within this approach, the velocity supremizer basis functions $L_{\text{sup}}$, with $L_{\text{sup}} =[{ {{
\bm
{\eta}_1}, \dots, {
\bm
{\eta}}_{N_{\text{sup}}^r}}}] \in \mathbb{R}^{N_u^h \times N_{\text{sup}}^r}, 
$ 
are computed and added to the 
	reduced velocity space which is transformed into 
	$\tilde{L}_u$:

\begin{equation}
  \nonumber
        {\tilde{L}_u = [{\bm{\varphi}_1,\dots,\bm{\varphi}_{N_u^r}},{ {
\bm
{\eta}_1,\dots,
\bm
{\eta}_{N_{\text{sup}}^r}}}] \in \mathbb{R}^{N_u^h \times (N_u^r+N_{\text{sup}}^r)} }.
\end{equation}
\subsection{Steady Stokes numerical experiments (SBM)} 
\begin{figure}
\centering (i)
\\
\centering
\begin{minipage}{\textwidth}
\centering
\begin{minipage}{0.24\textwidth}
  \includegraphics[width=\textwidth]{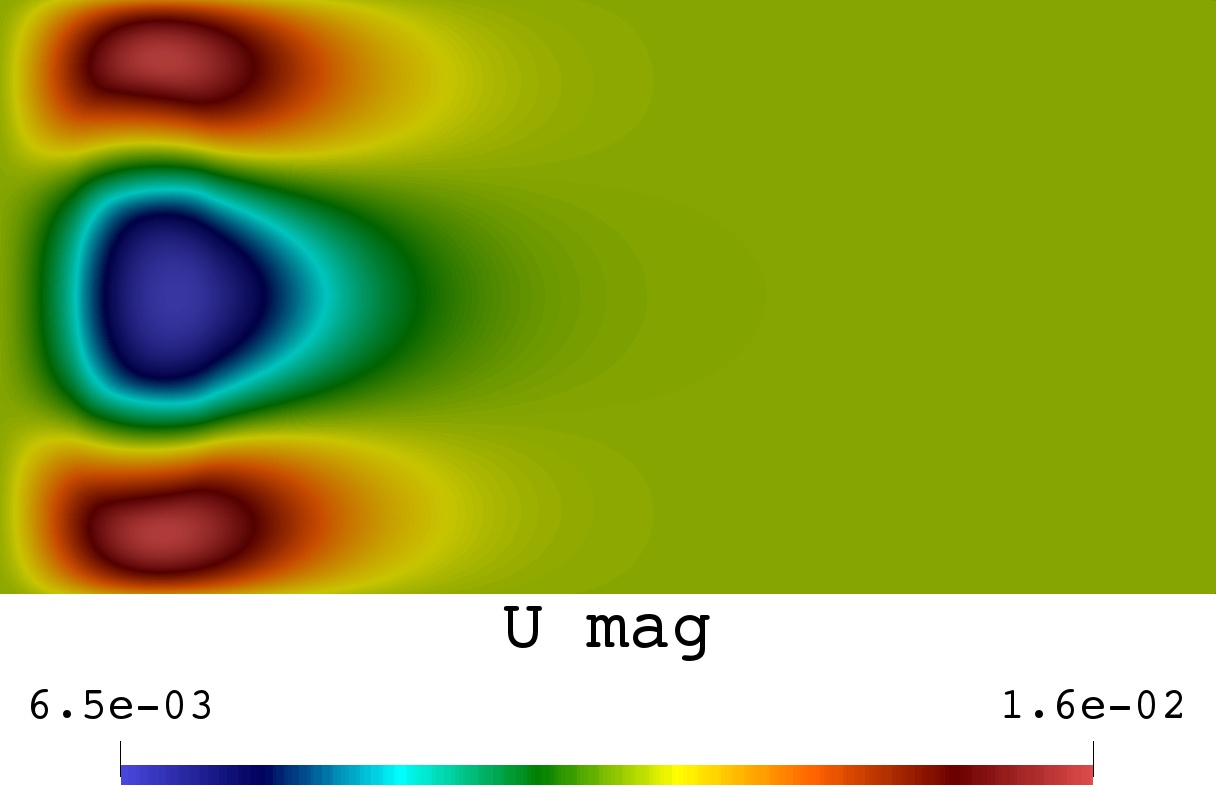}
\end{minipage}
\begin{minipage}{0.24\textwidth}
  \includegraphics[width=\textwidth]{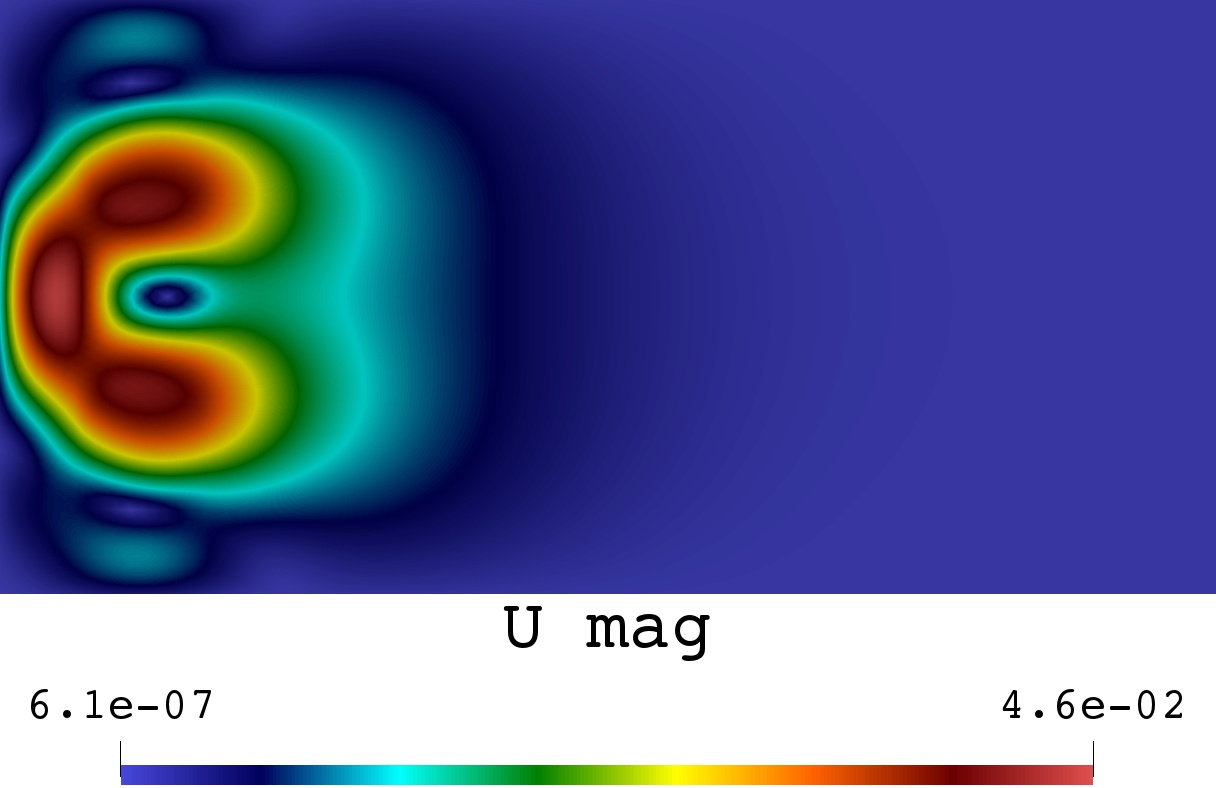}
\end{minipage}
\begin{minipage}{0.24\textwidth}
  \includegraphics[width=\textwidth]{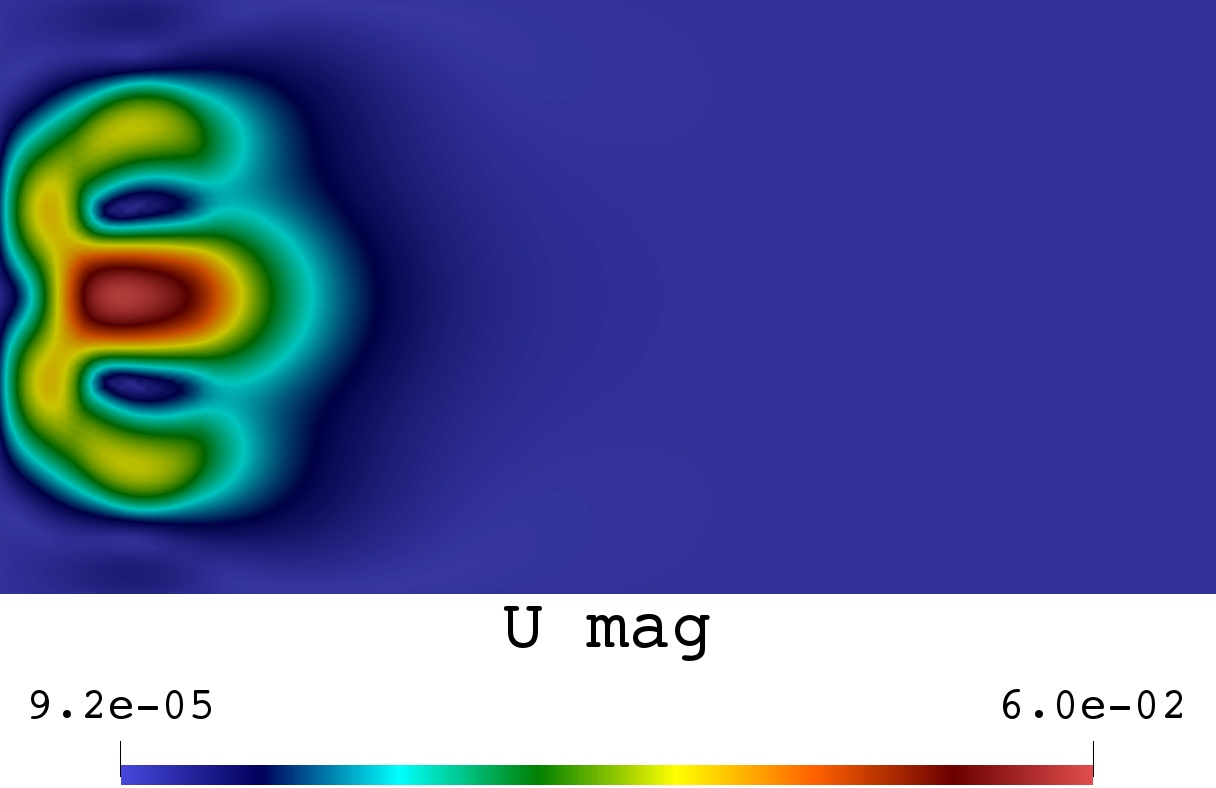}
\end{minipage}
\begin{minipage}{0.24\textwidth}
  \includegraphics[width=\textwidth]{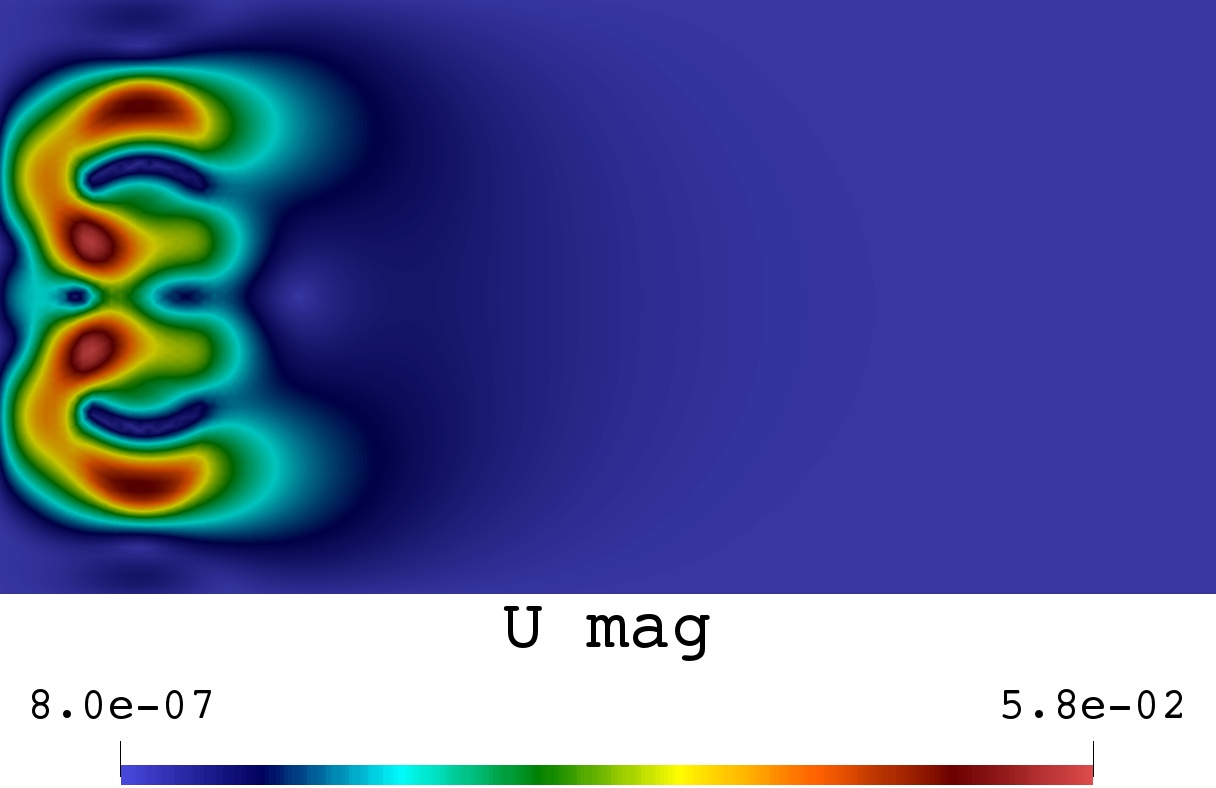}
\end{minipage}
\\
\centering (ii)
\\
\begin{minipage}{0.24\textwidth} 
  \includegraphics[width=\textwidth]{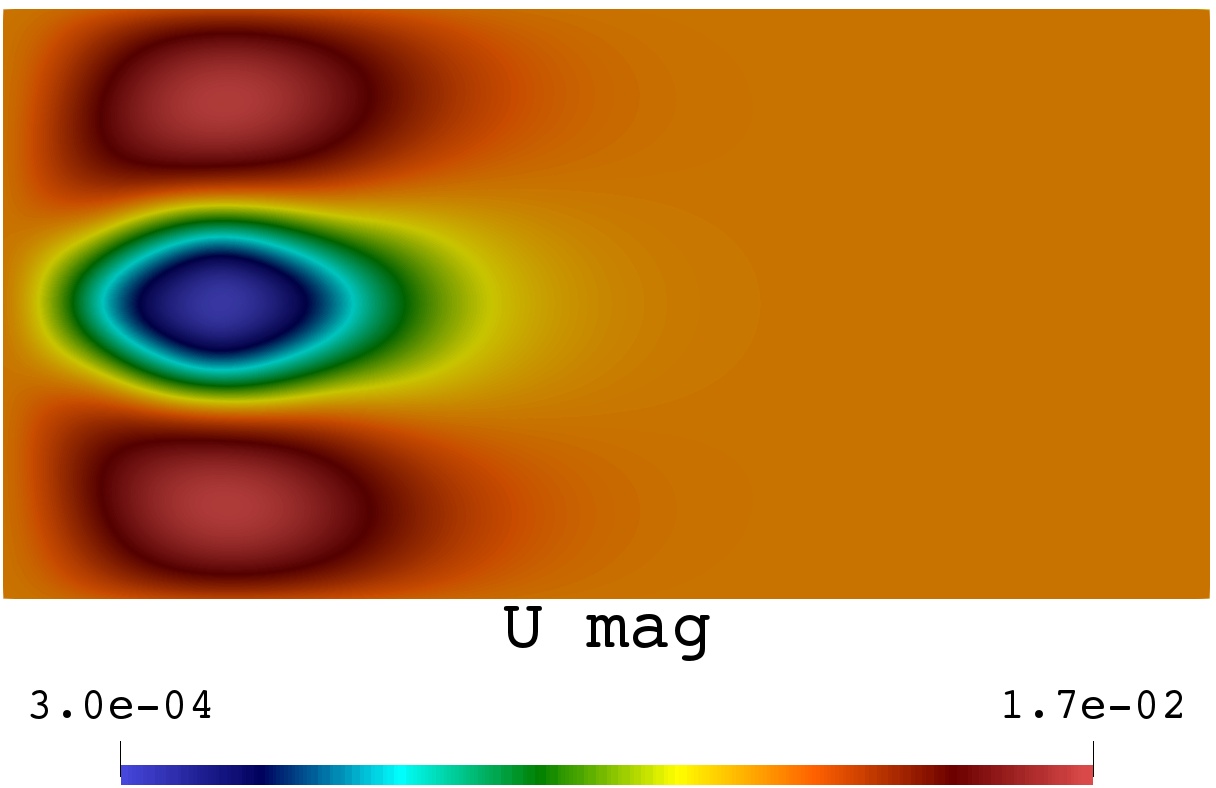}
\end{minipage}
\begin{minipage}{0.24\textwidth}
  \includegraphics[width=\textwidth]{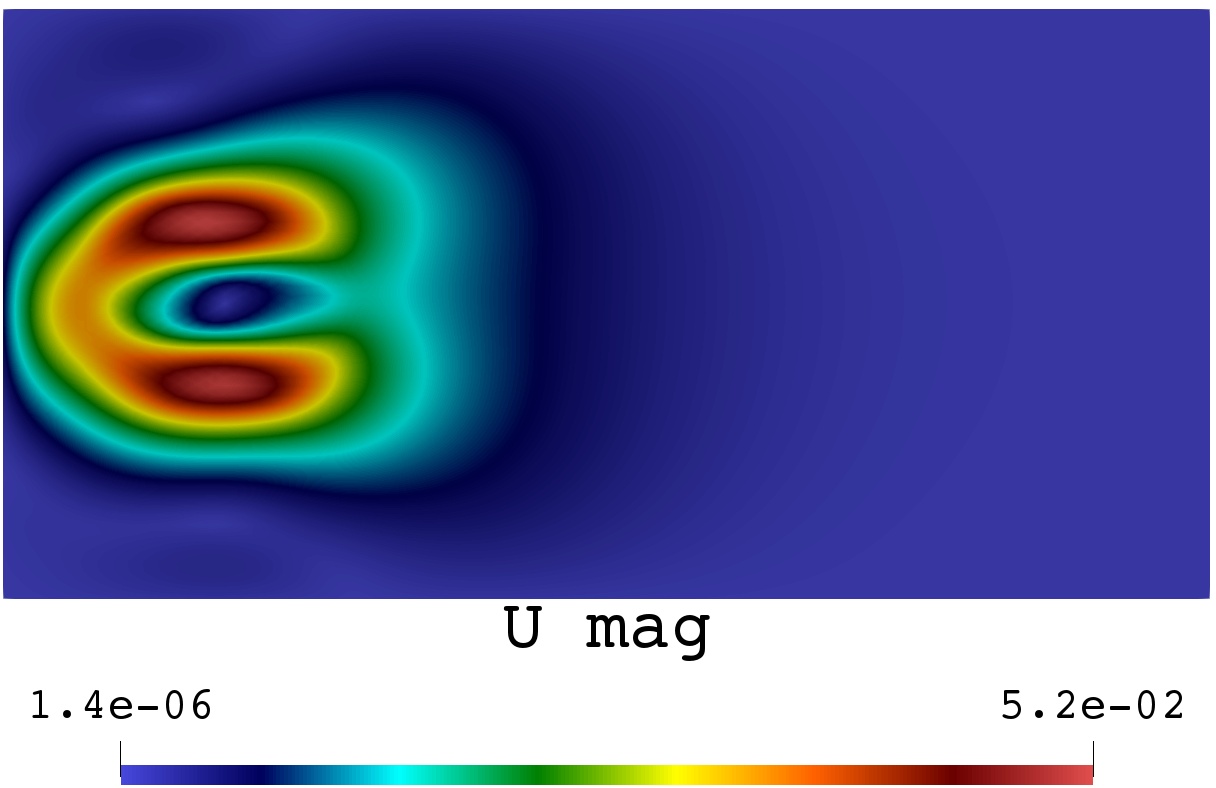}
\end{minipage}
\begin{minipage}{0.24\textwidth}
  \includegraphics[width=\textwidth]{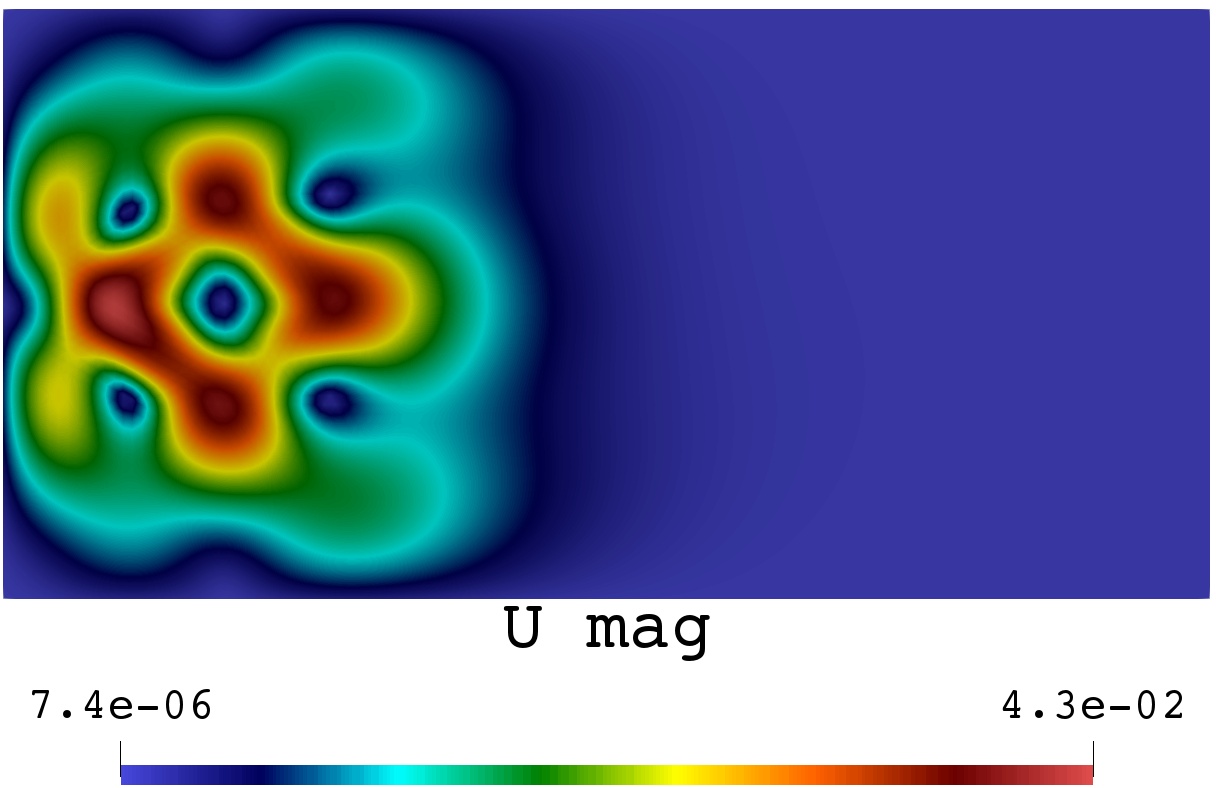}
\end{minipage}
\begin{minipage}{0.24\textwidth}
  \includegraphics[width=\textwidth]{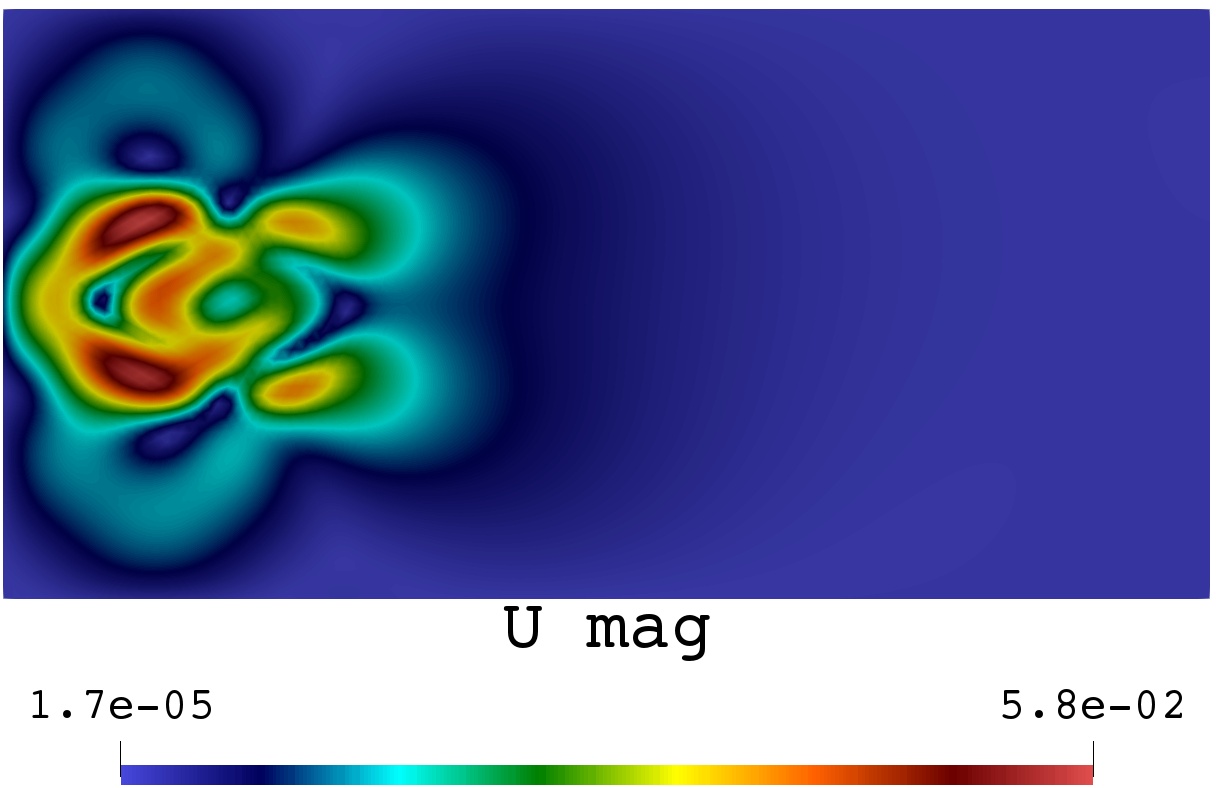}
\end{minipage}
\\
(i)
\\
\begin{minipage}{0.24\textwidth}
  \includegraphics[width=\textwidth]{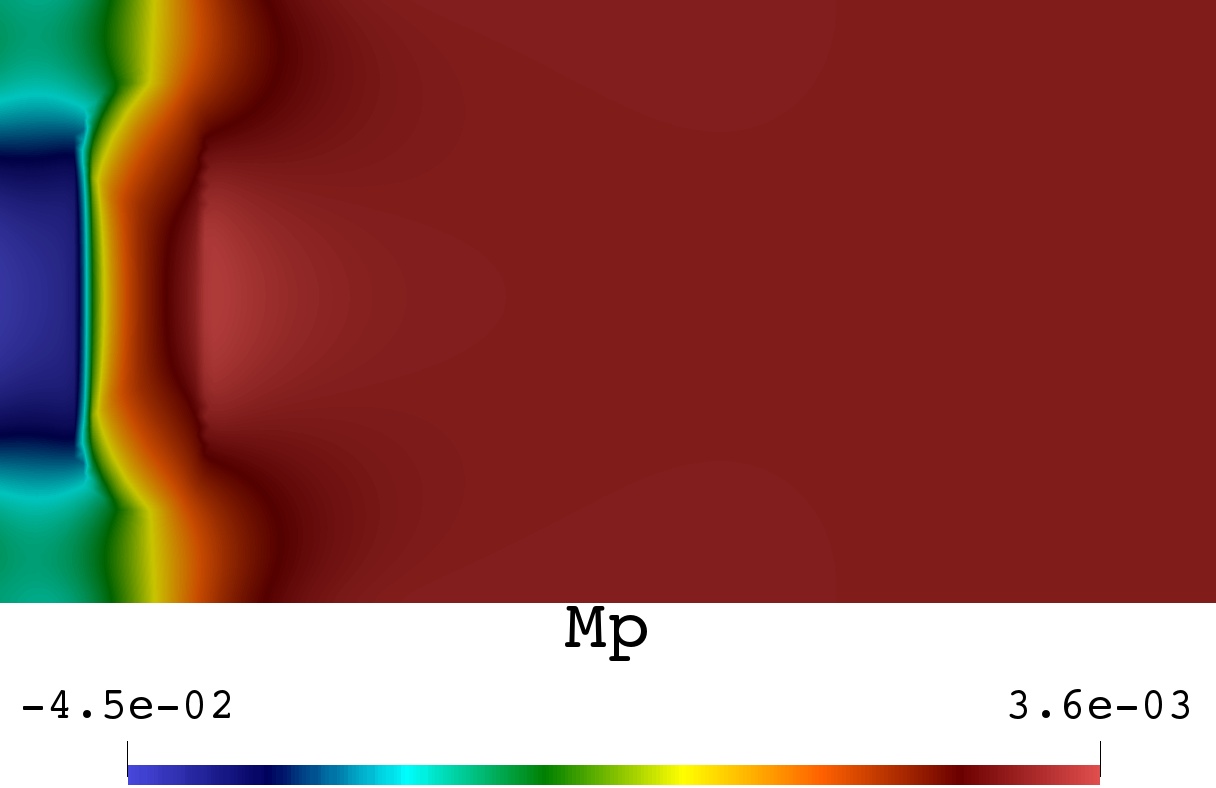}
\end{minipage}
\begin{minipage}{0.24\textwidth}
  \includegraphics[width=\textwidth]{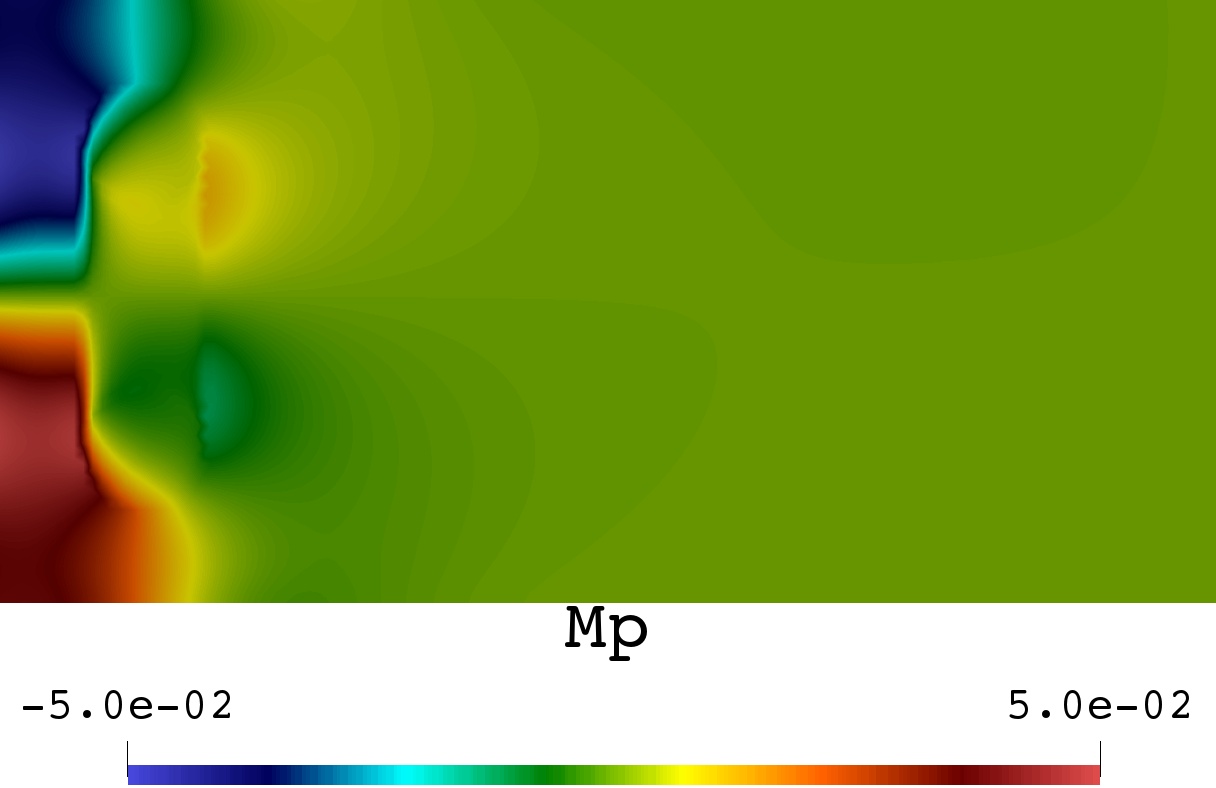}
\end{minipage}
\begin{minipage}{0.24\textwidth}
  \includegraphics[width=\textwidth]{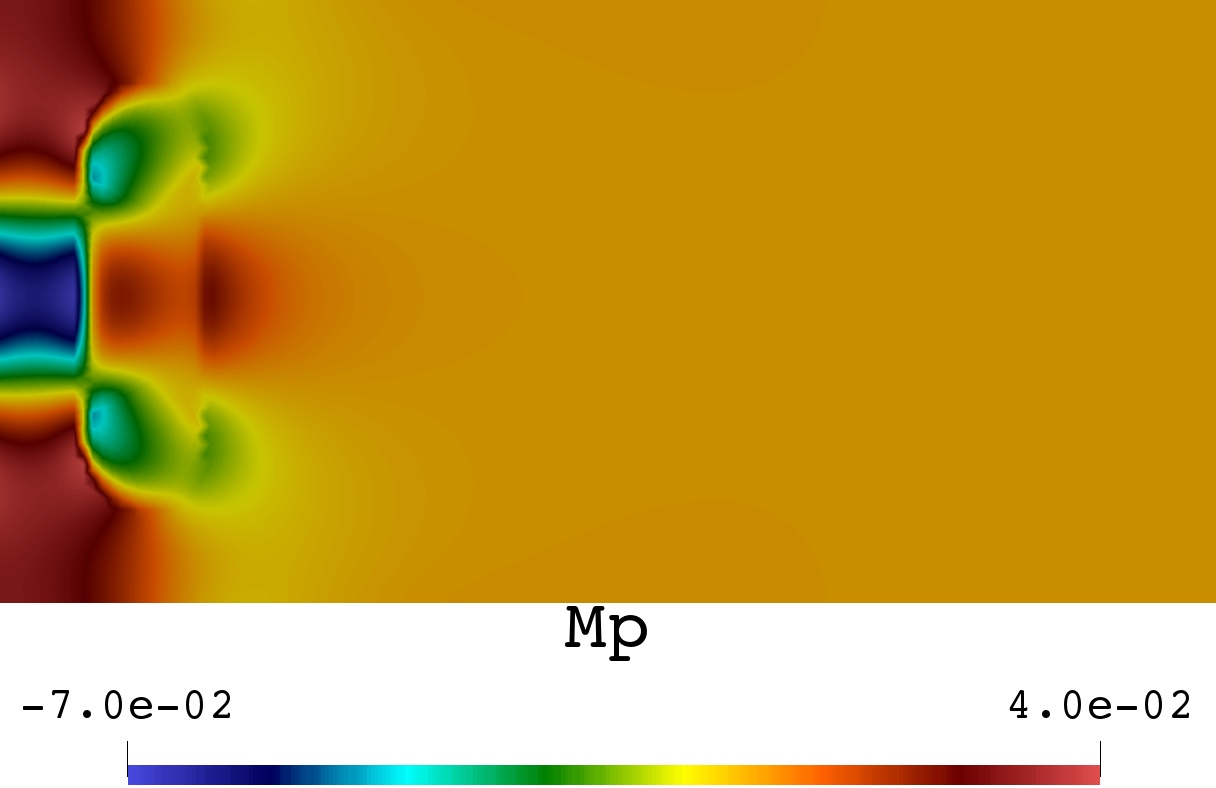}
\end{minipage}
\begin{minipage}{0.24\textwidth}
  \includegraphics[width=\textwidth]{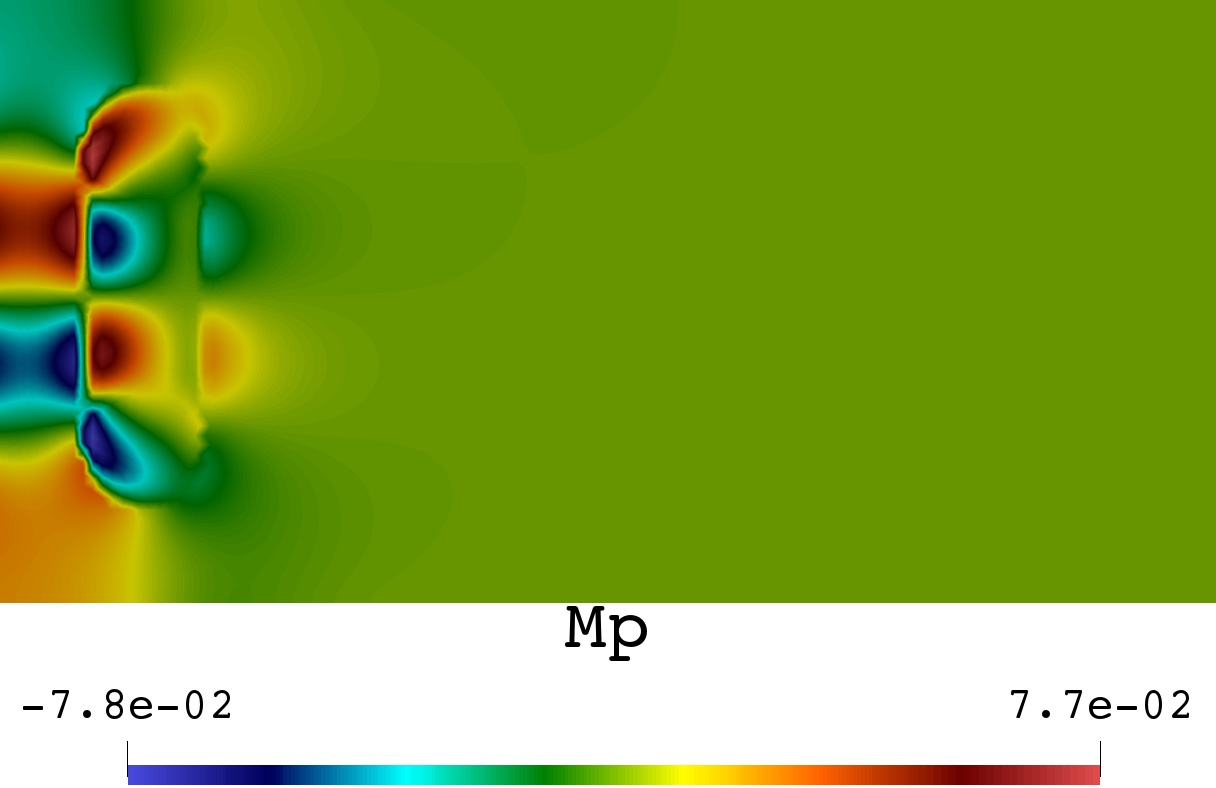}
\end{minipage}
\\
(ii)
\\
\begin{minipage}{0.24\textwidth}
  \includegraphics[width=\textwidth]{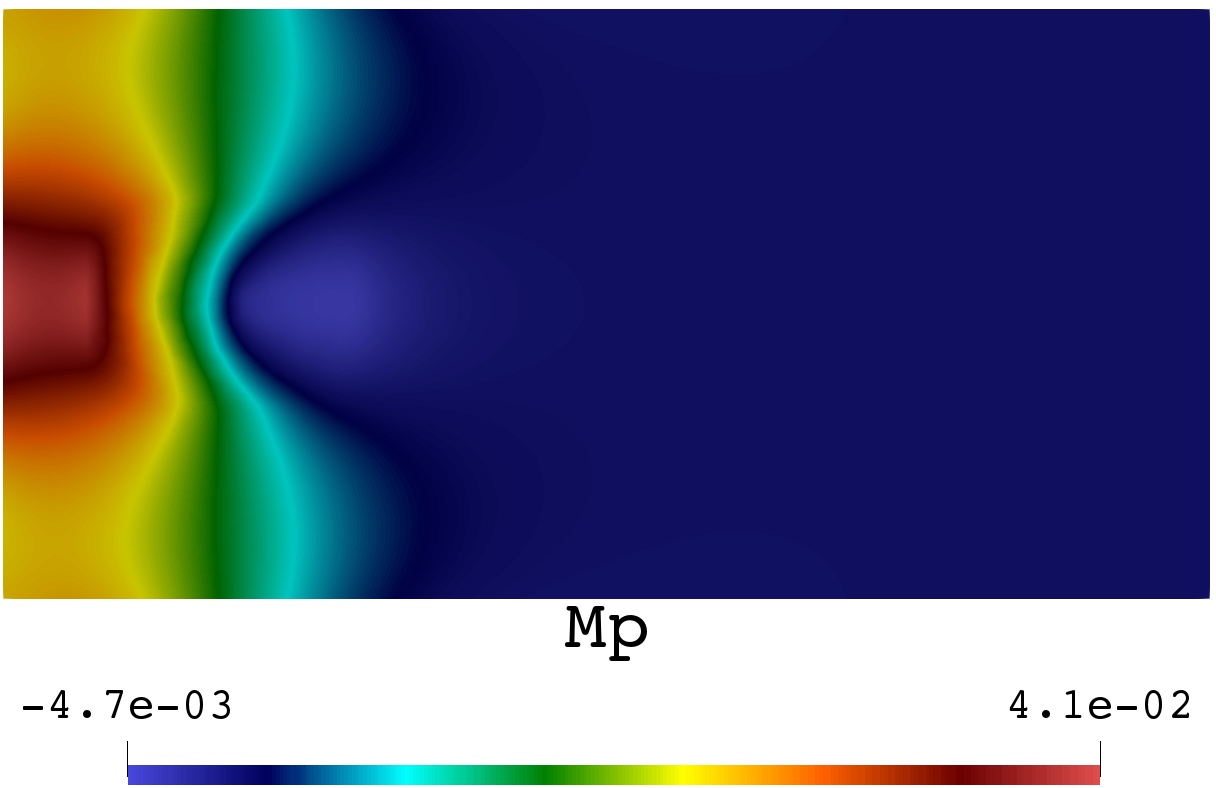}
\end{minipage}
\begin{minipage}{0.24\textwidth}
  \includegraphics[width=\textwidth]{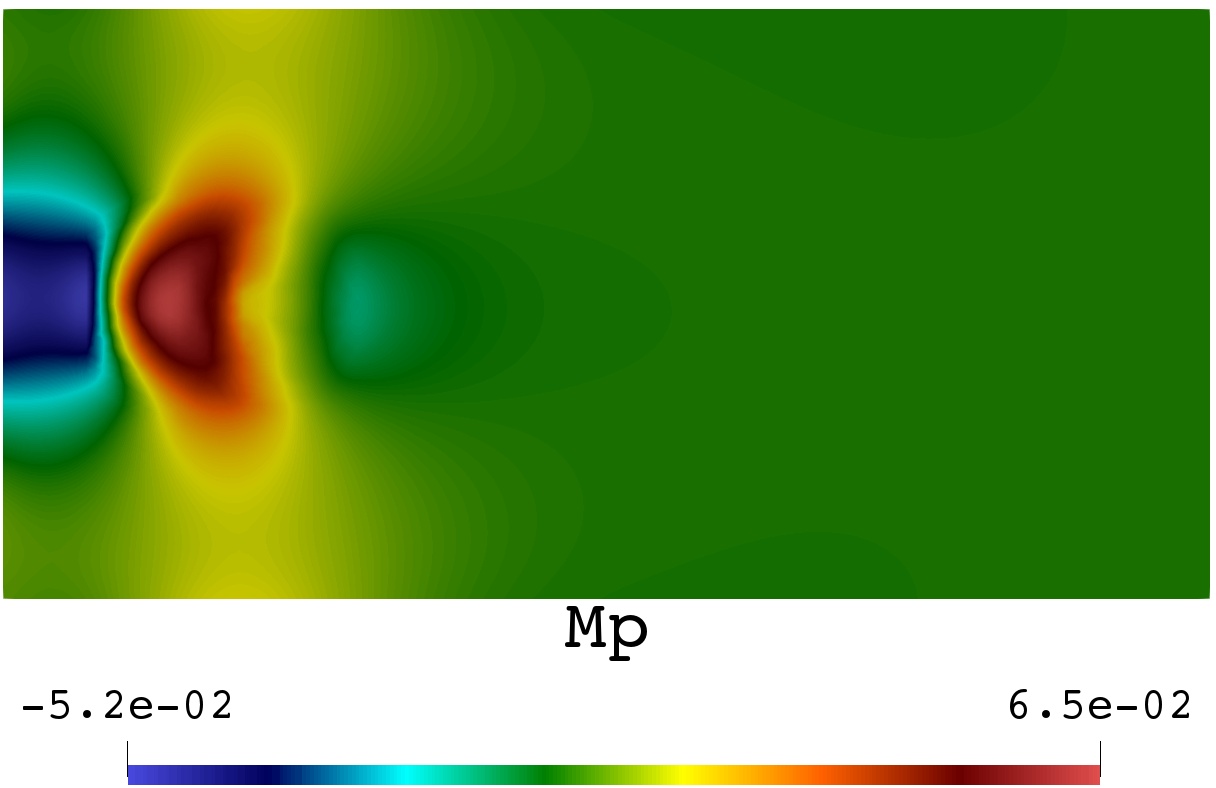}
\end{minipage}
\begin{minipage}{0.24\textwidth}
  \includegraphics[width=\textwidth]{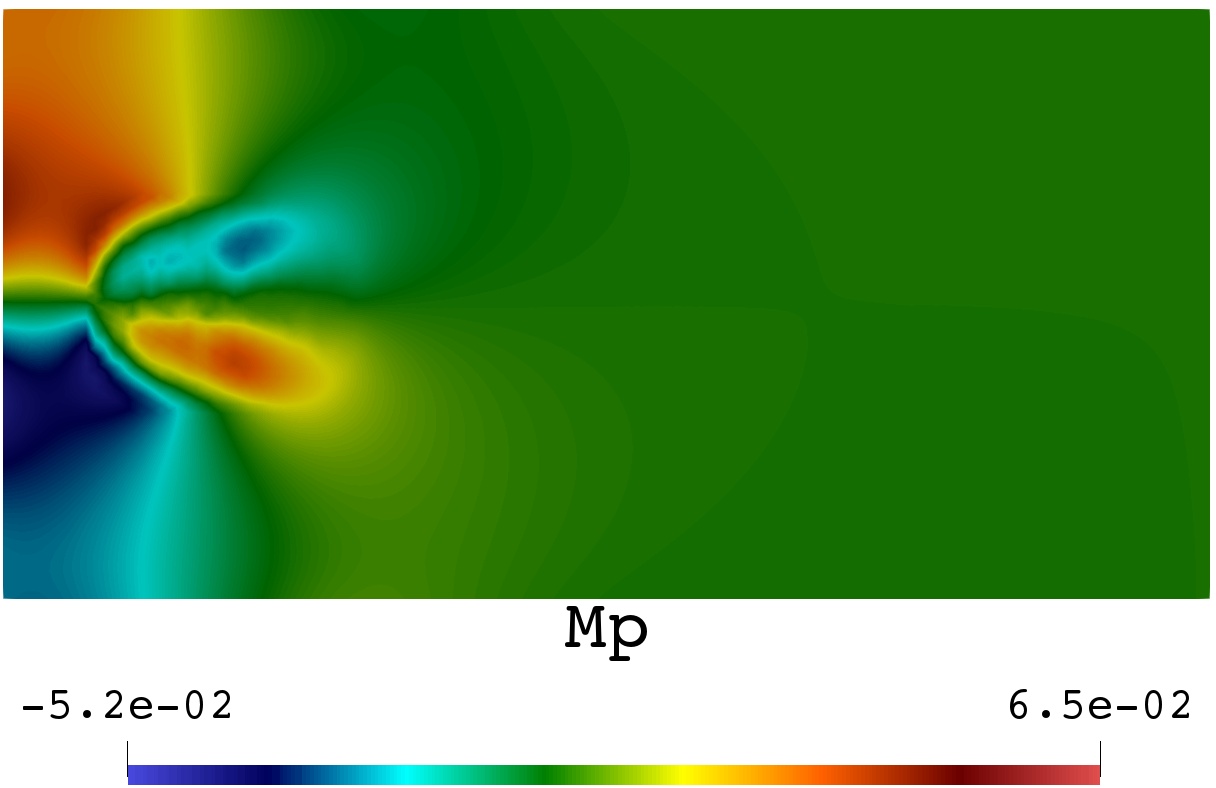}
\end{minipage}
\begin{minipage}{0.24\textwidth}
  \includegraphics[width=\textwidth]{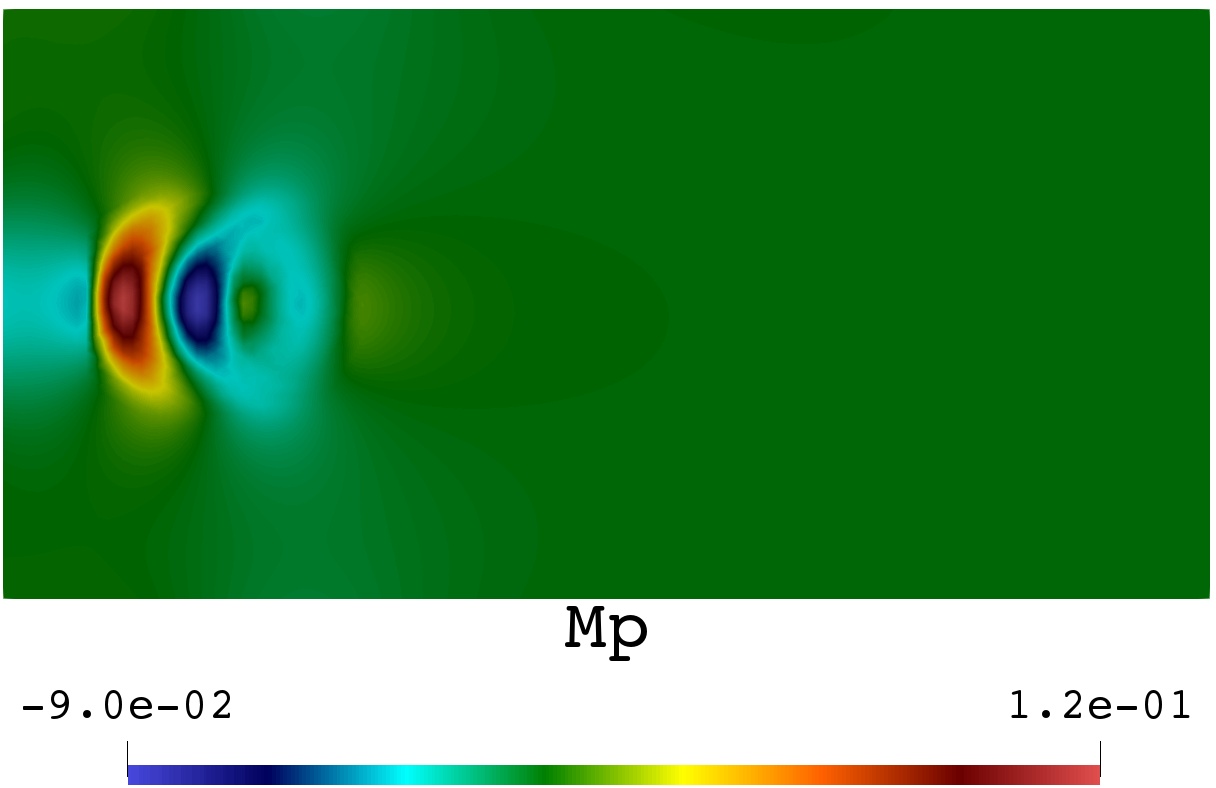}
\end{minipage}
\end{minipage}
\caption{Stokes system ROM-SBM  {POD velocity and pressure components for} (i) $\mu_1$-Geometry and (ii)  $(\mu_0,\mu_1)$-Geometry.}
\label{fig:2D_StokesSBM_modes}
\label{subsec:sup_enrich}
\end{figure}

\begin{figure}
\centering (i)\\
\begin{minipage}{0.32\textwidth}
  \includegraphics[width=\textwidth]{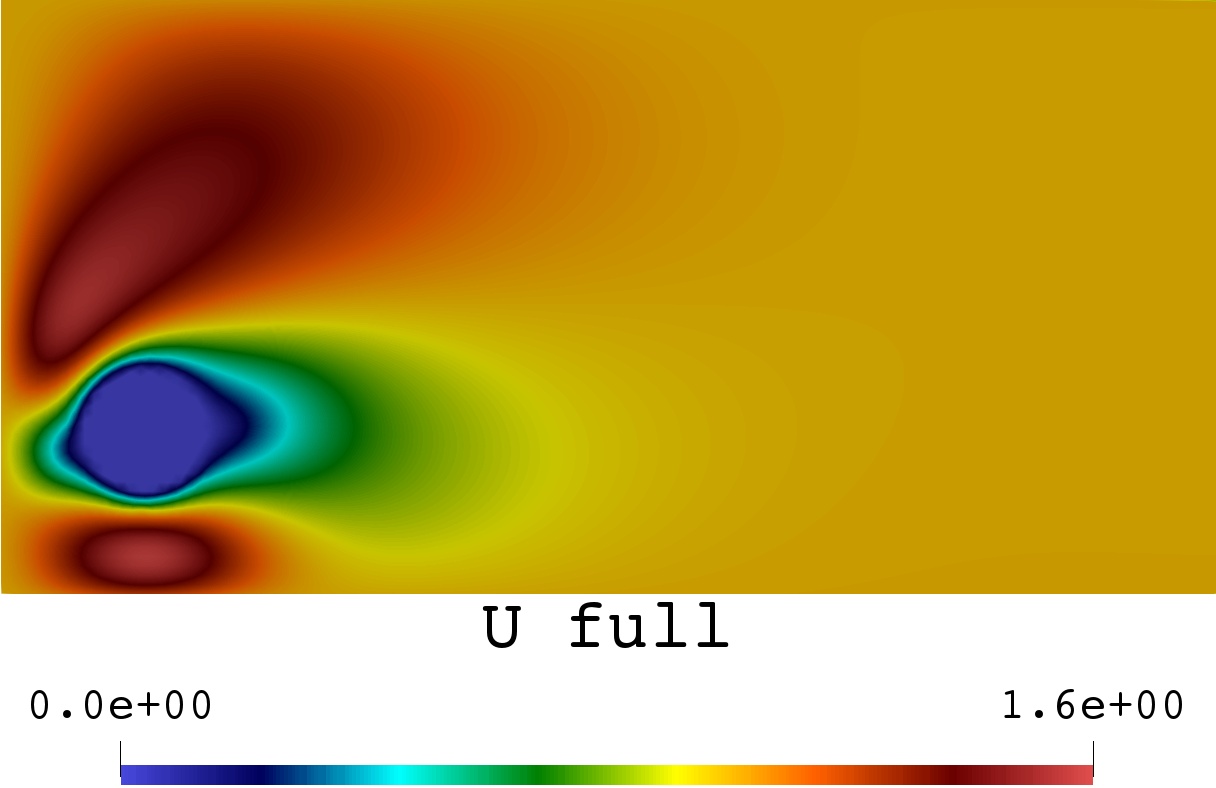}
\end{minipage}
\begin{minipage}{0.32\textwidth}
  \includegraphics[width=\textwidth]{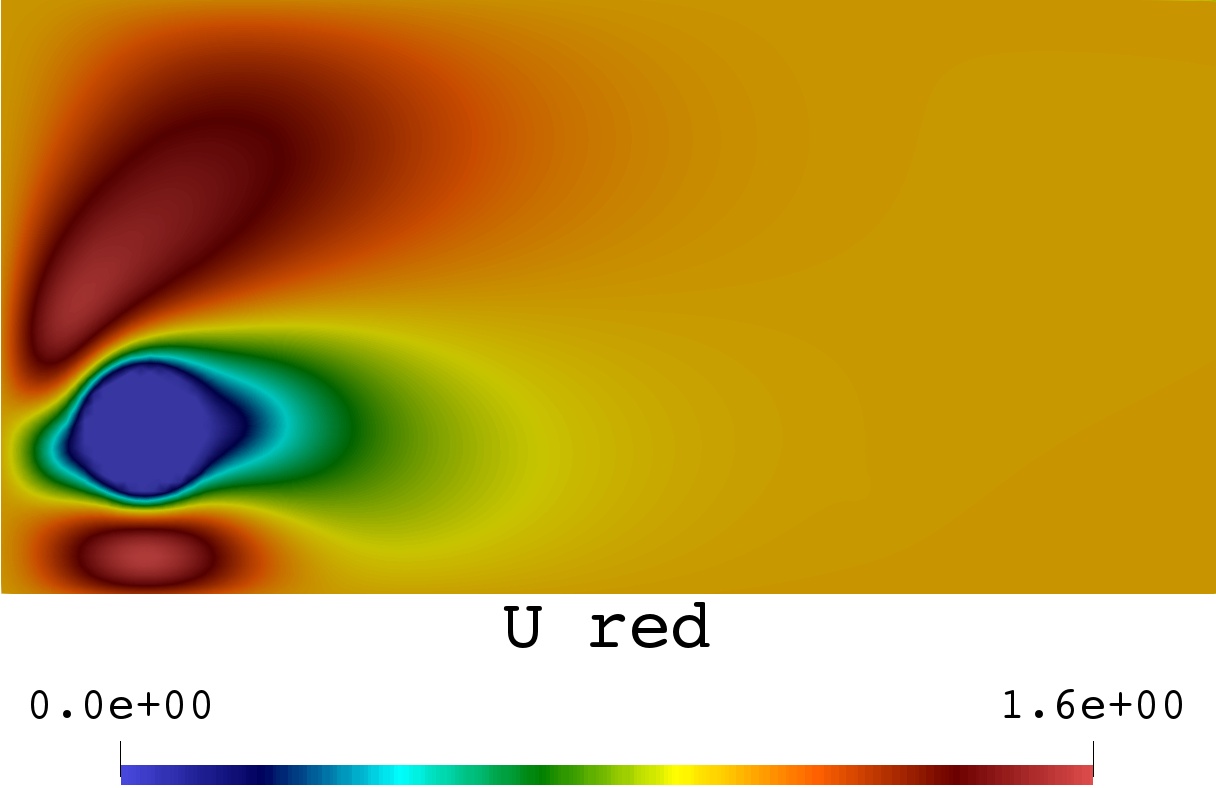}
\end{minipage}
\begin{minipage}{0.32\textwidth}
  \includegraphics[width=\textwidth]{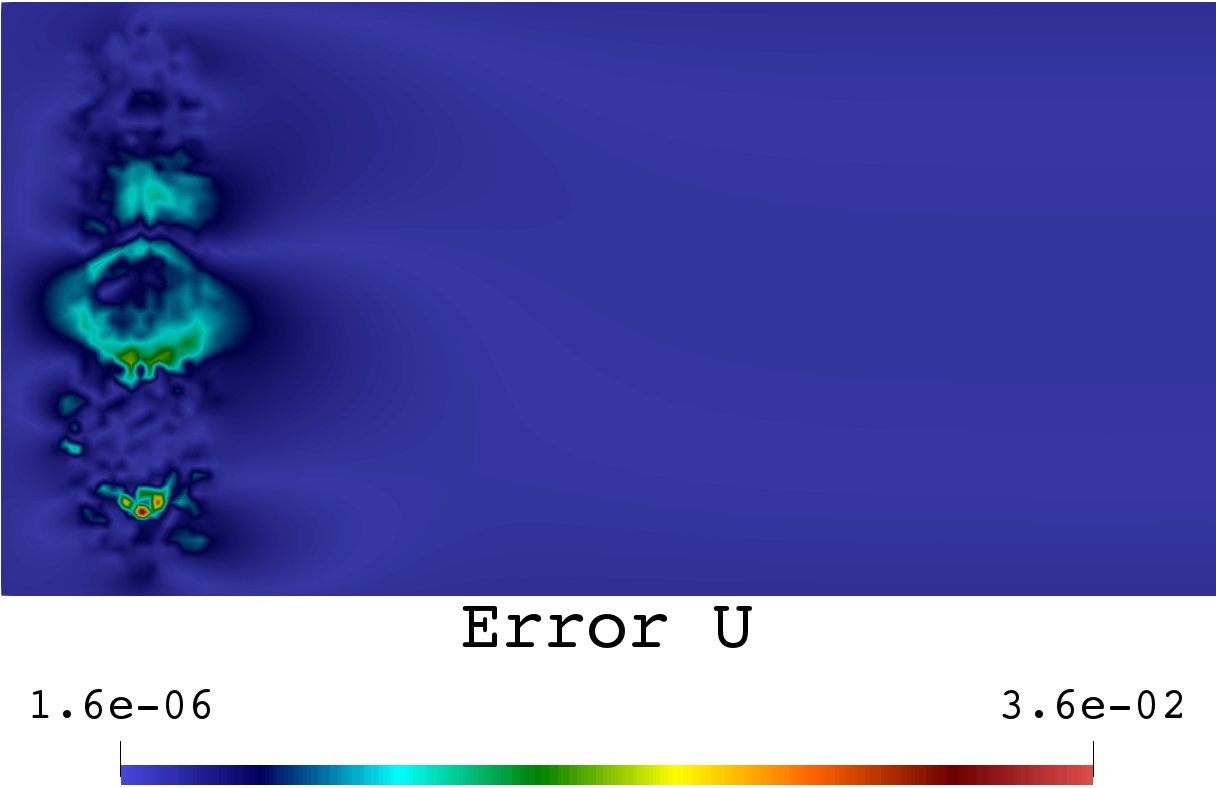}
\end{minipage}
\\\centering (ii)\\
\begin{minipage}{0.32\textwidth}
  \includegraphics[width=\textwidth]{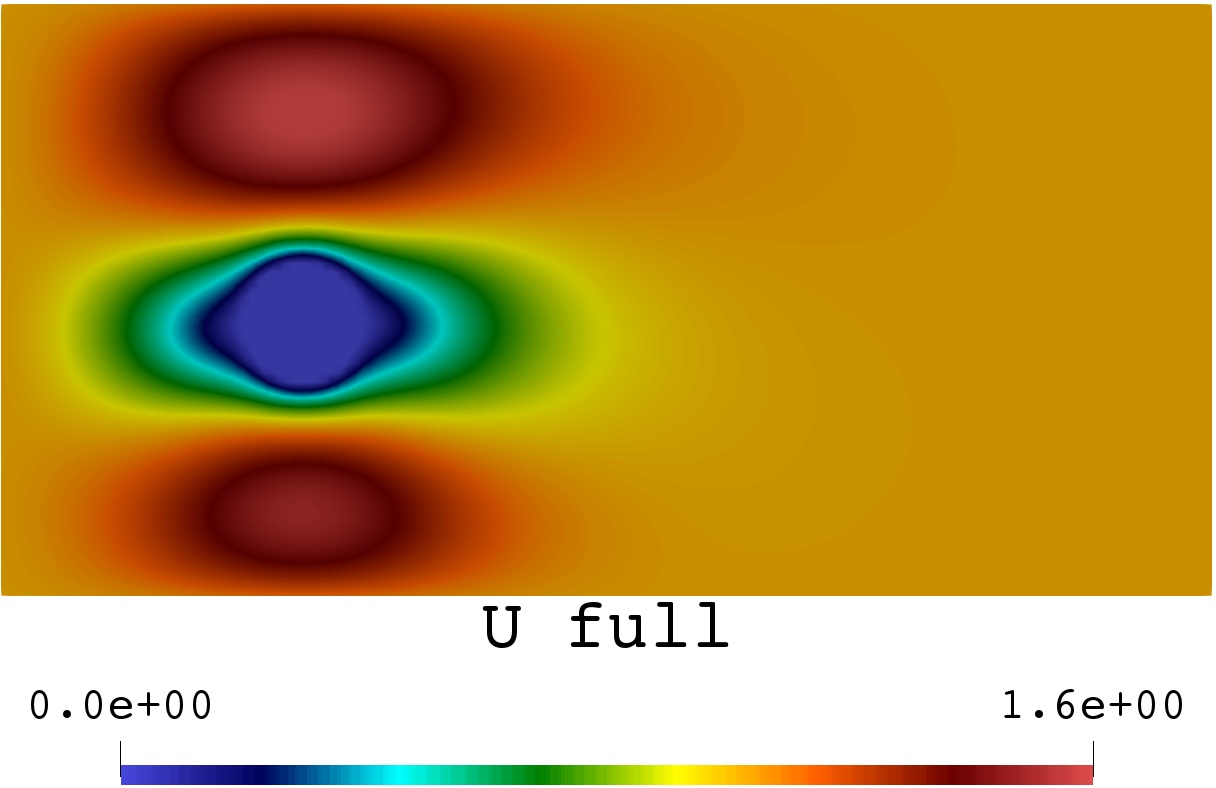}
\end{minipage}
\begin{minipage}{0.32\textwidth}
  \includegraphics[width=\textwidth]{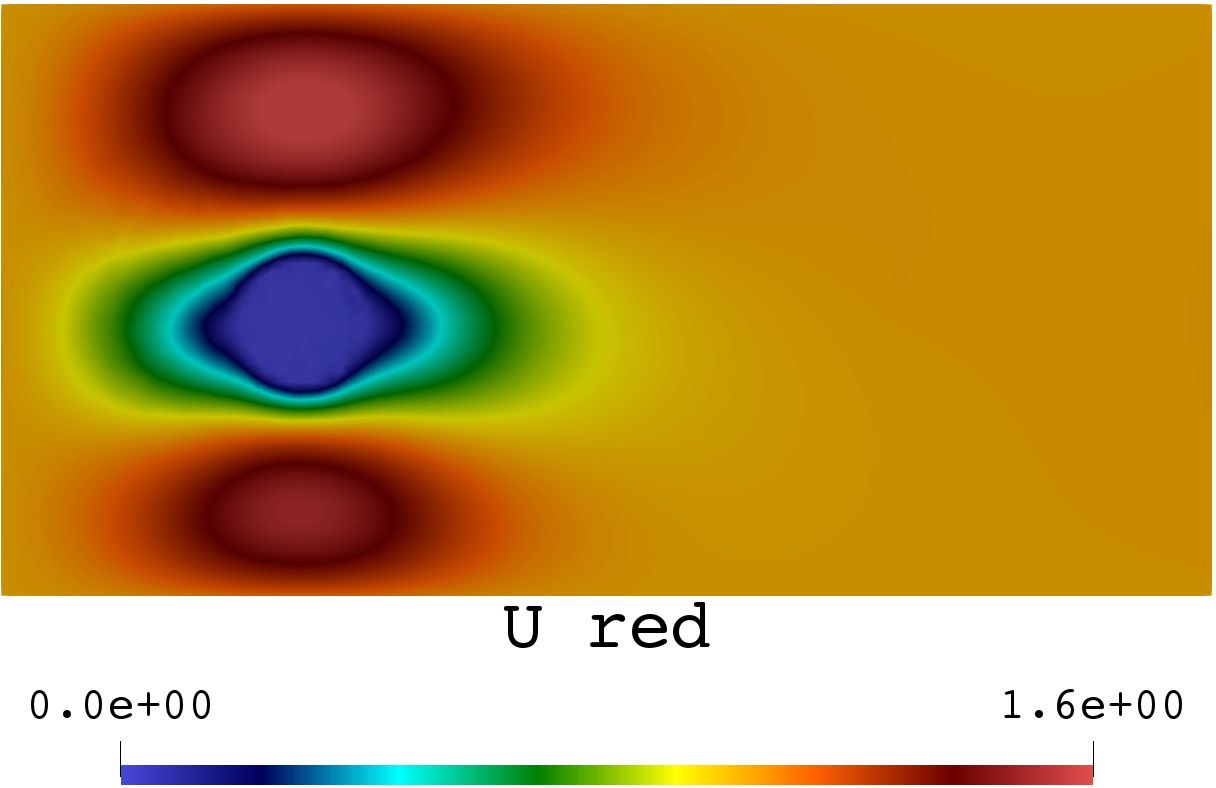}
\end{minipage}
\begin{minipage}{0.32\textwidth}
  \includegraphics[width=\textwidth]{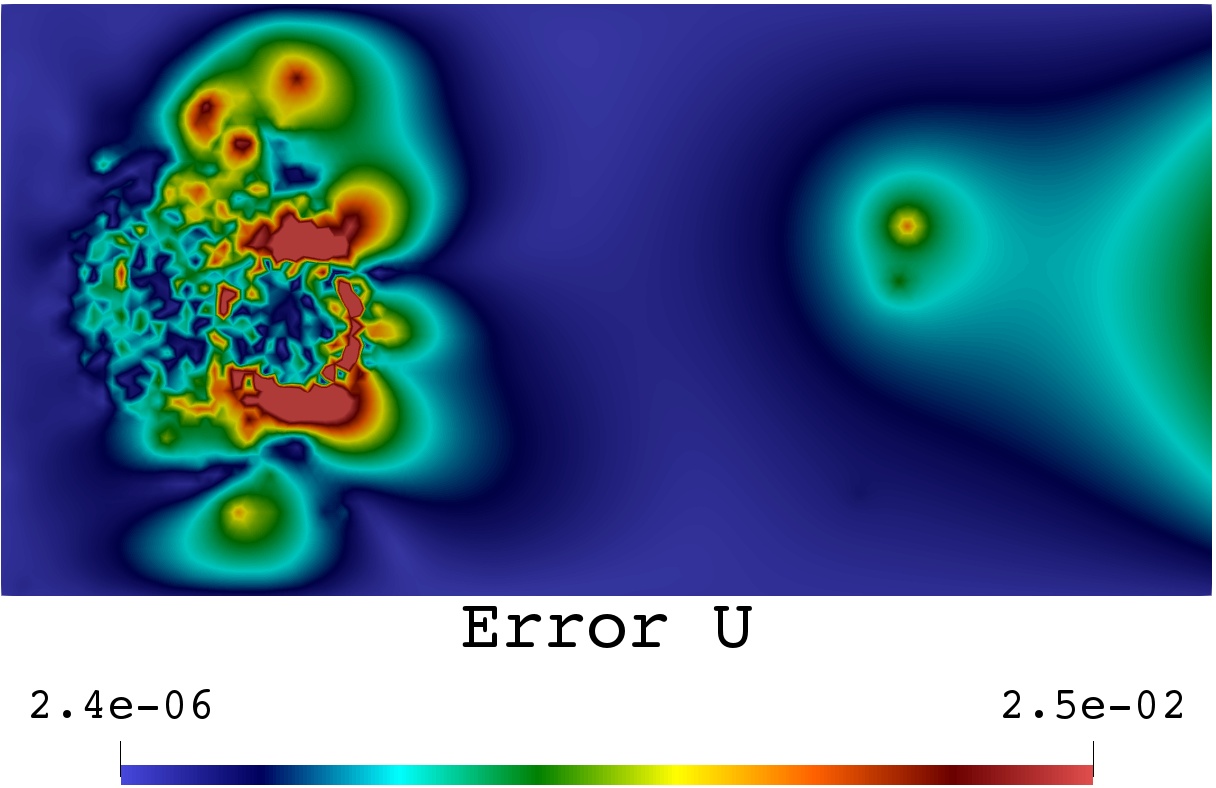}
\end{minipage}
\\
(i)
\\
\begin{minipage}{0.32\textwidth}
  \includegraphics[width=\textwidth]{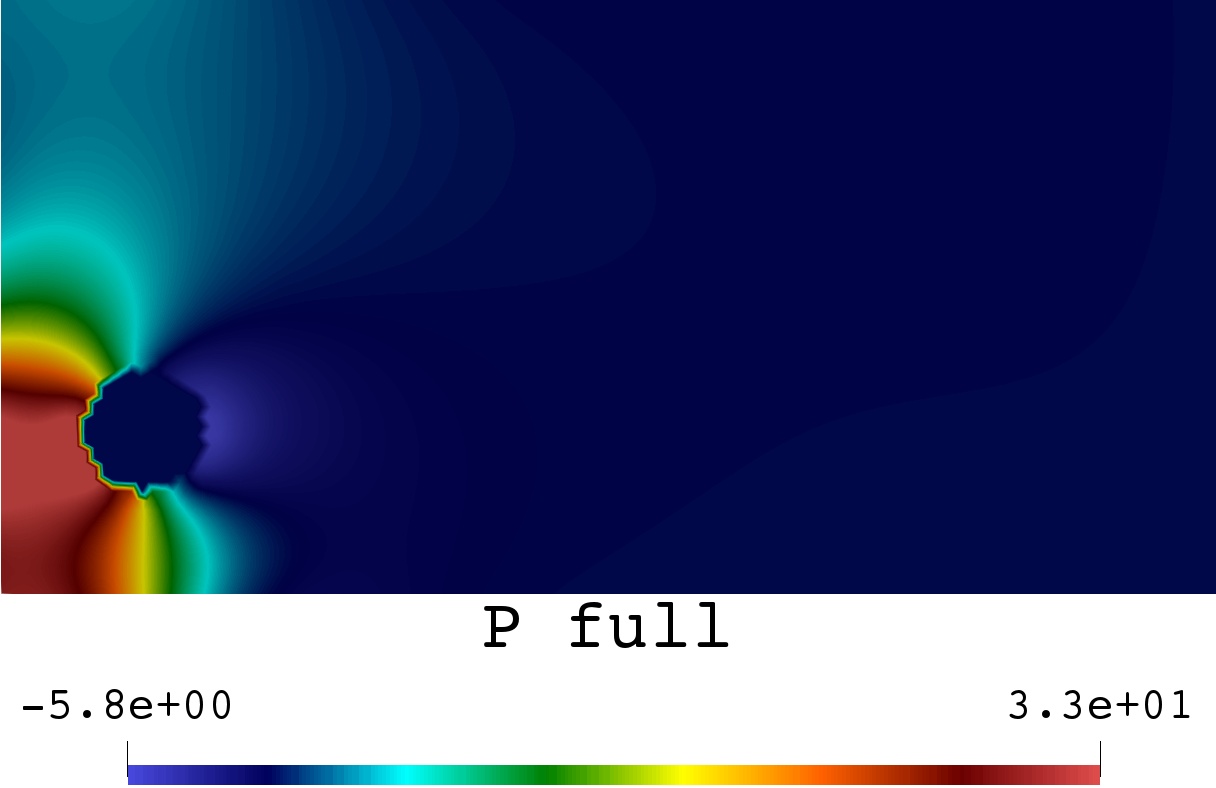}
\end{minipage}
\begin{minipage}{0.32\textwidth}
  \includegraphics[width=\textwidth]{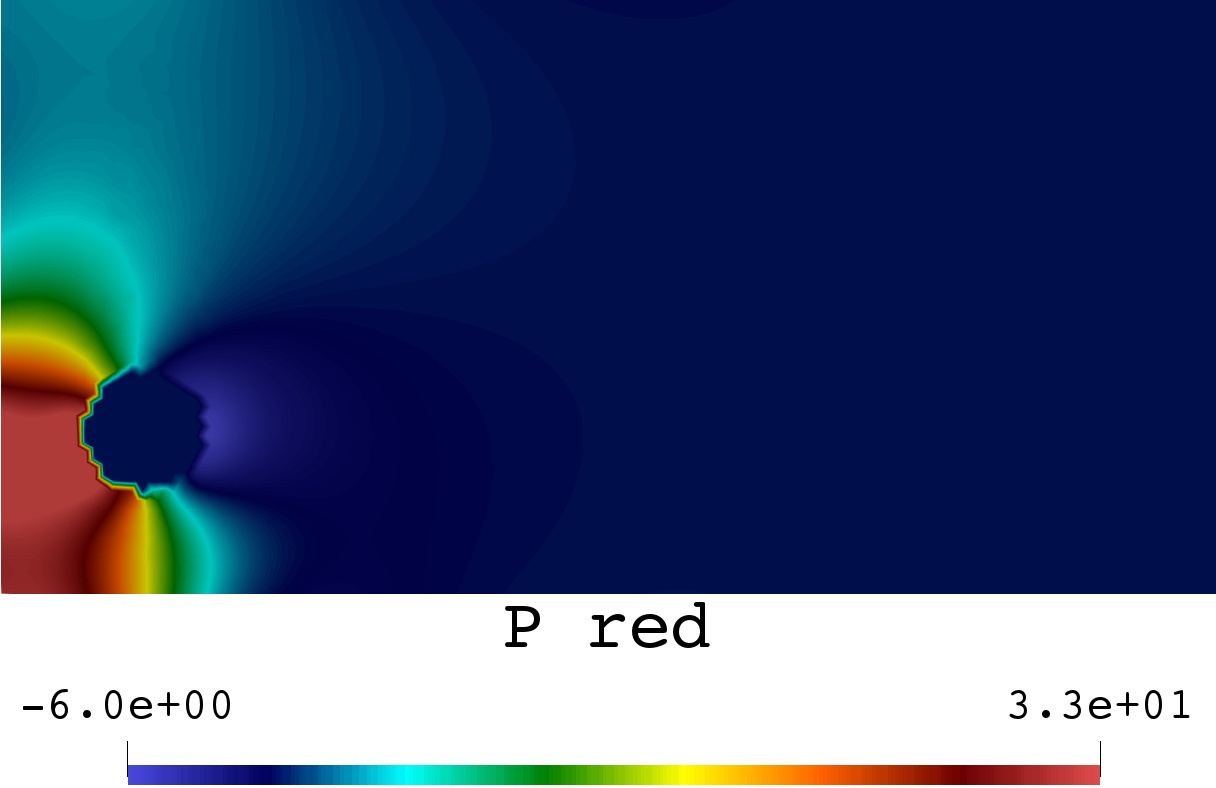}
\end{minipage}
\begin{minipage}{0.32\textwidth}
  \includegraphics[width=\textwidth]{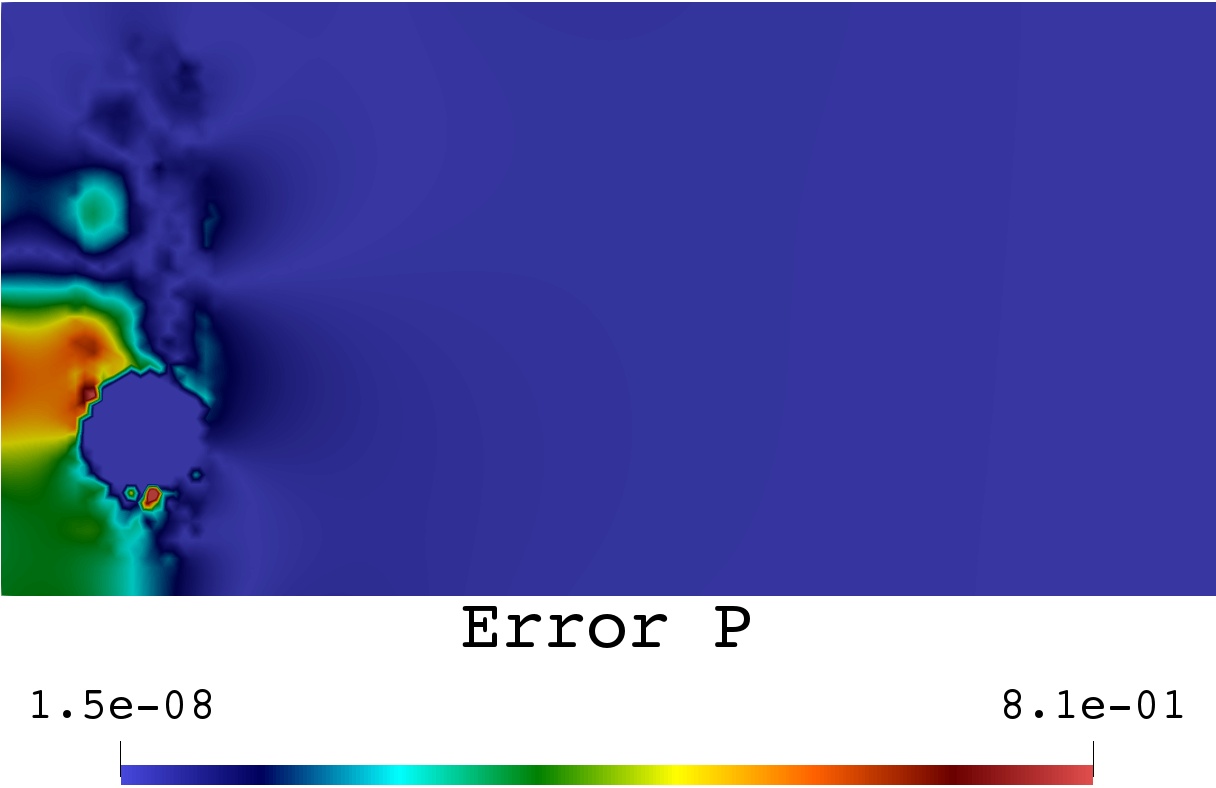}
\end{minipage}
\\
(ii)
\\
\begin{minipage}{0.32\textwidth}
  \includegraphics[width=\textwidth]{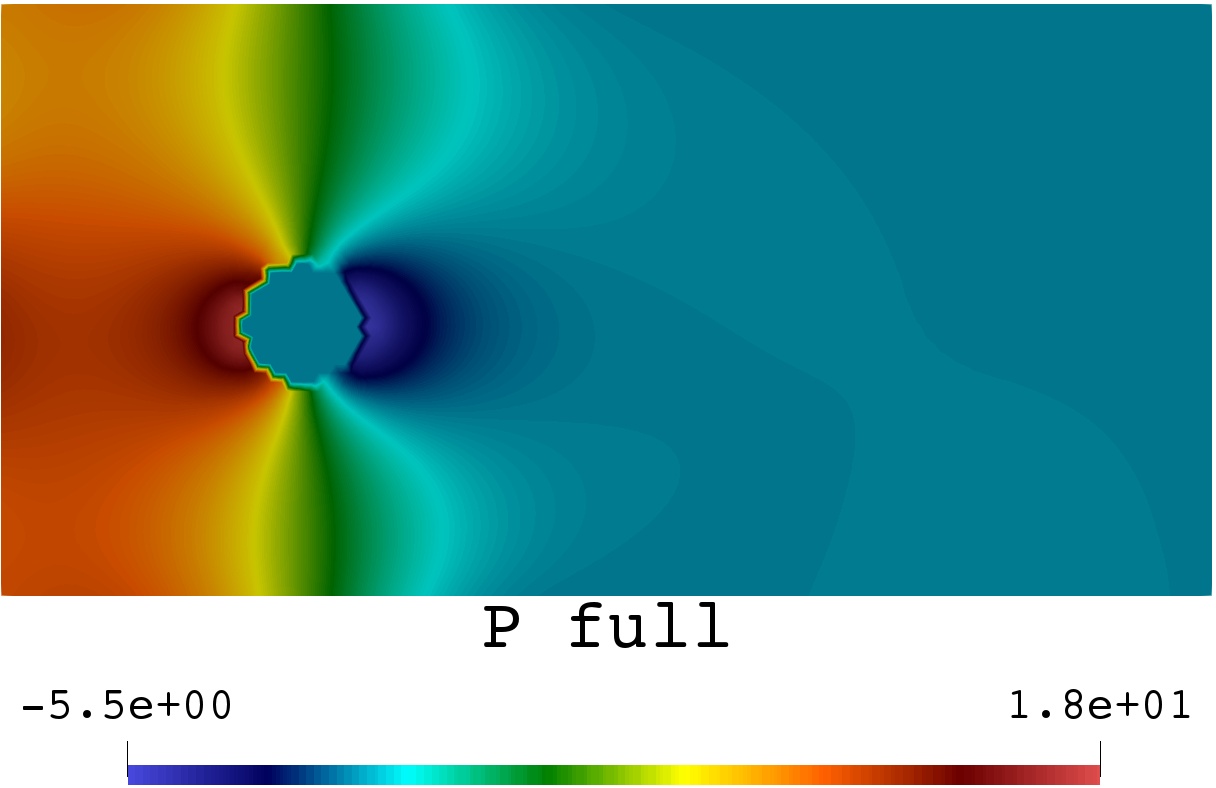}
\end{minipage}
\begin{minipage}{0.32\textwidth}
  \includegraphics[width=\textwidth]{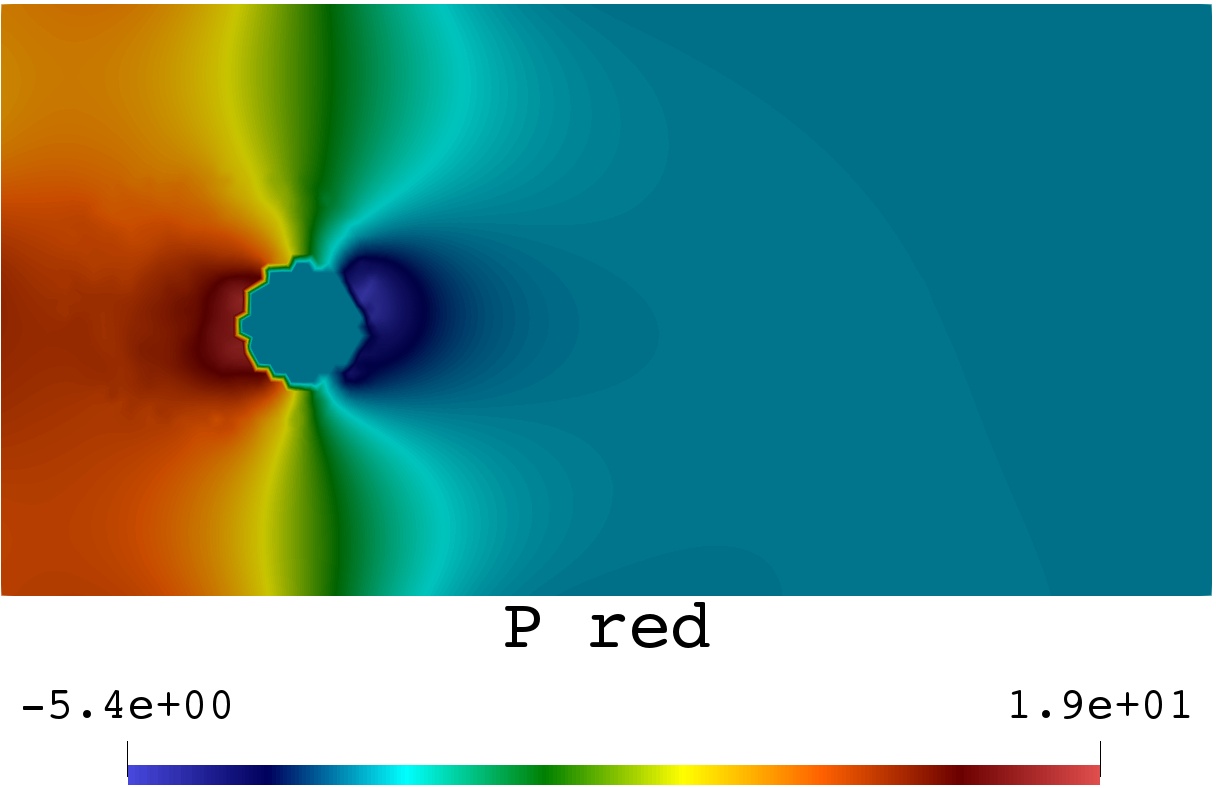}
\end{minipage}
\begin{minipage}{0.32\textwidth}
  \includegraphics[width=\textwidth]{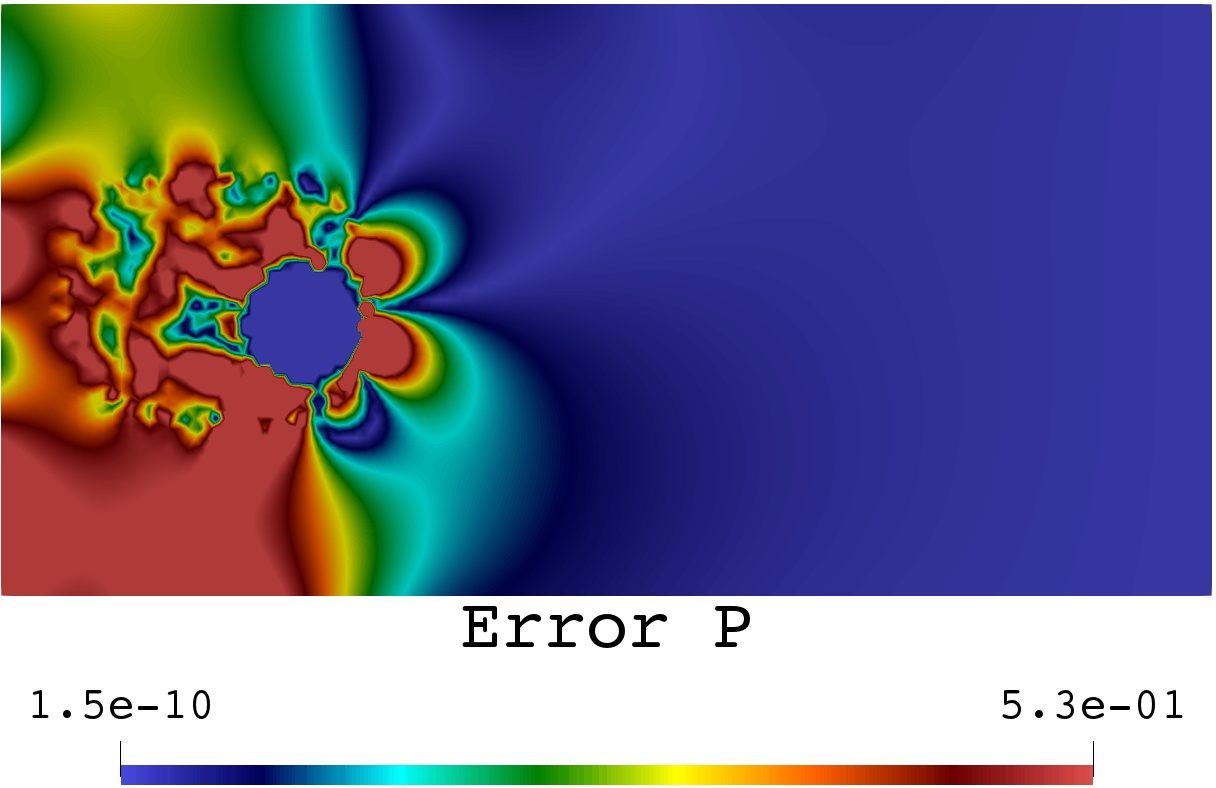}
\end{minipage}
\caption{
Stokes SBM: the full, the reduced solution and the absolute error  {for velocity and pressure for (i) 1D and (ii) 2D parametrization case}.
}
\label{FULL_RED_ERROR_2P}
\label{FULL_RED_ERROR_1P}
\end{figure}

\begin{figure} \label{fig:comp_time} 
\centering
\vskip-40pt
\begin{minipage}{\textwidth}
\centering
(i)$_a$
\begin{minipage}{0.49\textwidth}
\includegraphics[width=\textwidth]{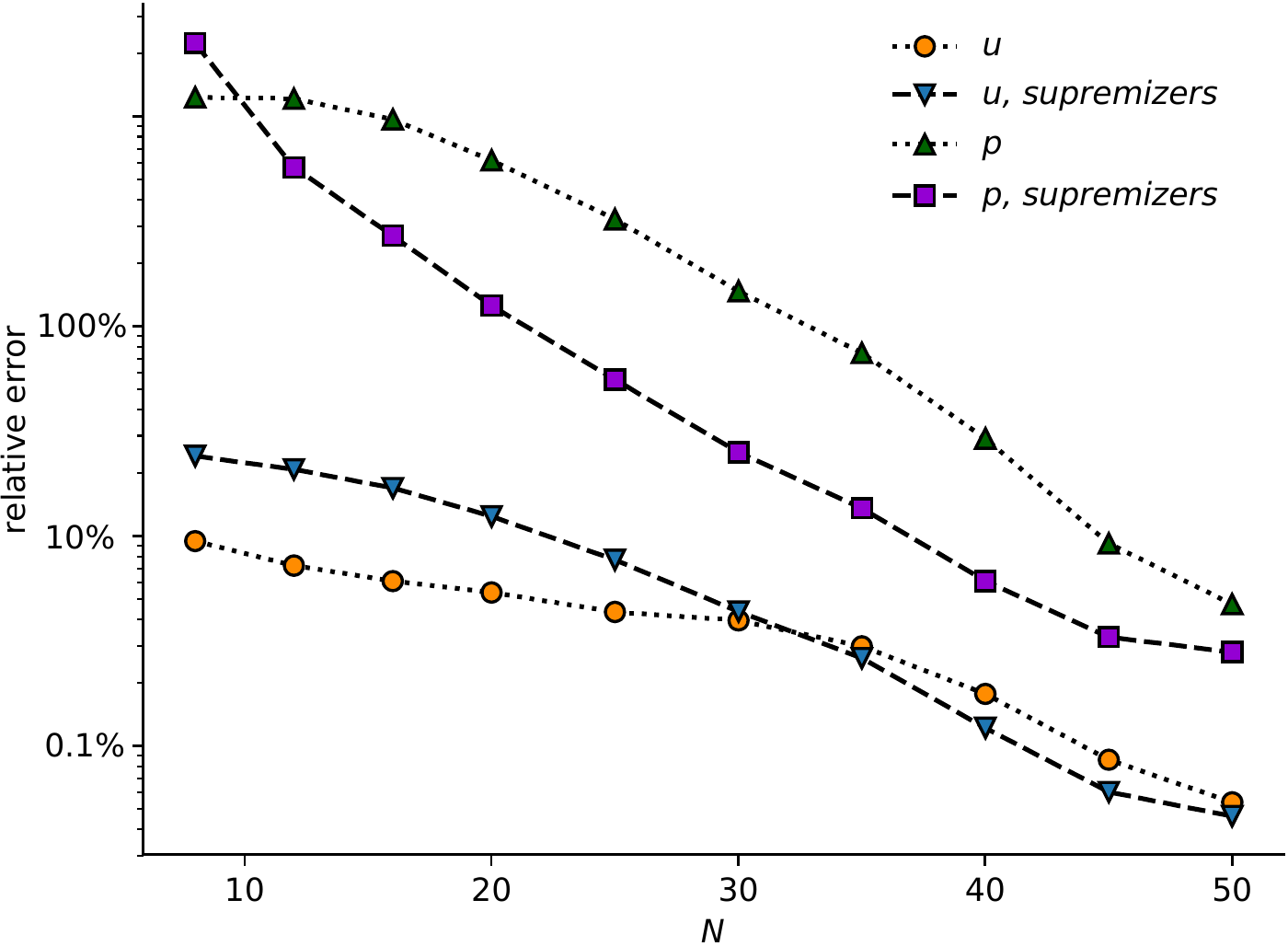}
\end{minipage}
\\
(i)$_b$
\begin{minipage}{0.49\textwidth}
\includegraphics[width=\textwidth]{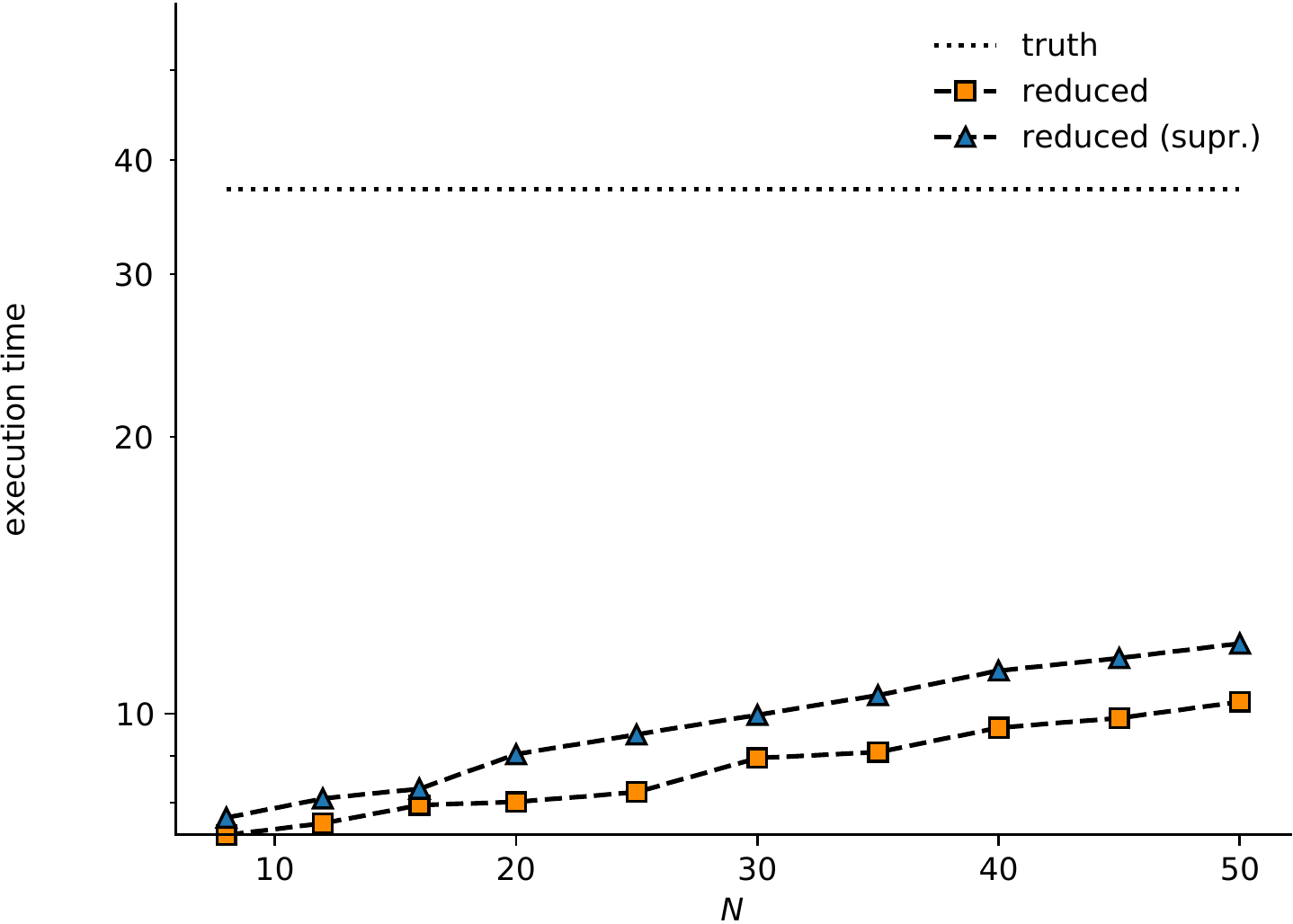}
\end{minipage}
\\
(ii)$_a$
\begin{minipage}{0.49\textwidth}
\includegraphics[width=\textwidth]{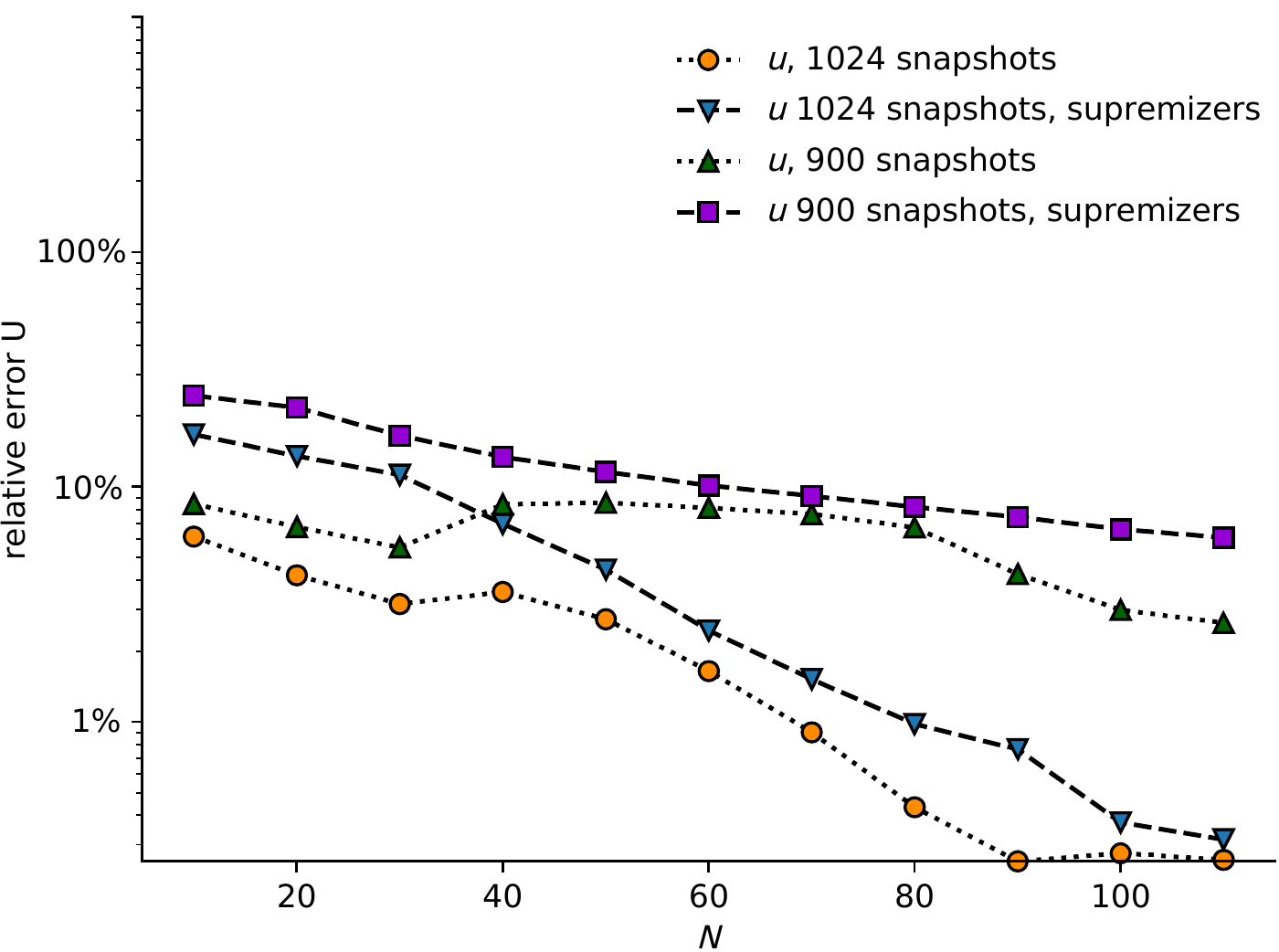}
\end{minipage}
\\
(ii)$_b$
\begin{minipage}{0.49\textwidth}
\includegraphics[width=\textwidth]{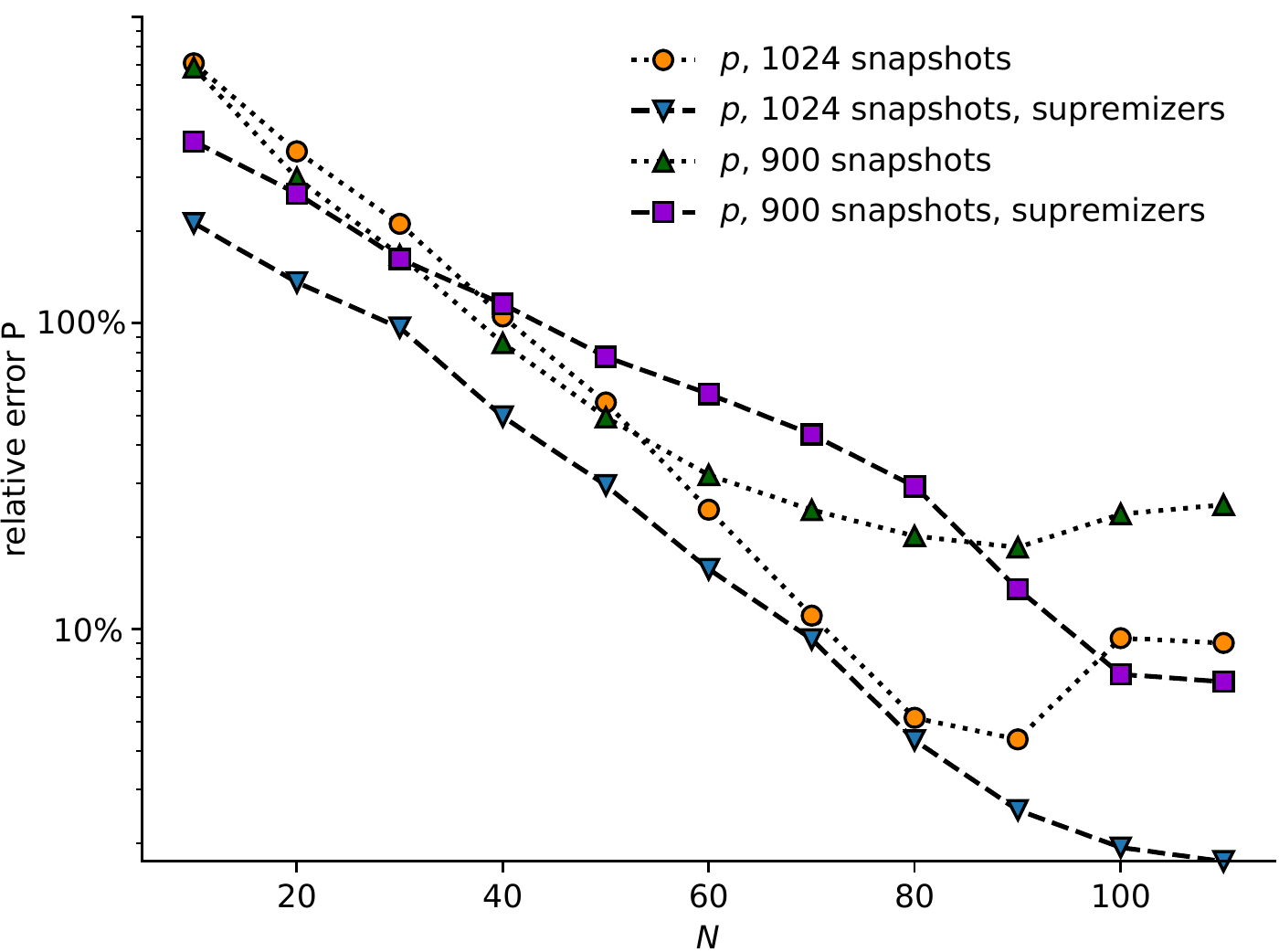}
\end{minipage}
\end{minipage}\vskip-8pt
\caption{Stokes (SBM) parametrization: (i) $\mu_1$-Geometry, the relative errors for velocity and pressure in (i)$_a$, and the execution times in (i)$_b$, 
(ii) $(\mu_0,\mu_1)$-Geometry, velocity and pressure field with and without supremizer stabilization.}
\label{fig:crop_up_stokes_SBM}
   \label{fig:timers} \label{fig:WithoughtSupremizer_solution_1000_100_1000_100_1024_j_3J_4j}\label{fig:Supremizer_solution_1000_100_1000_100_1024_j_3J_4j}
\end{figure}

The data of the present experiment consider a steady Stokes flow around an embedded circular cylinder, within the framework of a parametrized embedded domain described by the level set: 
$
(x-\mu_0)^2 +(y-\mu_1)^2 \le R^2.
$
We will consider two different geometrical parametrization test cases,  
 {a one dimensional} with parameter $\mu_1$,  
 and  {a} two dimentional  with two parameters  
 $\mu_0$ and $\mu_1$. 
The data of the problem are: 
an embedded cylinder with constant radius $0.2$,  
viscosity $\nu=1$, force $f=1$, in the $x$ direction, {$u_{\text{in}} = 1$} on the left side of the domain, an open boundary condition with {$p_{\text{out}} = 0$}  on the right. 
A slip no penetration boundary condition is applied on the top and bottom edges.
On the boundary of the embedded cylinder a no slip boundary condition is applied.
The mesh size is $h = 0.0350$ for the background mesh, using 15022 triangles, $\mathbb{P}1/\mathbb{P}1$ finite elements in space with stabilization terms,
with and without {supremizer} basis enrichment  {in the offline stage}.
\begin{figure} \centering
\includegraphics[width=0.4\textwidth]{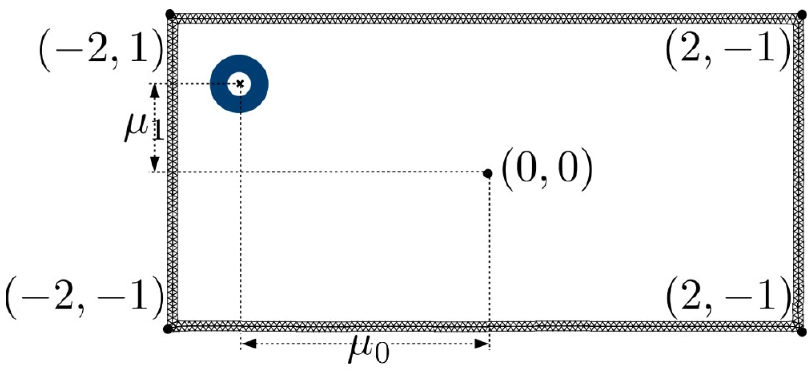}
 \caption{Stokes with Shifted Boundary Method experiment:  sketch mesh, the embedded domain, and the parameters considered in the numerical examples.}
 \label{fig:Stokes_SBM:background_mesh}
\end{figure}
Some reduced basis components derived with the proper orthogonal decomposition for the Stokes system discretized by the Shifted Boundary Method for the
1D and 2D geometrical parametrization with $\mu_1 \in [-0.65, 0.65]$ are visualized in Figure \ref{subsec:sup_enrich}.
In  {parallel, FOM}, ROM solution and absolute error values are visualized in Figure 
\ref{FULL_RED_ERROR_2P}  {for both experiments}. 
 {Relative} errors ($||u-u_r||_{L^2}/ ||u||_{L^2}$ and $||p-p_r||_{L^2}/ ||p||_{L^2}$), as well as, the execution times for the 1D parametrization  are reported in Tables \ref{table:errors_yes_no_supremizers}, \ref{table:timers}, Figure \ref{fig:timers} (i), while for the 2D case the relative errors in Tables \ref{table:2Dparametrization} and Figure \ref{fig:crop_up_stokes_SBM} (ii).
\begin{table} 
 \caption{ {Stokes (SBM): relative errors between the full order solution and the reduced basis solution, 1D geometrical parametrization.}}
\centering
\scalebox{0.85}{
  \begin{tabular}{|c|c|c|c|c|}
  \hline
      {Snapshots:} & \multicolumn{2}{c|}{1024} & \multicolumn{2}{c|}{1024}   \\
   \hline
     Suprem.:& \multicolumn{2}{c}{No} & \multicolumn{2}{|c|}{Yes}\\
    \hline
    Modes & rel. error u & rel. error  p & rel. error u & rel. error  p \\
    \hline
     8  & 0.0947158 &  12.309881  &  0.2406999 &  22.319781\\\hline
     12 & 0.0723268 &  12.133591  &  0.2078557 &  5.7159319\\\hline
     16 & 0.0610052 &  9.6652163  &  0.1692787 &  2.6962056\\\hline
     20 & 0.0538906 &  6.1692750  &  0.1243368 &  1.2535779\\\hline
     25 & 0.0434925 &  3.2331644  &  0.0770726 &  0.5568314\\\hline
     30 & 0.0396132 &  1.4693532  &  0.0437348 &  0.2504069\\\hline
     35 & 0.0298269 &  0.7455038  &  0.0262345 &  0.1356788\\\hline
     40 & 0.0177170 &  0.2918072  &  0.0121903 &  0.0611154\\\hline
     45 & 0.0085905 &  0.0923509  &  0.0060355 &  0.0330206\\\hline
     50 & 0.0053882 &  0.0473412  &  0.0046300 &  0.0279857\\
    \hline
  \end{tabular}
}
 \label{table:errors_yes_no_supremizers}
\end{table}
Focusing on the supremizer stabilization approach,
the improved results for pressure are obvious and the plots clearly show in a glance that the high fidelity and ROM solutions cannot be easily distinguished.
Related to the execution times investigation for the 1D geometrical parametrization we compare 
online stage execution times, 
against the full order computational times 
and FOM solutions. Namely, for  {all} $1024$ snapshots and supremizers the solutions costs time $1$h $42$m and $54$s. 
This is an expensive stage but fortunately it is executed once in the beginning. 
We clarify that in the online stage computation time, it is included: the assembling of the full order matrices, the production of the reduced order model and  its resolution (see Table \ref{table:timers}).
\begin{table} \centering
\caption{ {Execution time, at the reduced order level, for the case with 1D geometrical parametrization.}}
{\scalebox{0.85}{
 \begin{tabular}{|c|c|c|c|c|}
    \hline
     & no supr. &  supr.   \\
     Modes       & exec. time (sec) &  exec. time (sec) \\
    \hline
     8   &  7.3858961  &   7.710907   \\\hline
     12 &  7.6042165  &   8.091225  \\\hline
     16 &  7.9584049  &   8.290780   \\\hline
     20 &  8.0206915  &   9.036709    \\\hline
     25 &  8.2229143  &   9.495323   \\\hline
     30 &  8.9529275  &   9.972288  \\\hline
     35 &  9.0867916  &   10.47633   \\\hline
     40 &  9.6555775  &   11.13931    \\\hline
     45 &  9.8934008  &   11.49422   \\\hline
     50 &  10.302459  &   11.92024   \\
    \hline
  \end{tabular}
}}
\label{table:timers}
\end{table}
Ten different values of the input parameter are considered, with reference time the time execution at FOM level which  {for each one parameter solution} is equal to $\approx 37$ sec. 
Obviously, the ROM leads to a considerable speed-up for the all different analyzed configurations and for both cases with and without supremizer enrichment.  
The interest in the latter experiments was more into testing the feasibility and the accuracy of a reduced order model, constructed starting from a shifted boundary FOM, and we did not employed any hyper reduction  {technique, which} means that, also at the reduced order level, we assembled the full order discretized differential operators. 
Focused on the $( {\mu_0},\mu_1)$-Geometry experiment with range  $\mu=(\mu_0,\mu_1)=[-1.5,-1.0]\times[-0.15,0.15]$, 
which is a more complex and demanding scenario, 
we can easily notice that the supremizer enrichment is pointed out necessary for convergence stability and reliable pressure results, 
and without a supremizers enrichment and $900$ snapshots  
the best achieved -and  {disappointing}- relative error were equal to $0.0263821$ and $0.1854861$ for velocity and pressure.
Figure \ref{fig:crop_up_stokes_SBM} and Table \ref{table:2Dparametrization} also give an overview of the number of snapshots dependence in two choices of  $900$ and $1024$ as they are used in the offline stage.
\begin{table} 
\caption{ {Stokes SBM $(\mu_0,\mu_1)$-Geometry: Supremizer basis enrichment and the relative error.}}
\centering
\scalebox{0.75}{
  \begin{tabular}{|c|c|c|c|c|
}
   \hline
     Snapshots: & \multicolumn{2}{c|}{900} & \multicolumn{2}{c|}{1024}   \\
       \hline
      {Suprem.:}& \multicolumn{2}{c}{Yes} & \multicolumn{2}{|c|}{Yes}\\
    \hline
    Modes & rel. error u & rel. error  p & rel. error u & rel. error  p 
\\
    \hline
     10 & 0.2448511 &  3.9240637   &   0.1672753   &   2.1243228 
\\\hline
     20 & 0.2175821 &  2.6531343   &   0.1353706   &   1.3611011 
\\\hline
     30 & 0.1652331 &  1.6234701   &   0.1124619   &   0.9680506 
\\\hline
     40 & 0.1340978 &  1.1560352   &   0.0696437   &   0.4958605 
\\\hline
     50 & 0.1158443 &  0.7777786   &   0.0444991   &   0.2958338 
\\\hline
     60 & 0.1013961 &  0.5876048   &   0.0244793   &   0.1574037 
\\\hline
     70 & 0.0914650 &  0.4335489  &    0.0151749   &   0.0928402 
\\\hline
     80 & 0.0822658 &  0.2933336  &    0.0097848   &   0.0434299
\\\hline
     90 & 0.0744696 &  0.1355488   &   0.0076431   &   0.0257060
\\\hline
     100& 0.0660493 &  0.0714350   &   0.0037280   &   0.0194051  
\\\hline
     110& 0.0609040 &  0.0675720 &     0.0031577   &   0.0174815
\\
    \hline
  \end{tabular}
}
  \label{table:2Dparametrization}  \label{table:errors_supremizers}
\end{table}
\section{Searching a better ROM with CutFEM}
\begin{figure} \centering
\includegraphics[width=0.99\textwidth]{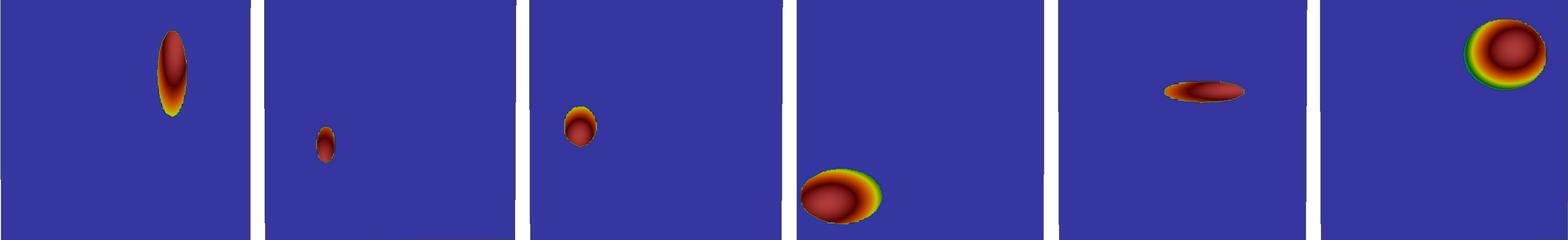}
\\ \vskip2pt
\includegraphics[width=0.3\textwidth]{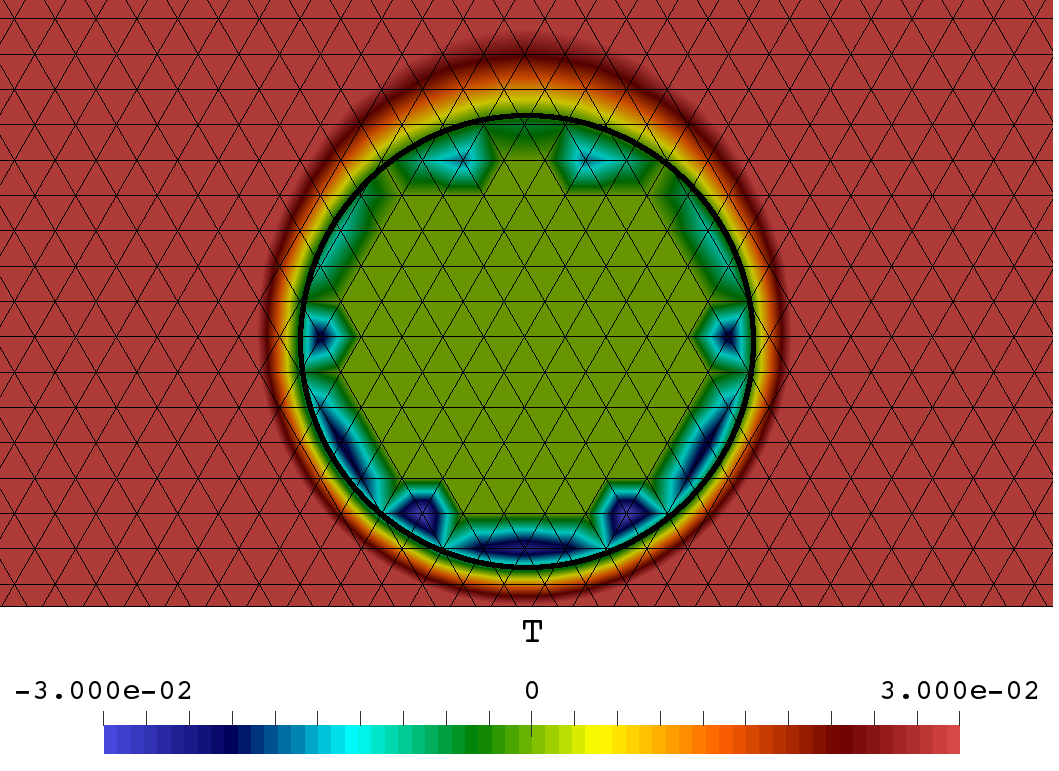}
\label{fig:poisson_zoom}
  \caption{Shape parametrization with large deformations and a zoom into the embedded cylinder visualizing the extended solution.}
\label{fig:extreme_def}
\end{figure}
In embedded methods, the reduced order basis is constructed on the whole background domain and
great care  {is needed in} the manipulation  of the out of interest -outside- the truth geometry 
area, the so called ''ghost area``. 
In this  {section} approach, with CutFEM,  again we employ 
the solution values of the surrogate boundary as they were computed after using the natural smooth extension
from the truth to the surrogate domain,
that allows a smooth extension of the solution to the neighboring ghost elements with values which are decreasing smoothly to zero.   
The latter is guaranteeing a regular ``solution'' in the background domain and provides the construction of a good reduced order basis and a much more promising good reduced order approximation.
Next we present the weak formulation for the Poisson equation introducing the Cut Finite Element basics. 

For any real parameter vector $\mu$   
we seek   $u(\mu)$ in $$V_{g_D}(\mu)=\left\{w\in H^1\left(\mathcal D(\mu)\right) \text{ with } w|_{\Gamma
(\mu)}=g_D(\mu)\right\}$$
such that $\forall v(\mu)\in V_0(\mu)$
\begin{eqnarray*}
  \left( \nabla u(\mu) ,\nabla \upsilon(\mu)\right) = \left(g(\mu) , \upsilon(\mu) \right)
  + \left(g_N(\mu) , \upsilon(\mu) \right)_{\Gamma_N(\mu)}.\label{P:ContWeakForm}
\end{eqnarray*}

The 
boundary value problem is
 formulated on a domain $\mathcal{T}(\mu)$ that contains {${\mathcal D}(\mu) \subset \mathcal{T}(\mu)$},
 its mesh $\mathcal{T}_h(\mu)$ is not fitted to the domain boundary,
 {${G}_h(\mu):=\{K\in \mathcal{T}_h(\mu): K\cap\Gamma(\mu)\neq\emptyset\}$} is the set of elements that are intersected by the interface, 
 ${\mathcal D}_{\mathcal{T}}(\mu) := \{K\in \mathcal{T}_h(\mu): K\cap\mathcal{D}(\mu)\} \cup {G}_h(\mu)$,
 background domain { {is typified by}  $\mathcal{B}$,  {while} its corresponding mesh is denoted by $\mathcal{B}_h$, such that ${\mathcal D}_{\mathcal{T}}(\mu) \subset \mathcal{B}$ and $\mathcal{T}_h(\mu) \subset \mathcal{B}_h$ for all $\mu \in \mathcal{K}$.
 See also Figure \ref{SurrogateMesh} (ii).
We remark that 
 $\mathcal{T}_h(\mu)$, ${G}_h(\mu)$ and ${\mathcal D}_{\mathcal{T}}(\mu)$ depend on $\mu$ through $\mathcal{D}(\mu)$ or its boundary, 
 while 
 the background domain $\mathcal{B}$ and its mesh $\mathcal{B}_h$ \emph{do not} depend on $\mu$.
Furthermore, the 
 the set of element faces $\mathcal{F}_G(\mu)$ associated with ${G}_h(\mu)$, is defined as follows: 
for each face {$F\in \mathcal{F}_G(\mu)$}, there exist two simplices $K \neq K'$ such that {$F=K\cap K'$ and at least one of the two is a member of ${G}_h(\mu)$}.
Note that
the boundary faces of $\mathcal{T}_h(\mu)$ are excluded from $\mathcal{F}_G(\mu)$.
On a face $F \in \mathcal{F}_G(\mu)$, $F=K\cap K'$, the jump of the gradient of $v \in C^0(\overline{\mathcal{D}}_{\mathcal{T}})$ is defined by $\ldb{{\bf n}_F}\cdot\nabla  v\rdb={\bf n}_F\cdot\nabla v |_K -{\bf n}_F\cdot \nabla v|_{K'}$, where 
${\bf n}_F$ denotes the outward pointing unit normal vector  {with respect to $K$}. 
So the CutFEM discretization is: 
we seek a discrete solution $u_h(\mu)$ in the finite element space
\begin{eqnarray}
\nonumber
V_h(\mu)=\left \{ \upsilon \in C^0({\overline {\mathcal D}}_{\mathcal T}(\mu))\,:\,\upsilon |_K  \in P^1(K), \, \forall K\in \mathcal{T}_h(\mu)    \right \},
\nonumber
\end{eqnarray}
such that $\forall \upsilon_h(\mu)\in V_h(\mu)$
\begin{eqnarray*}\label{CutFem_Discrete_Form_omega}
\alpha^h(u_h(\mu),v_h(\mu))&=&\ell^h(v_h(\mu))
\nonumber \\
\alpha^h(u_h(\mu),v_h(\mu))&=&{\left( \nabla u_h(\mu) ,\nabla v_h(\mu)\right)_{{\mathcal D}(\mu)}}{-\left({\bf n}_\Gamma\cdot \nabla u_h(\mu),\upsilon_h(\mu)\right)_{\Gamma
(\mu)}}
  \nonumber \\
 && {- \left(u_h(\mu),{\bf n}_\Gamma \cdot \nabla\upsilon_h(\mu)\right)_{\Gamma
(\mu)}}
+\left(\gamma_Dh^{-1}u_h(\mu),v_h(\mu)\right )_{\Gamma
(\mu)}
  \nonumber \\
&&   
{
   +\left(\gamma_N h {\bf n}_\Gamma\cdot \nabla u_h(\mu),{\bf n}_\Gamma\cdot \nabla v_h(\mu) \right)_{\Gamma _N(\mu)}}
   + { {j(u_h(\mu),v_h(\mu))}  }
   \nonumber \\
  \ell^h(v_h(\mu))&=&{
   \left(g
  {(\mu)}, \upsilon_h(\mu)\right )_{{\mathcal D}{(\mu)}}}
{- \left(g_D,{\bf}n_\Gamma \cdot \nabla\upsilon_h(\mu)\right)_{\Gamma
(\mu)}}
  \\&&
  \hskip-20pt
    {+\left(\gamma_Dh^{-1}g_D,v_h(\mu)\right )_{\Gamma
(\mu)}}+\left({g_N(\mu),v_h(\mu)}
  {+\gamma_N h {\bf n_\Gamma}\cdot \nabla v_h(\mu)}\right)_{{\Gamma _N{(\mu)}}},\nonumber
\end{eqnarray*}
 where the stabilization term
 \begin{eqnarray*}
  {j(u_h(\mu),v_h(\mu))=\sum_{F\in \mathcal{F}_G}\left(\gamma_1 h\ldb{{\bf n_F}}\cdot\nabla  u_h(\mu)\rdb,\ldb{{{\bf n_F}}\cdot\nabla} v_h(\mu)\rdb\right)_F,}
  \end{eqnarray*}
 extends  the coercivity from the physical domain ${\mathcal D}(\mu)$ to the $\mu$- {dependent}  mesh domain ${\mathcal D} _{\mathcal{T}}$, $\gamma_D$, $\gamma_N$, and $\gamma_1$ are positive penalty parameters.
The coefficients $\gamma_D$ and $\gamma_N$ account for a Nitsche weak imposition of boundary conditions.
We  {set up} some experiments again for  embedded finite element methods and ROMs emphasizing on improvements in cases with
large geometrical, deformations,
The strongly nonlinear parametrized domain $\mathcal{D}(\mu) \subset \mathbb{R}^2$ is an ellipse,  
 described by the level set:
$\phi(x,y;\mu_1, \mu_2, \mu_3, \mu_4) = \mu^2_2(x-\mu_3)^2 +\mu^2_1(y-\mu_4)^2 -\mu^2_1 \mu^2_2 R$,
where the 
reference radius $R = 0.05$, 
the length of the axes of the ellipse: $(\mu_1, \mu_2) \in [0.3, 1.8]^2$, while 
 the position of the center of the ellipse: $(\mu_3, \mu_4) \in [-0.85, 0.85]^2$. 
A corresponding 
 background domain $\mathcal{B} = [-1.2, 1.2]^2$,
is chosen so that
 the ellipse is strictly contained in $\mathcal{B} $
 for any 
$\mu = (\mu_1, \mu_2, \mu_3, \mu_4)$ in the parametric range $\mathcal{K} = [0.3, 1.8]^2 \times [-0.85, 0.85]^2$.
 The value $\overline{\mu} = (1, 1, 0, 0)$, corresponding to a circle of radius $R$ centered in the origin, is chosen for what concerns the transport method. 
The data of the problem  {described by Equation} \ref{P:Poisson} are the 
force $g(x, y; \mu) = 20$  and the Dirichlet boundary force $g_D(x, y; \mu) = 0.5 +x y$.

\begin{figure} \centering 
{(i) Standard POD modes (no preprocessing)}\\
\includegraphics[width=0.22\textwidth]{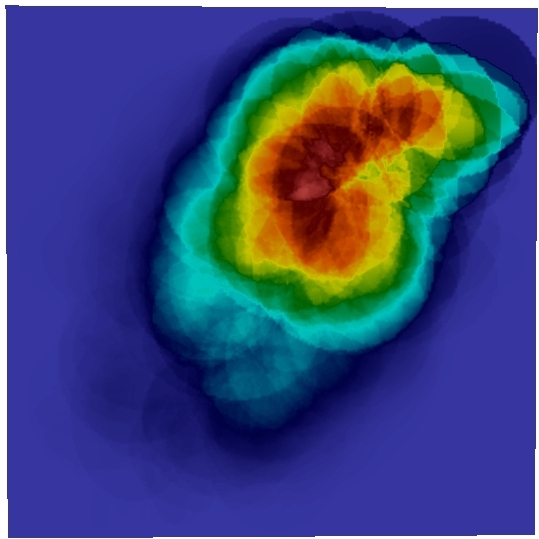} 
\includegraphics[width=0.22\textwidth]{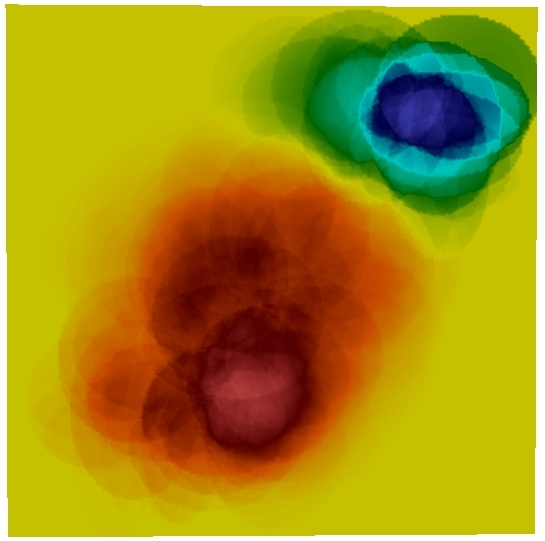}
\includegraphics[width=0.22\textwidth]{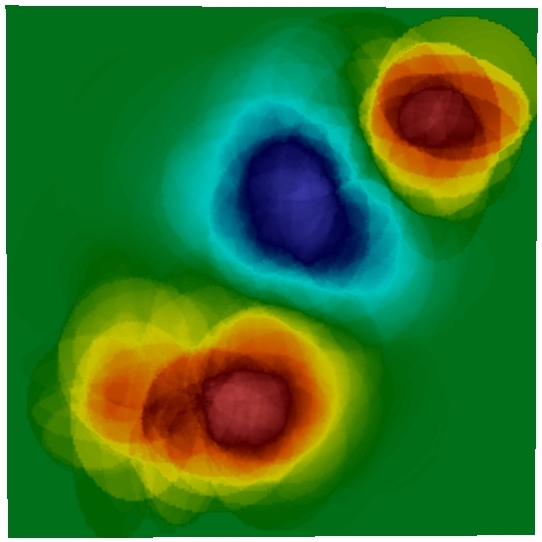}
 \\
\includegraphics[width=0.22\textwidth]{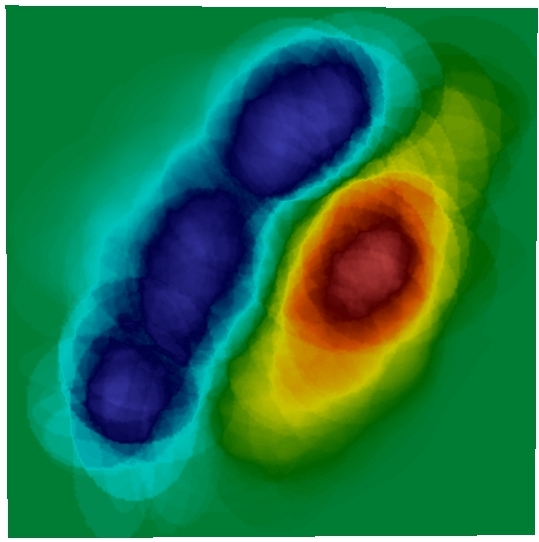}
\includegraphics[width=0.22\textwidth]{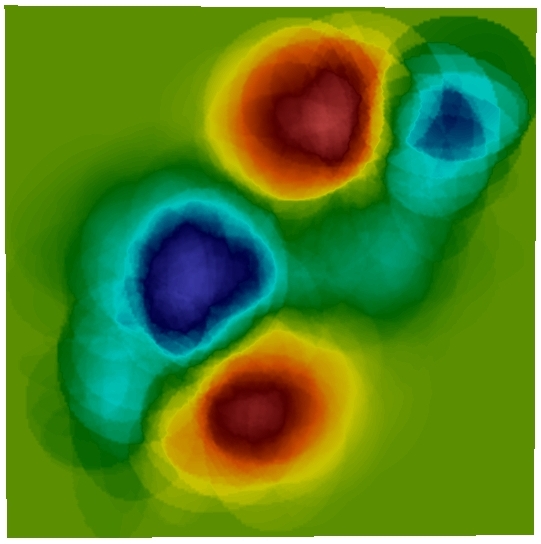}
\includegraphics[width=0.22\textwidth]{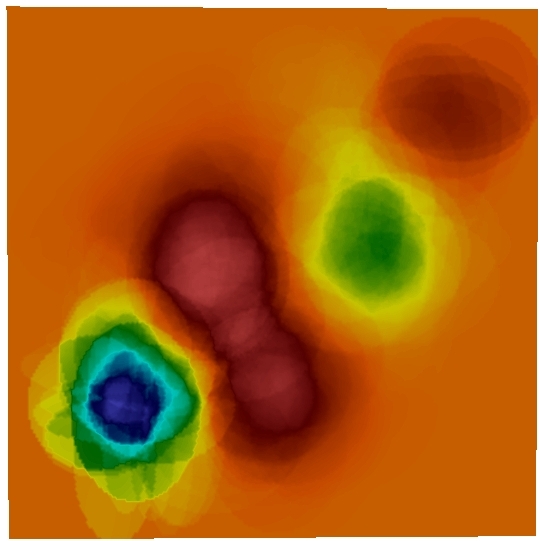}
\\
{(ii)POD modes with preprocessing}
\\
\includegraphics[width=0.22\textwidth]{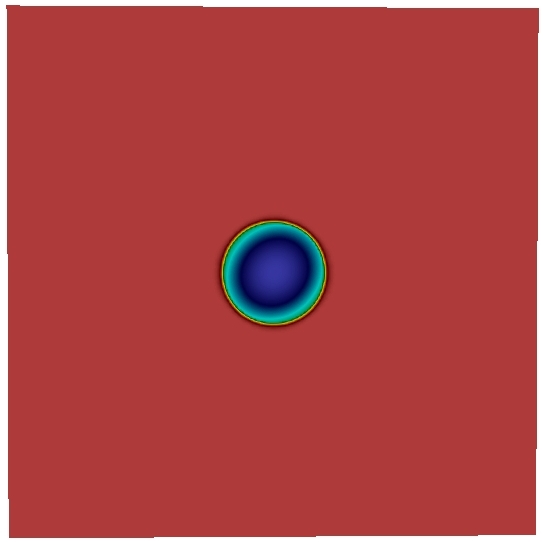}
\includegraphics[width=0.22\textwidth]{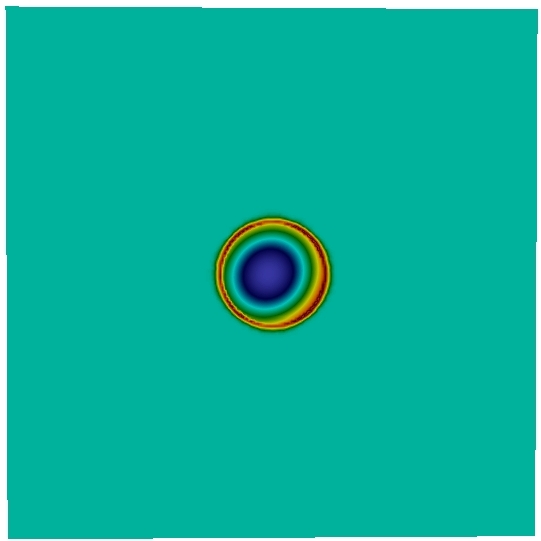}
\includegraphics[width=0.22\textwidth]{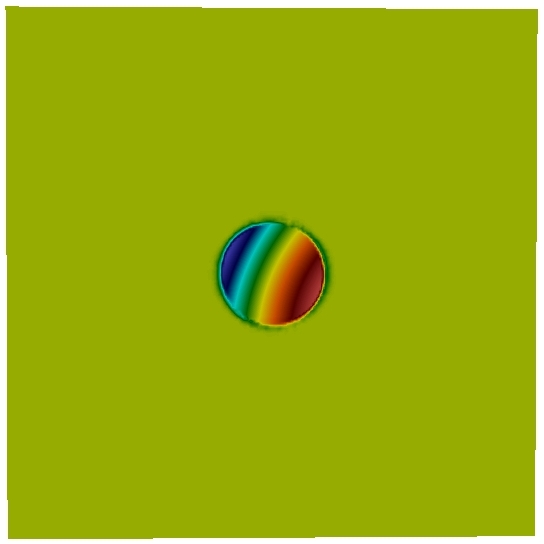}
\\
\includegraphics[width=0.22\textwidth]{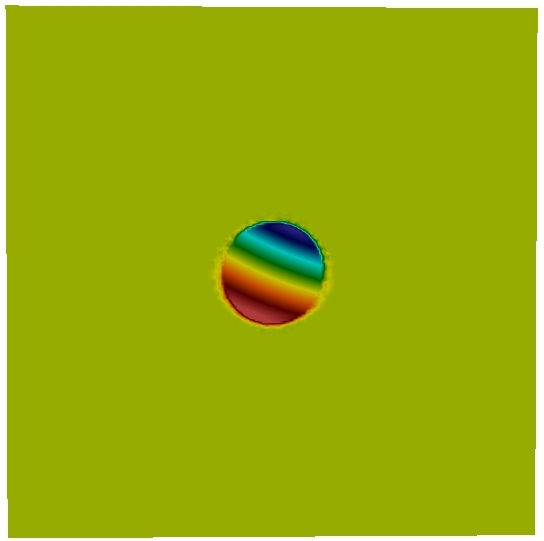}
\includegraphics[width=0.22\textwidth]{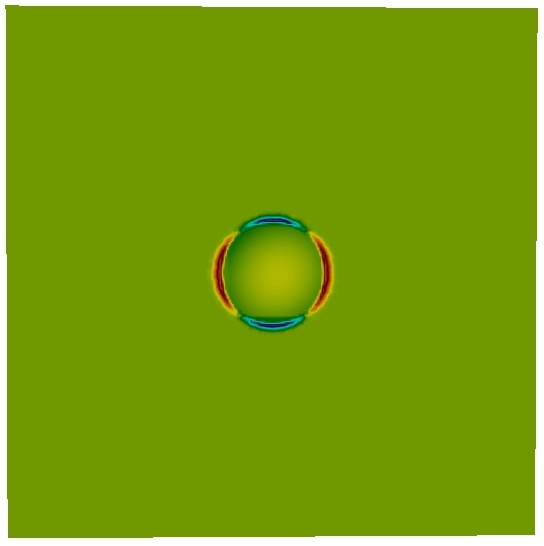}
\includegraphics[width=0.22\textwidth]{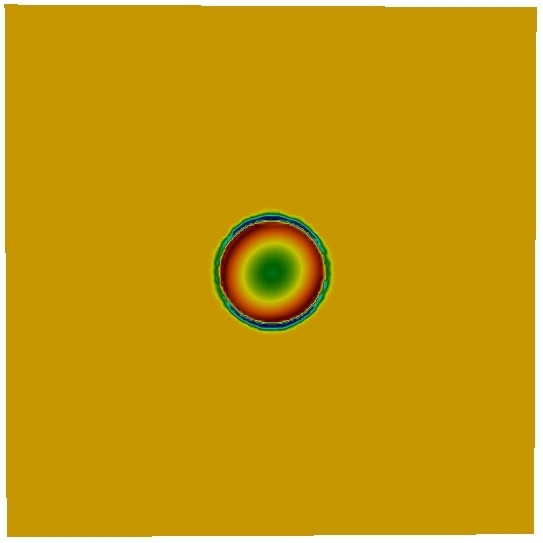}
  \caption{Six classical (i) and six improved (ii) POD modes.}
  \label{fig:yes_no_preprocessing}
\end{figure}
\begin{figure} \centering \vskip-40pt
\begin{minipage}{\textwidth}\centering
(i)
\begin{minipage}{0.65\textwidth}
   \includegraphics[width=\textwidth]{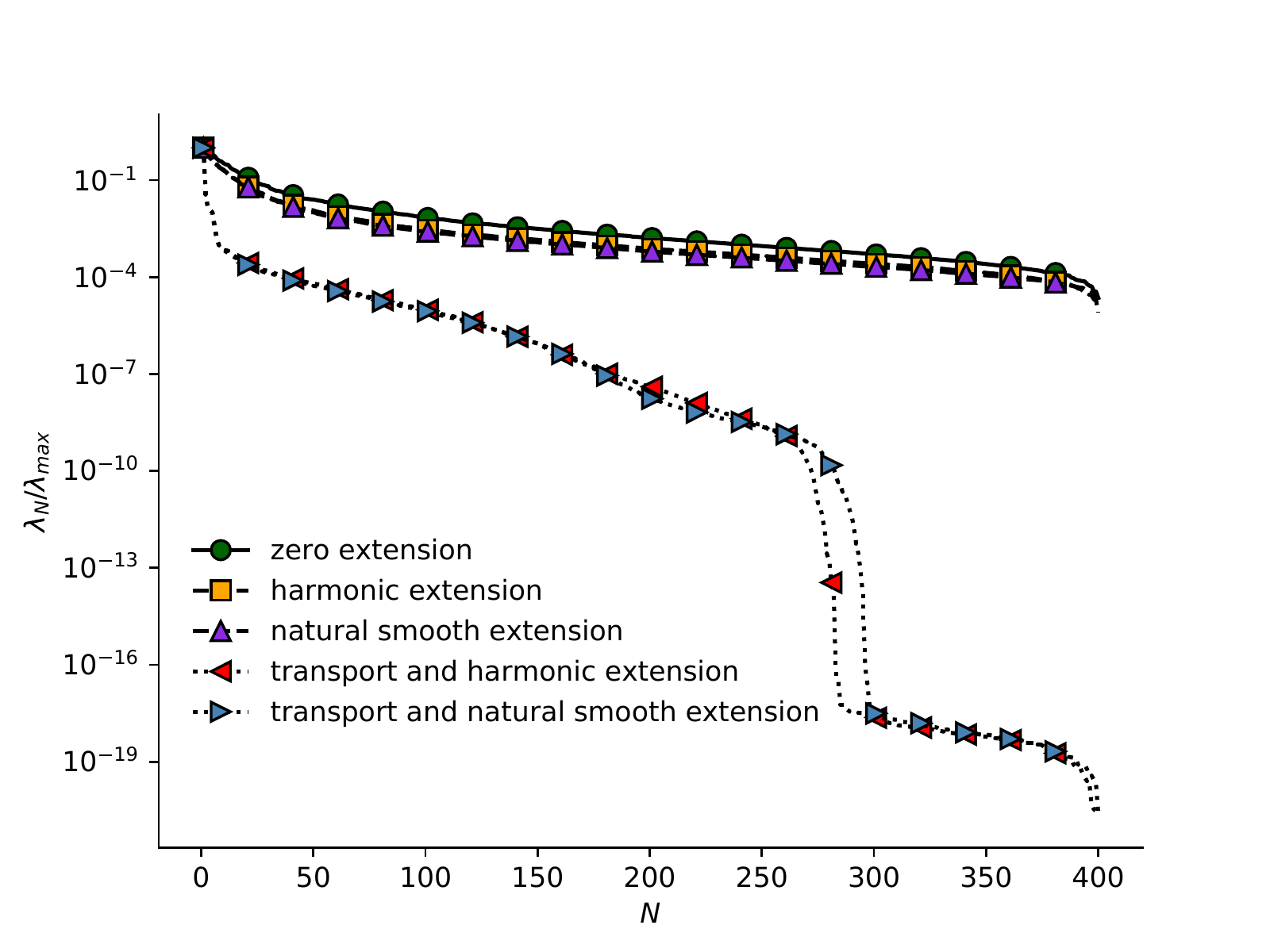}
\end{minipage}
\begin{minipage}{0.65\textwidth}
  \hskip7pt \includegraphics[width=\textwidth]{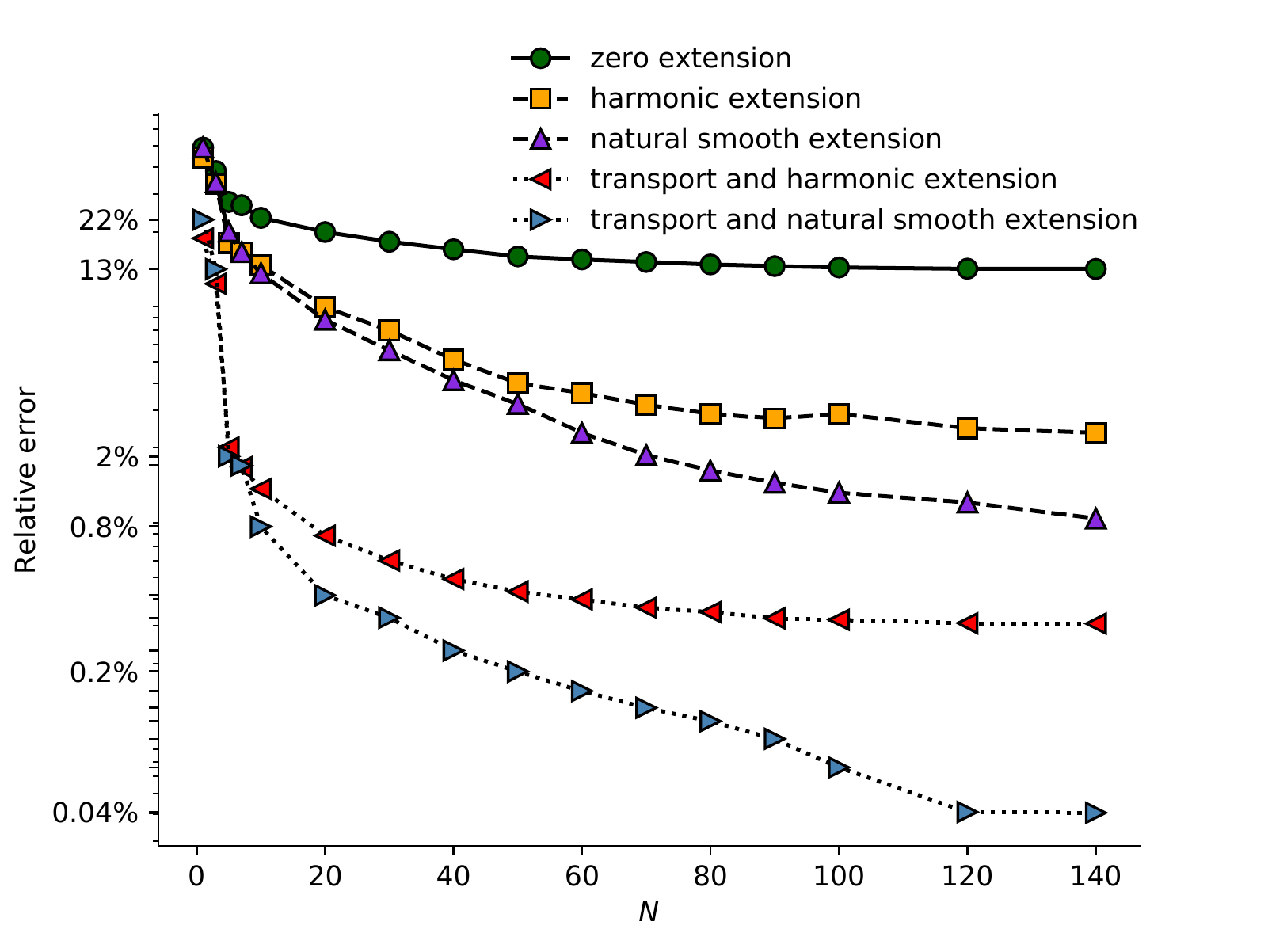}
\end{minipage}
 \end{minipage}
\vskip-20pt
(ii)
\begin{minipage}{0.65\textwidth}\centering
\includegraphics[width=\textwidth]{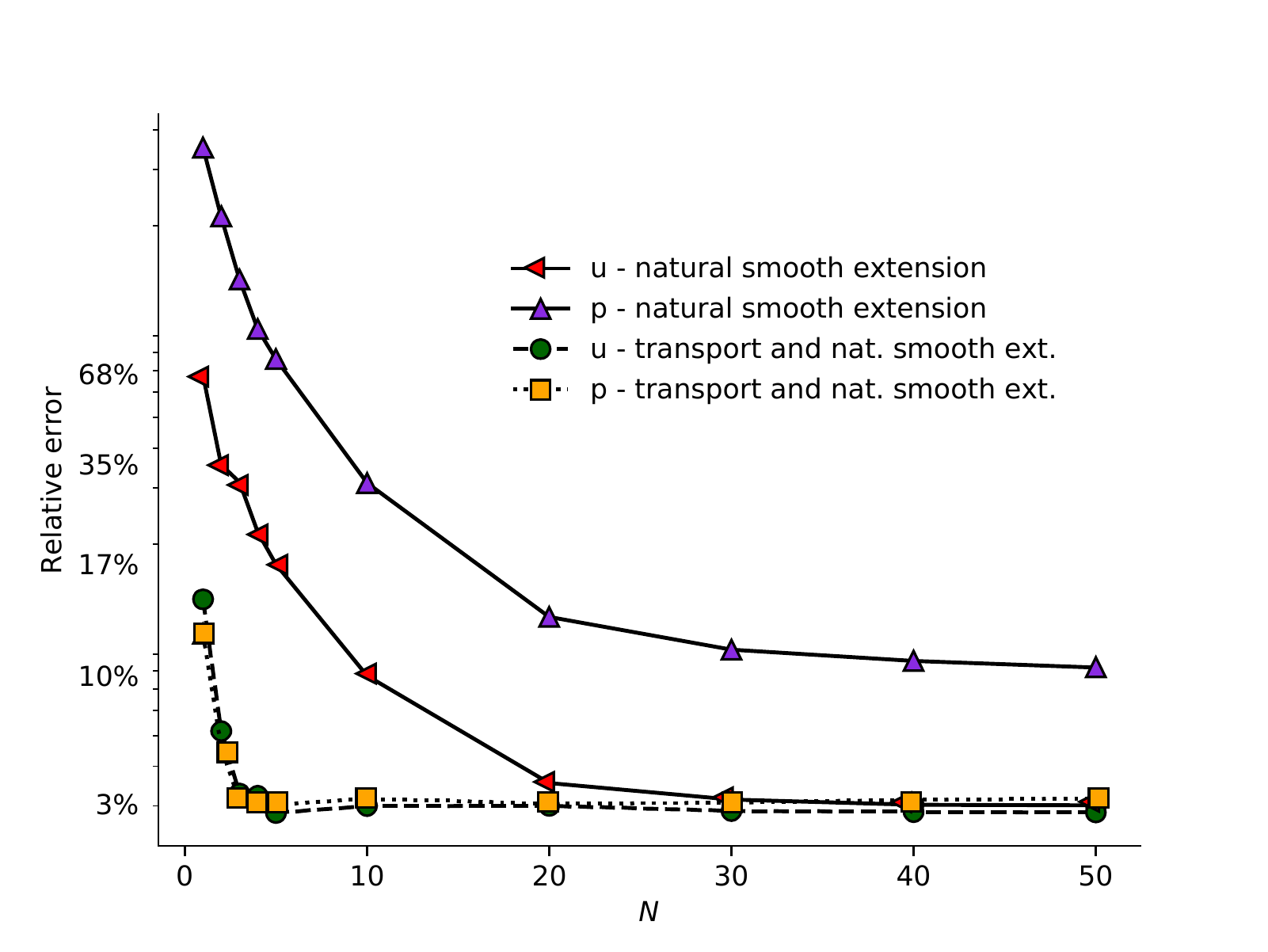}
  \end{minipage}
  \caption{
  (i) Poisson system (CutFEM): Eigenvalues decay and error analysis between reduced order and high fidelity approximations with and without transport, (ii) Steady Stokes (CutFEM): Relative errors with and without transport.
  }
  \label{fig:Stokes_transport}
  \label{fig:Poisson results}
\end{figure}
Explaining the results in Figures \ref{fig:Poisson results} (i) for the CutFEM Poisson system, the  reduced solution obtained from the  zero extension, is inaccurate even for $N=140$ modes, being affected by relative errors of the order of $10^{-1}$. A non-zero extension is beneficial, resulting in  relative errors of the order of $10^{-2}$, for the maximum value of $N$. 
The combination with inverse transportation allows to further improve results, up to errors of $10^{-4}$ for $N=140$ modes in the case of POD basis obtained from transport and natural smooth extension.
Thus, the pivotal role of snapshots transportation can be inferred from these results, being capable of improving the results of almost three orders of magnitude compared to the simplest zero extension.
Nonetheless, 
all methods reach a plateau after which no further improvement is shown. We claim that this is due to integration errors occurring on $\Gamma(\mu)$ and Nitsche weak imposition of Dirichlet boundary conditions, e.g. the maximum values of the error are consistently attained on the boundary. 

Further investigation have shown similar good behaviour for fluid flow systems, namely Steady Stokes and RB with transportation and CutFEM, 
 \cite{KaBaRO18},  with  relative errors   without and with transport improvements to be visualized in Figure \ref{fig:Stokes_transport} (ii). Once more, the transported snapshots and the inverse transported modes appear beneficial for both velocity and pressure with a need of very few numbers of modes and in particular it is managed to reach best accuracy with only three basis components.
\section{ROM and a fourth order evolutionary non-linear system}
In this last section, we consider the Cahn-Hilliard model problem 
describing the phase flow time evolution. An unknown function $u$ indicates the perturbation of the concentration of one of the phases of e.g. fluid components constituting a liquid mixture which contains a binary fluid. 
As first suggested by \cite{DeGroot62} and thereafter extended in \cite{CaHi58}, and if we assume that the mobility is equal to $1$ and $\varepsilon$ is a measure of the size of the interface of two fluids, then the mass flux is given by
${{\bf{J}}(\mu)}=-\nabla \left( \frac{{1}}{\varepsilon^2}F(u(\mu)) - {{\varepsilon^2}} \Delta u(\mu)\right)$,
where 
$F$ denotes the chemical potential difference between the two species. 
From the latter we can derive that the
Ginzburg--Landau energy is
$
E(u(\mu))=\int_{{\mathcal{D}}}\left(
F(u(\mu))+ \frac{{\varepsilon^2}}{2}|\nabla u(\mu)|^2\right)d{\bf{x}},
$ 
and that the equilibrium state of the considered mixture minimizes the above Ginzburg--Landau energy,
subject to the mass conservation: 
$\frac{\partial}{\partial t}  \int_{{\mathcal{D}}(\mu)}u(\mu)d{\bf{x}} = 0$.
Hence, the parametrized
Cahn--Hilliard system, 
can be described as:
\begin{eqnarray}\label{Eq:4thOrder}
\frac{\partial{u}(\mu)}{\partial{t}} = -\varepsilon^2 \Delta ^2 u(\mu) + \frac{1}{\varepsilon^2}\Delta F'(u(\mu)) ,  &&\text{ in } {{\mathcal{D}}}(\mu)\times[0,T],
\\\label{Eq:4thOrder_boundary}
\partial_n u(\mu) = \partial _n (-\varepsilon ^2 \Delta u(\mu) + \frac{1}{\epsilon^2}F'(u(\mu))) = g_N(\mu), 
&&\text{ on } {\Gamma}(\mu)\times[0,T],
\\ \label{Eq:4thOrder_initial}
u(\cdot, 0)= u_0 (\cdot),    &&\text{ in } {\mathcal{D}}(\mu),
\end{eqnarray}
where $n$ is the unit outer normal vector of $\Gamma$, 
 $F$ is a double well
function of $u$ usually polynomial of fourth power:  
\begin{equation}\label{double-well}
 F(u(\mu)) = \gamma_2 \frac{u^4(\mu)}{4} + \gamma_1 \frac{u^3(\mu)}{3} +\gamma_0\frac{u^2(\mu)}{2} \text{ with } \gamma_2 > 0.
\end{equation}
Starting from the above initial form and based on a splitting method and Nitsche boundary enforcement and with efficient CutFEM stabilization, we manage a proper weak form of $H^1$ space regularity, for more details we refer to \cite{KaratzasRozzaCH20}.
The CutFEM mesh and level set geometry can be seen in Figure \ref{fig:circular_embedded_geometry},
\begin{figure} 
\centering
  \includegraphics[width=0.25\textwidth]{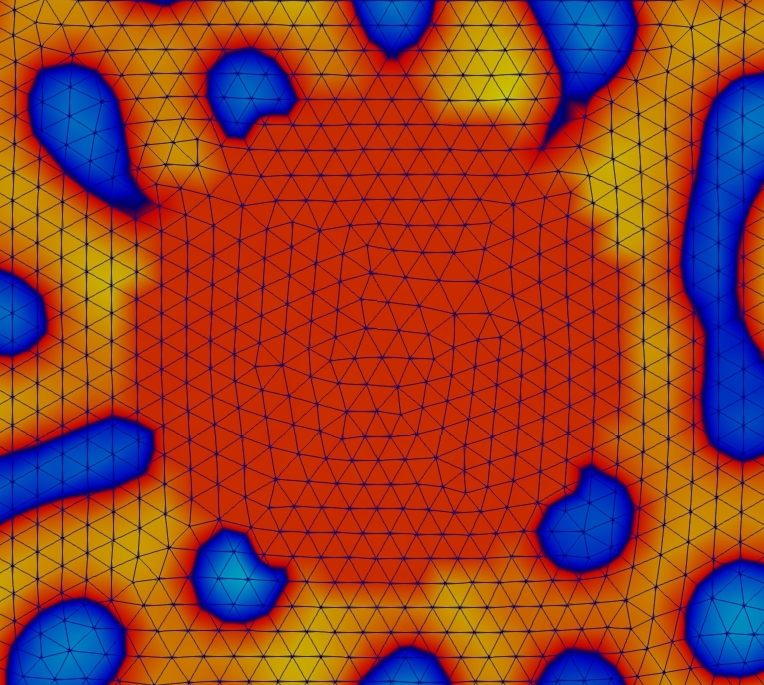}
  \caption{Circular embedded geometry (parametrized).} 
\label{fig:circular_embedded_geometry}
\end{figure} 
while some ROM basis compononents can be seen in Figure \ref{fig:CH_modes}, and the concentration field in the full-order, the reduced, absolute error level after parametrization of the embedded circle,  
is visualized for the time instances $t=[28,36,46,100]dt$ and for a randomly selected parameter $\mu=0.4261$ in Figure \ref{fig:CH_FULL_RED_ERROR_1P}. 
\begin{figure} 
\centering
\begin{minipage}{\textwidth}
\centering
\begin{minipage}{0.16\textwidth}
  \includegraphics[width=\textwidth]{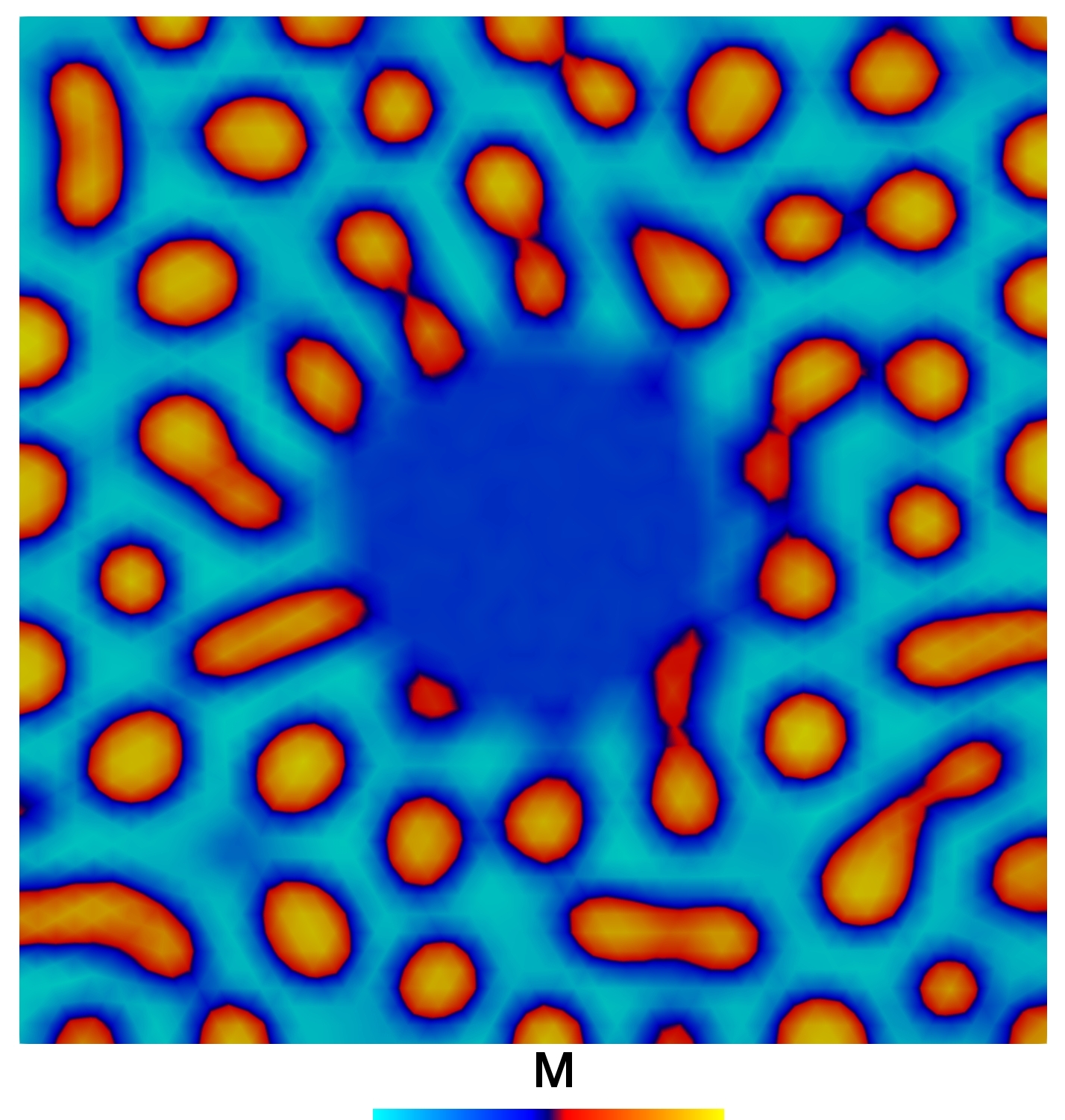}
\end{minipage}
\begin{minipage}{0.16\textwidth}
  \includegraphics[width=\textwidth]{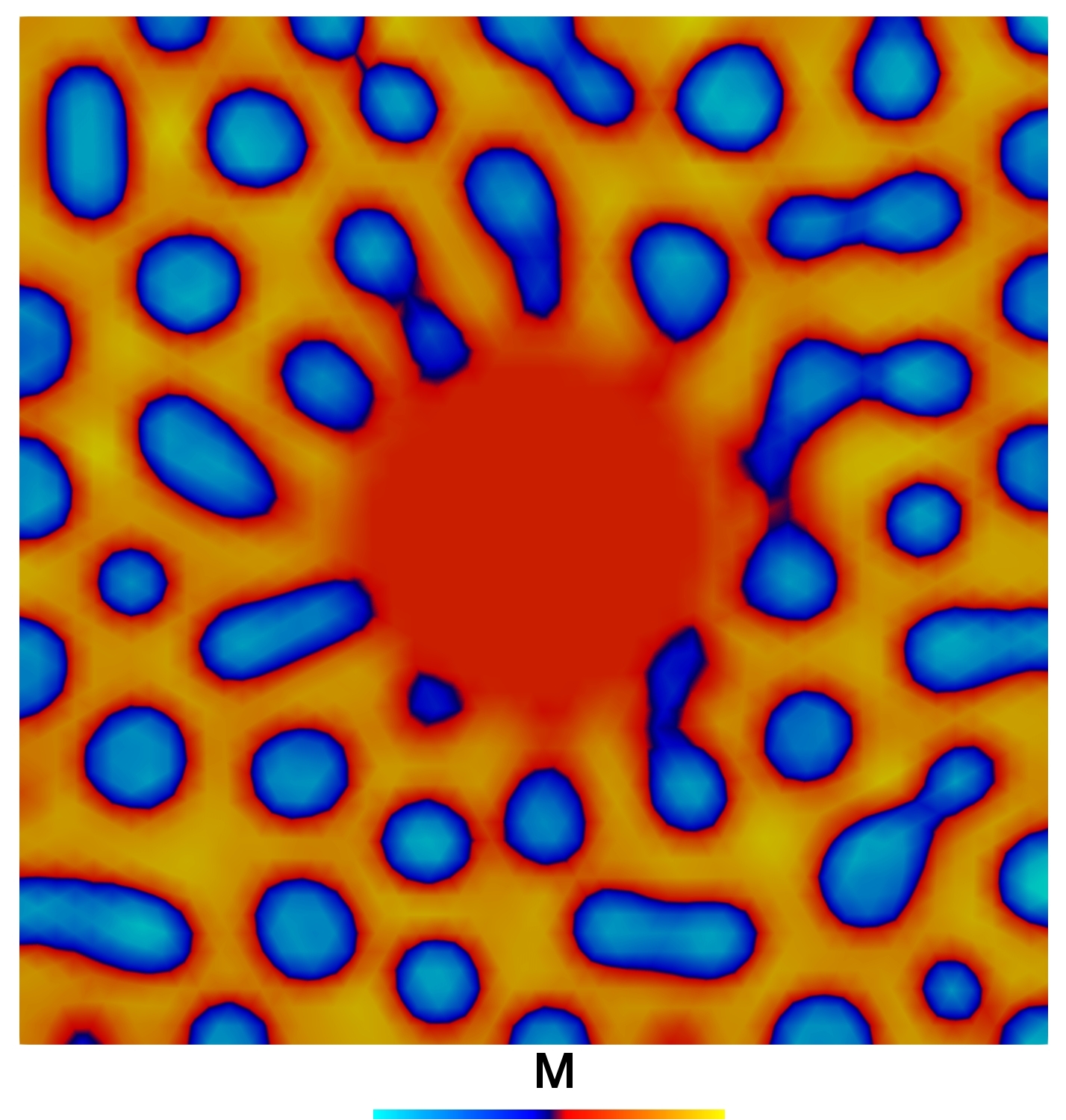}
\end{minipage}
\begin{minipage}{0.16\textwidth}
  \includegraphics[width=\textwidth]{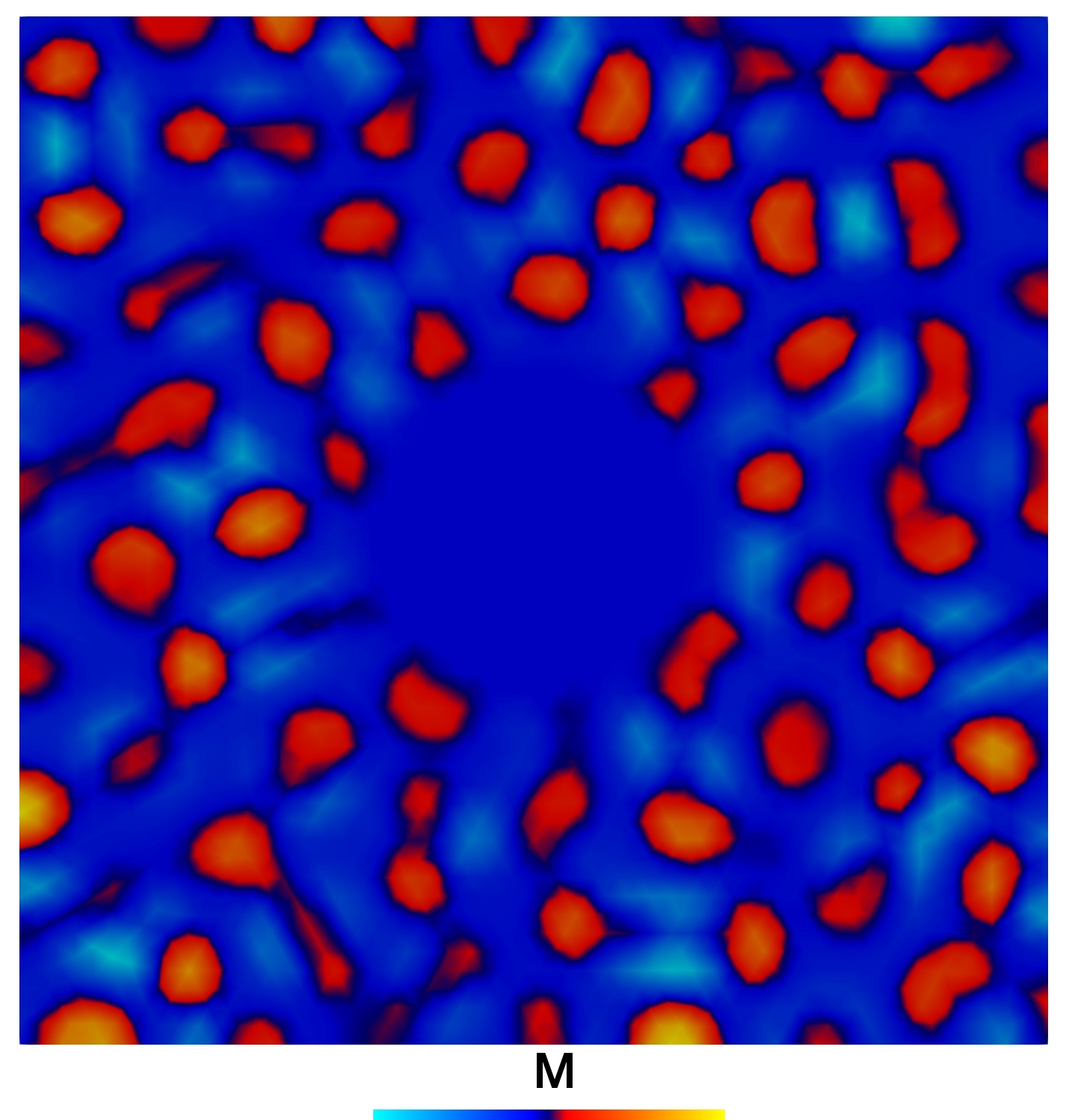}
\end{minipage}
\begin{minipage}{0.16\textwidth}
  \includegraphics[width=\textwidth]{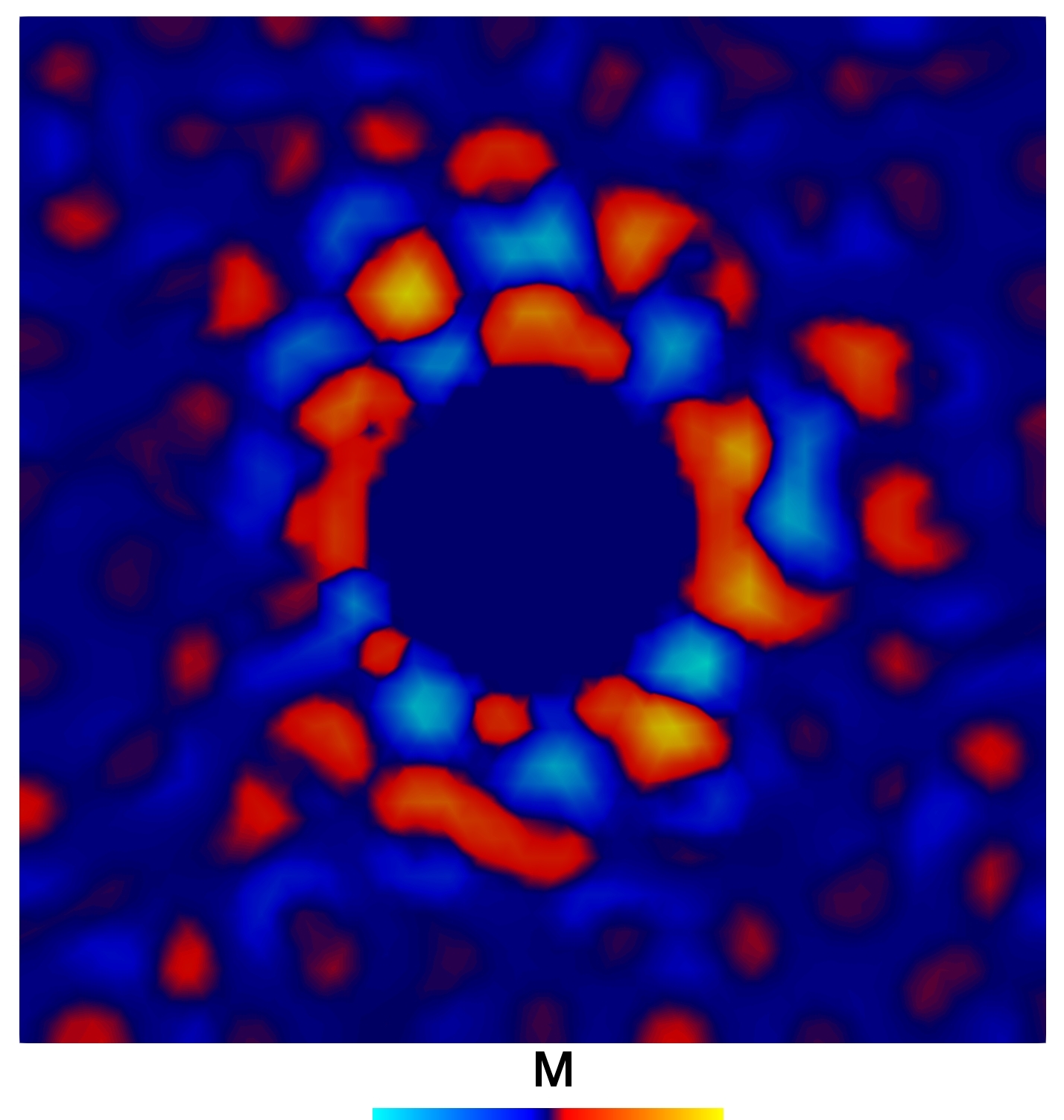}
\end{minipage}
\begin{minipage}{0.16\textwidth}
  \includegraphics[width=\textwidth]{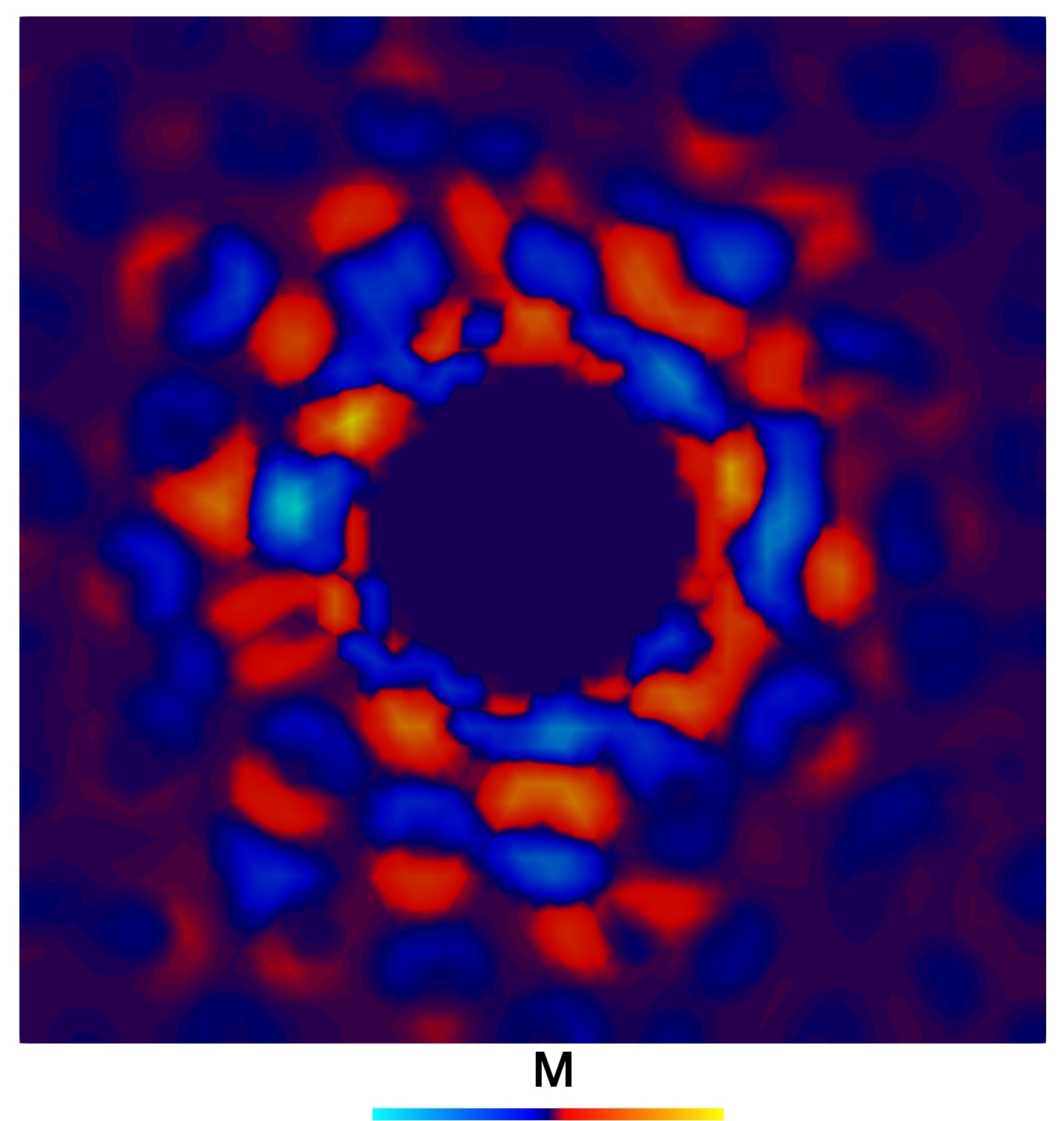}
\end{minipage}
\begin{minipage}{0.16\textwidth}
  \includegraphics[width=\textwidth]{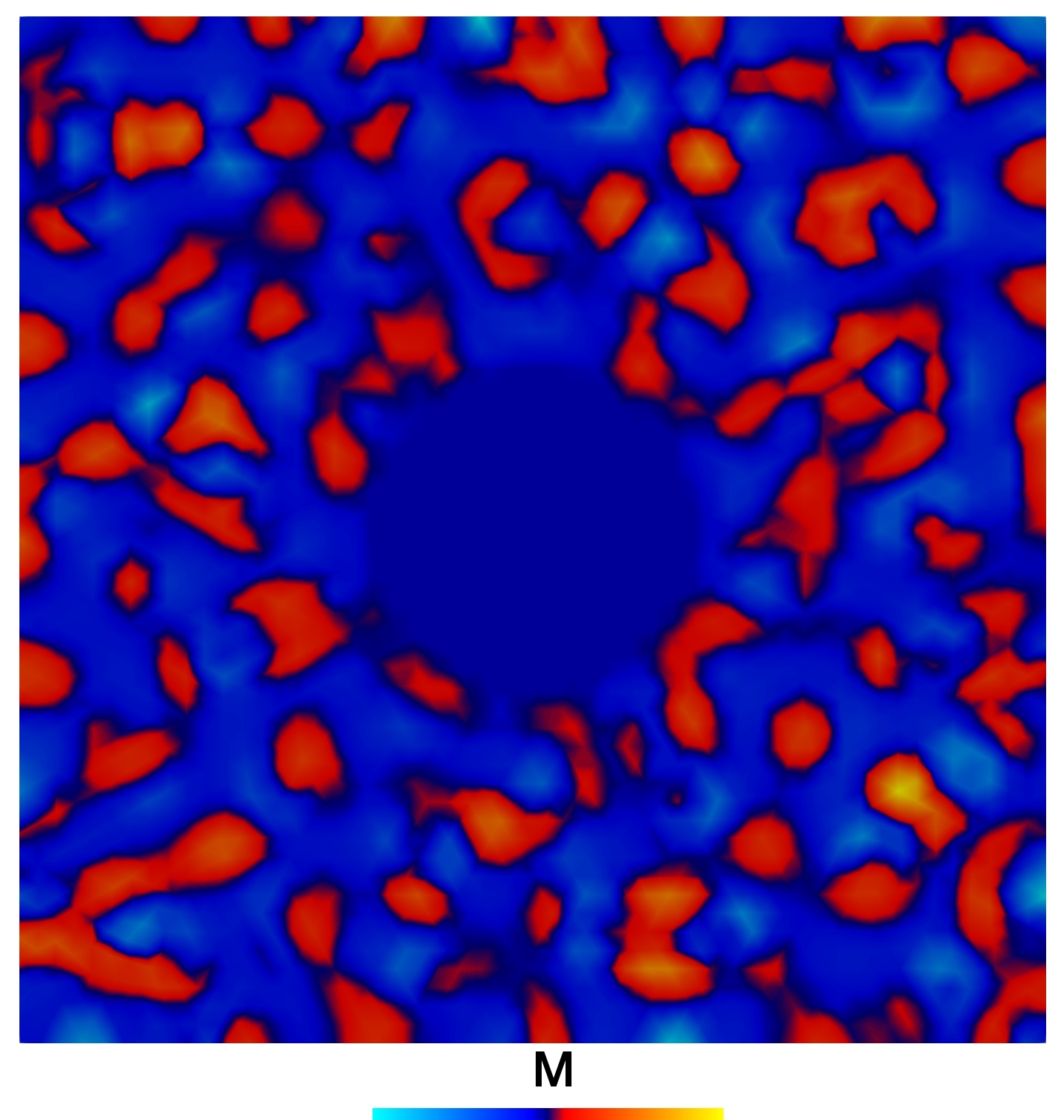}
\end{minipage}
\end{minipage}
\caption{Cahn-Hilliard (CutFEM): ROM basis results, the first six modes.
} 
\label{fig:CH_modes}
\end{figure} 
\begin{figure} 
\centering
\begin{minipage}{\textwidth}
\centering
\begin{minipage}{0.24\textwidth}
  \includegraphics[width=\textwidth]{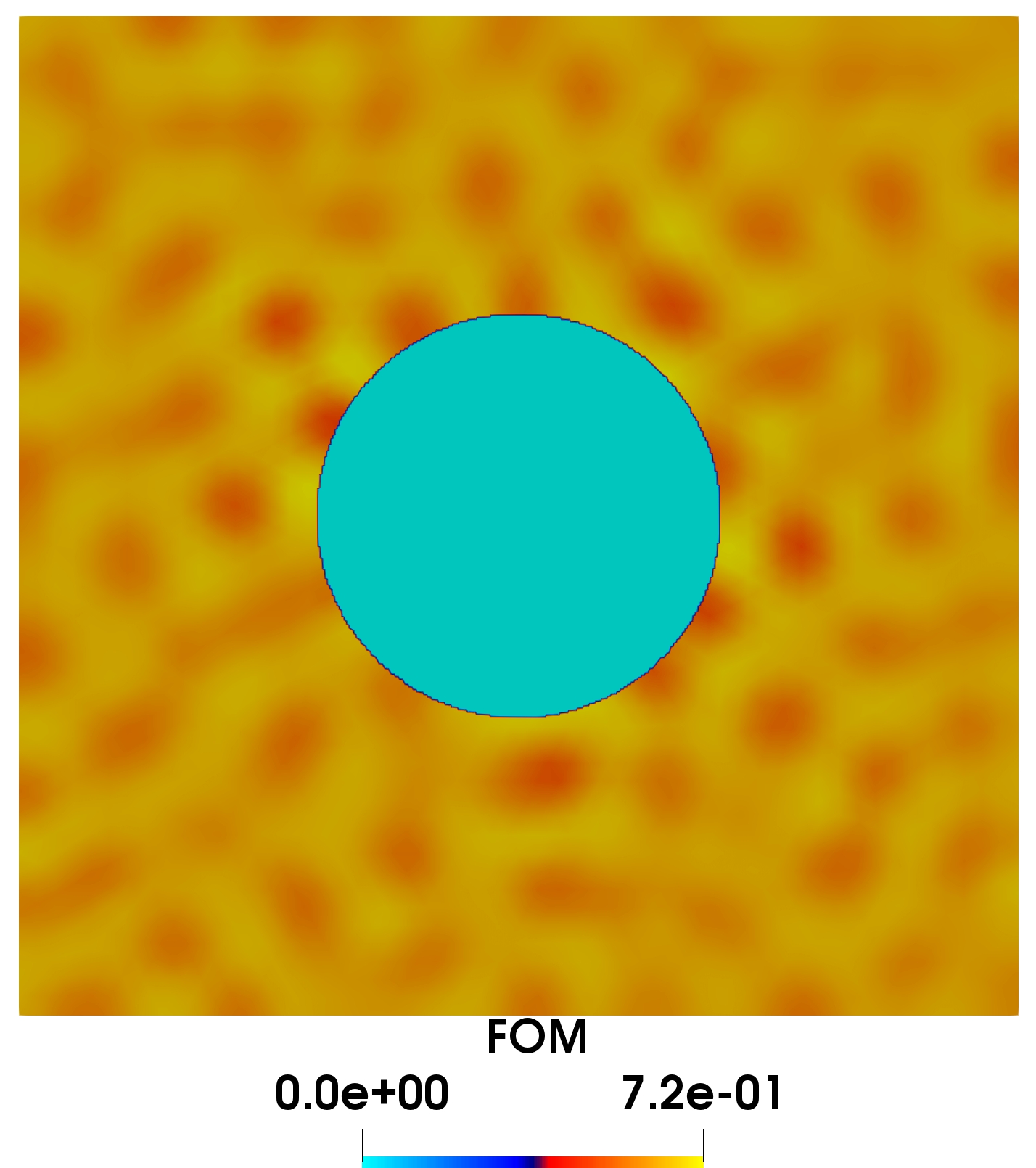}
\end{minipage}
\begin{minipage}{0.24\textwidth}
  \includegraphics[width=\textwidth]{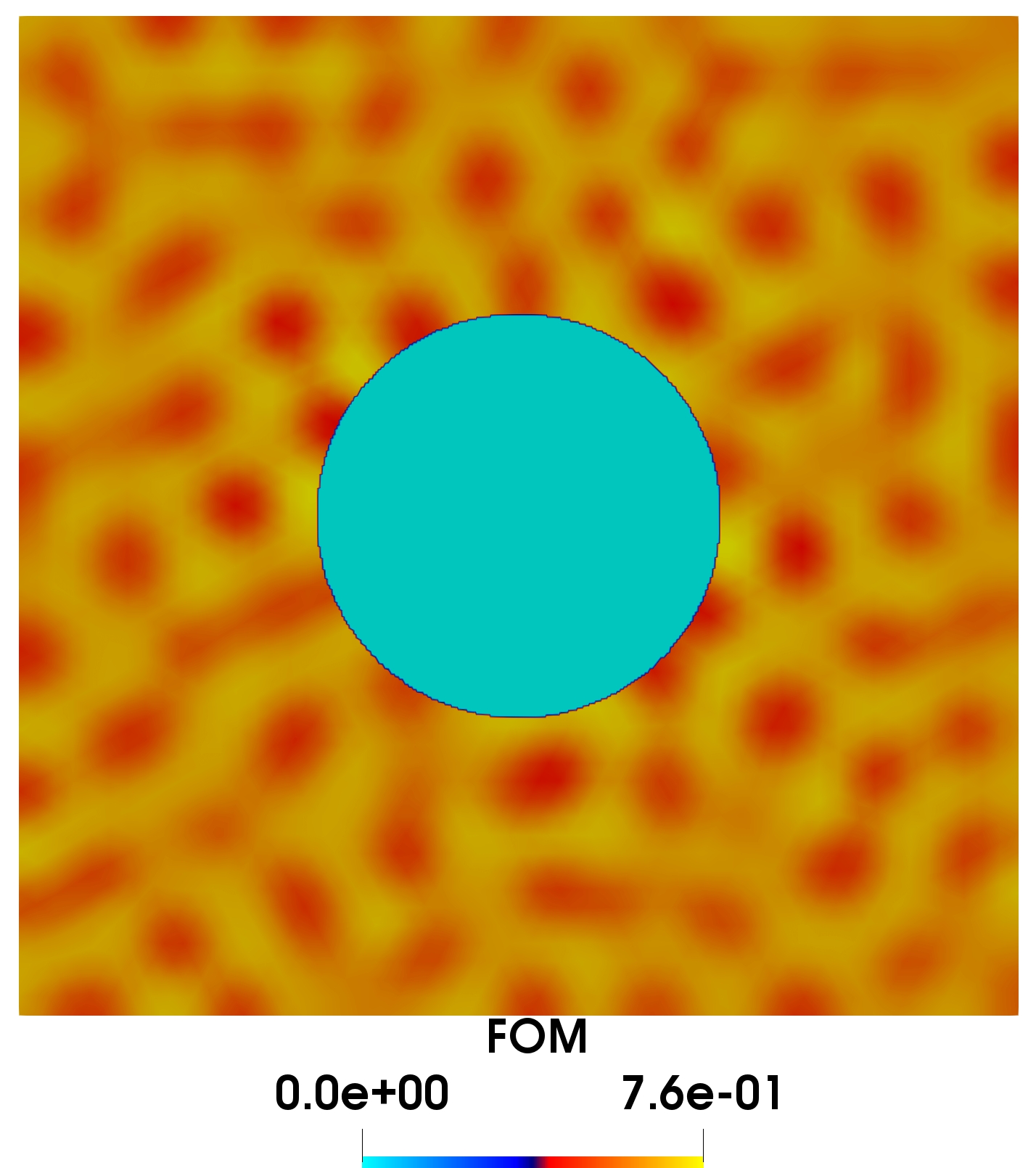}
\end{minipage}
\begin{minipage}{0.24\textwidth}
  \includegraphics[width=\textwidth]{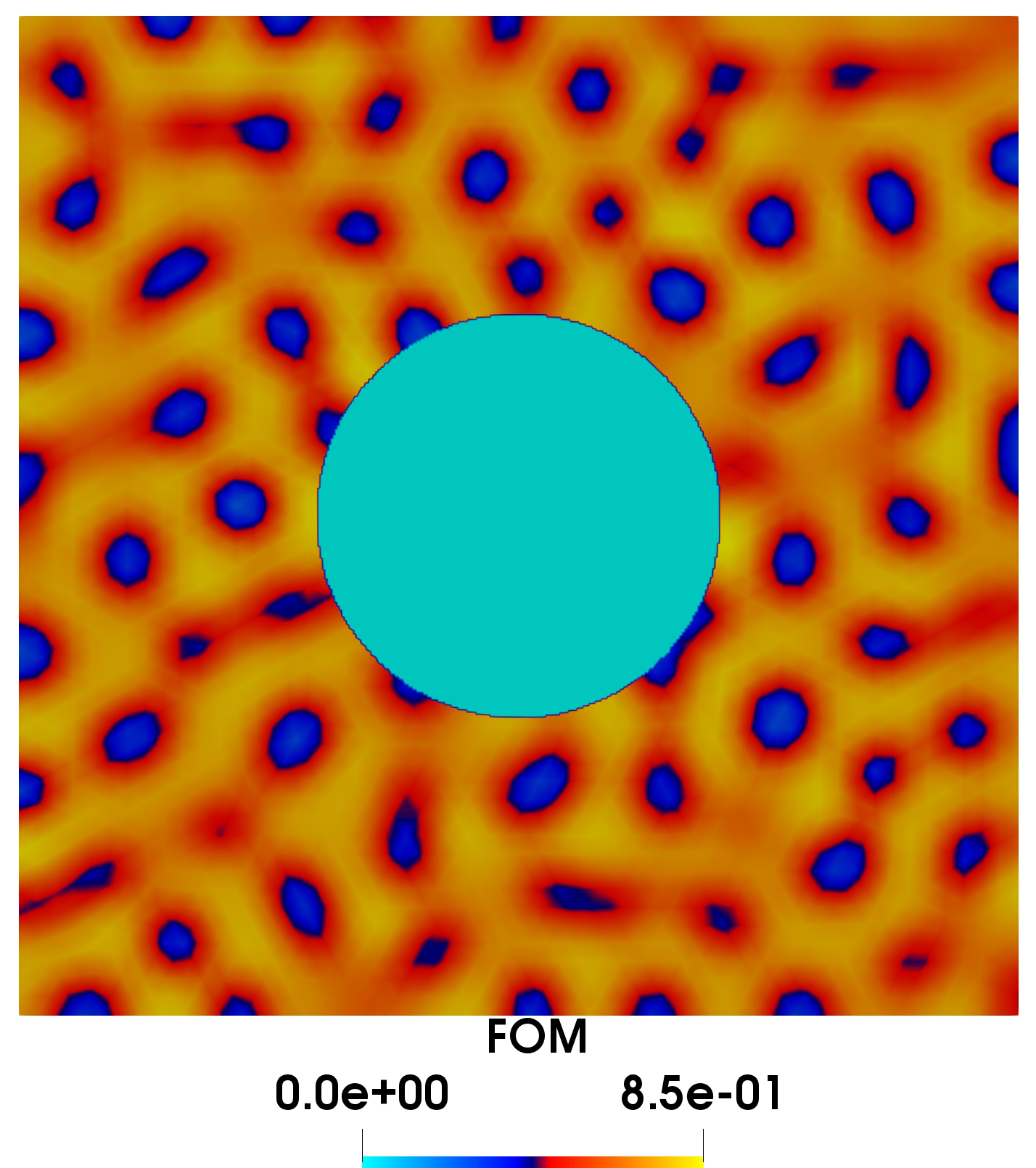}
\end{minipage}
\begin{minipage}{0.24\textwidth}
  \includegraphics[width=\textwidth]{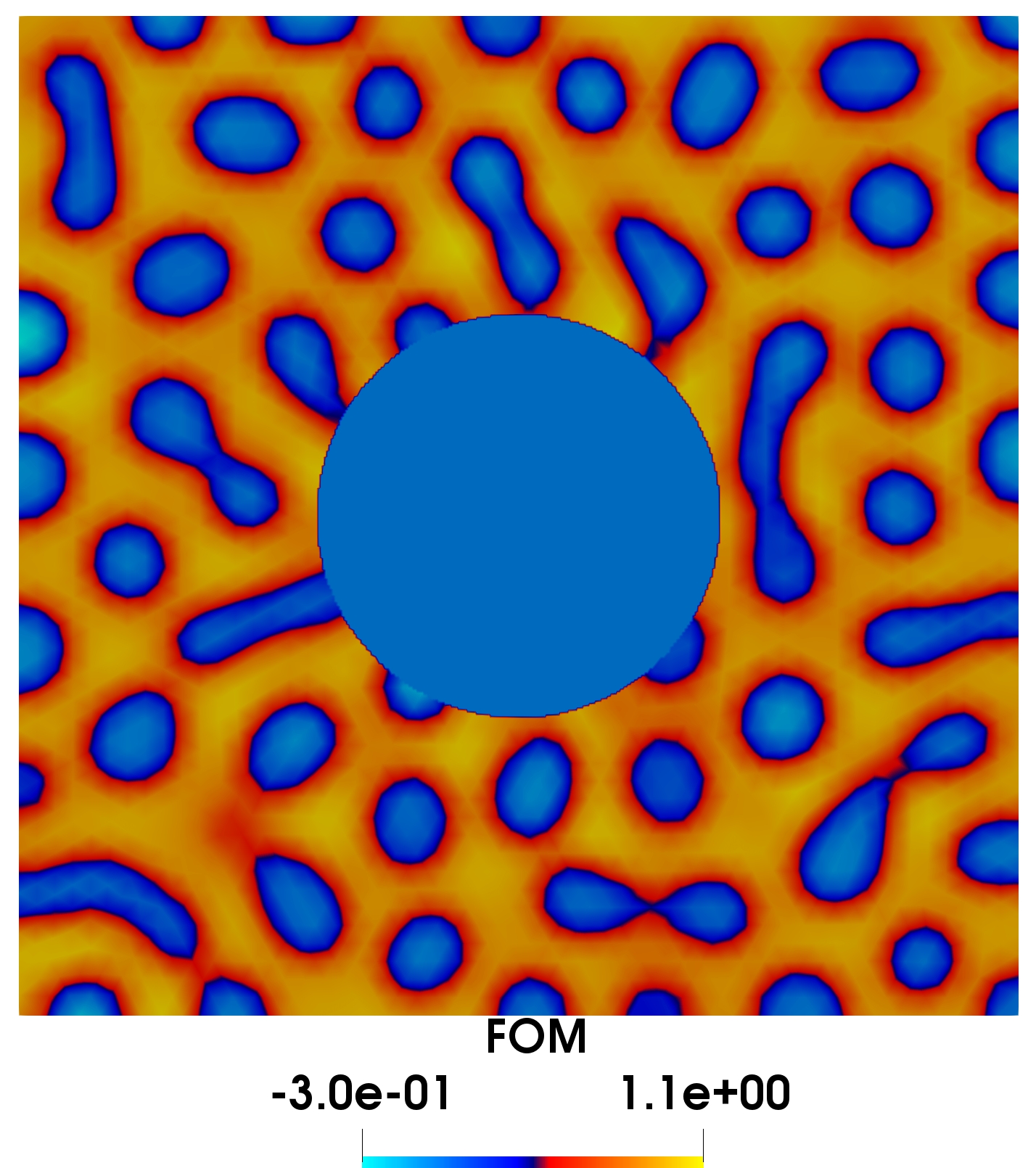}
\end{minipage}
\\
\begin{minipage}{0.24\textwidth}
  \includegraphics[width=\textwidth]{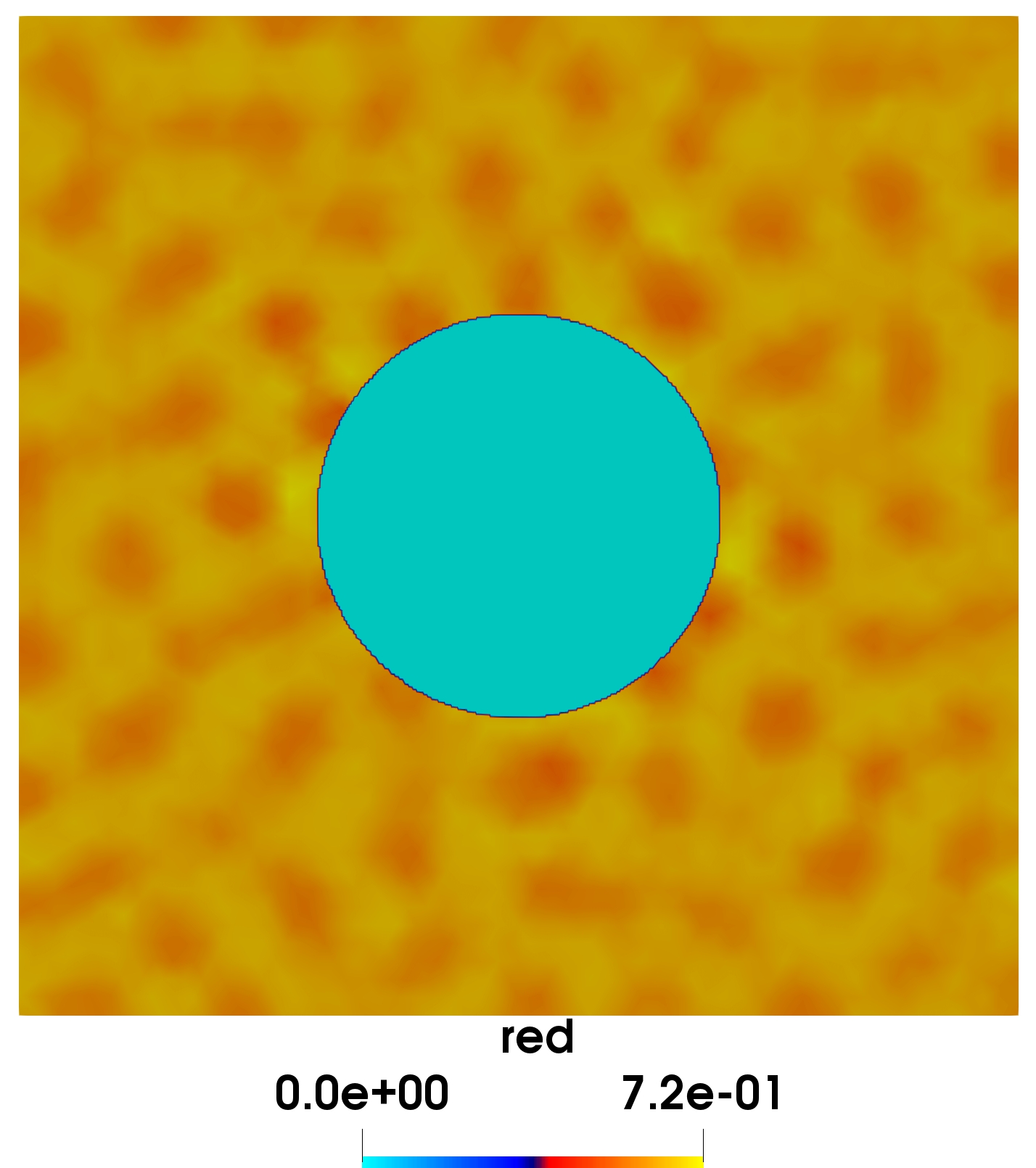} 
\end{minipage}
\begin{minipage}{0.24\textwidth}
  \includegraphics[width=\textwidth]{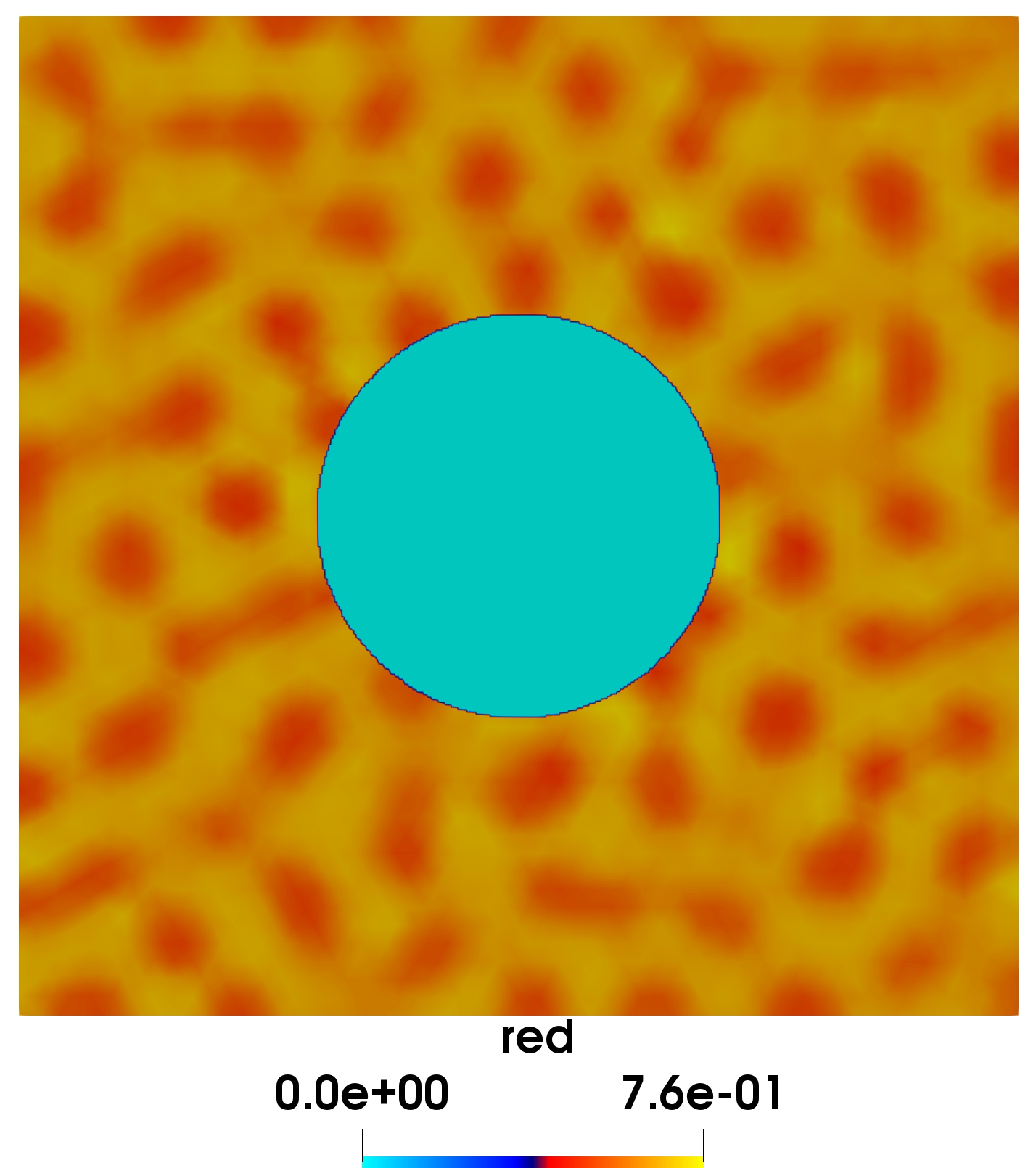} 
\end{minipage}
\begin{minipage}{0.24\textwidth}
  \includegraphics[width=\textwidth]{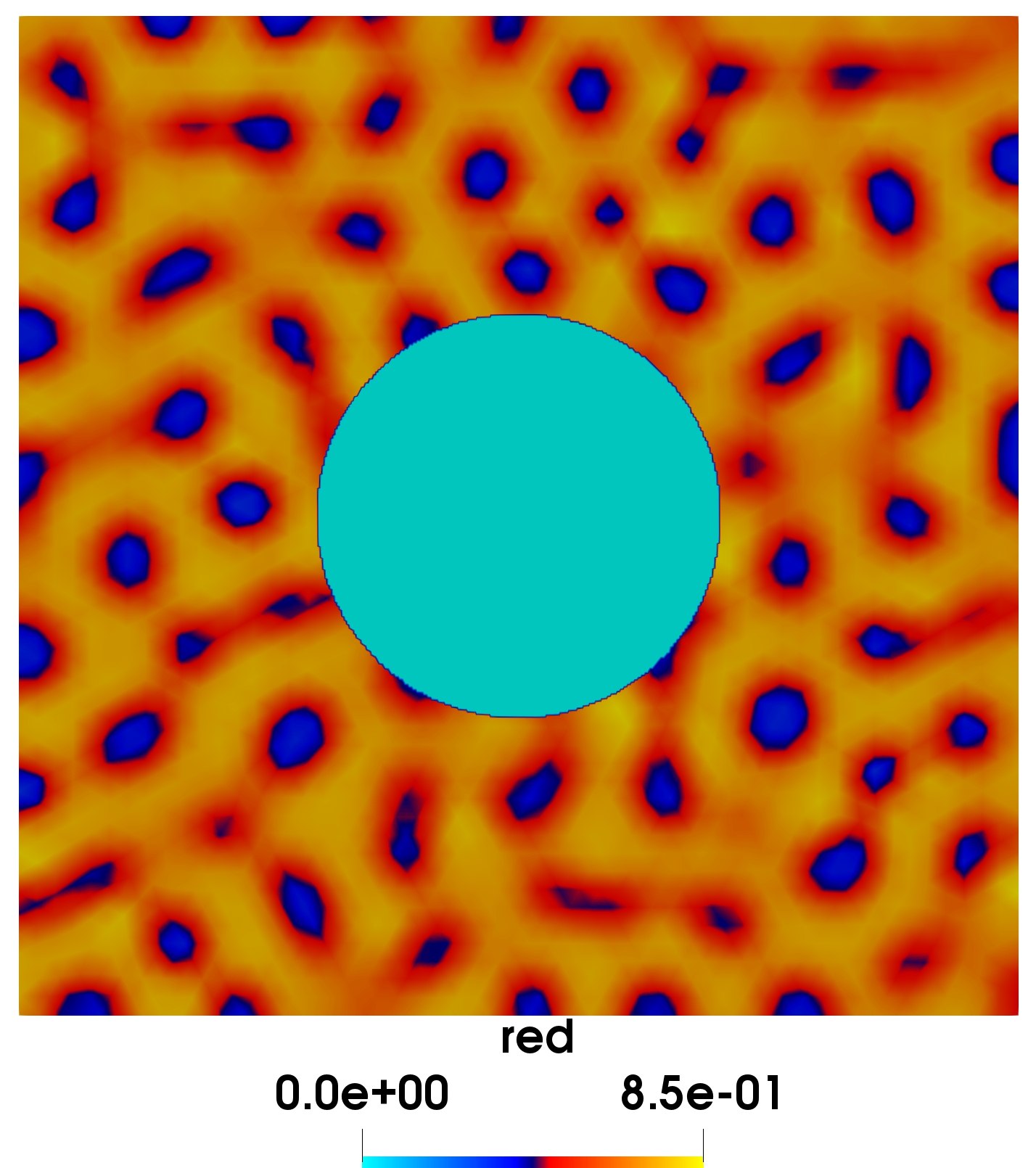}
\end{minipage}
\begin{minipage}{0.24\textwidth}
  \includegraphics[width=\textwidth]{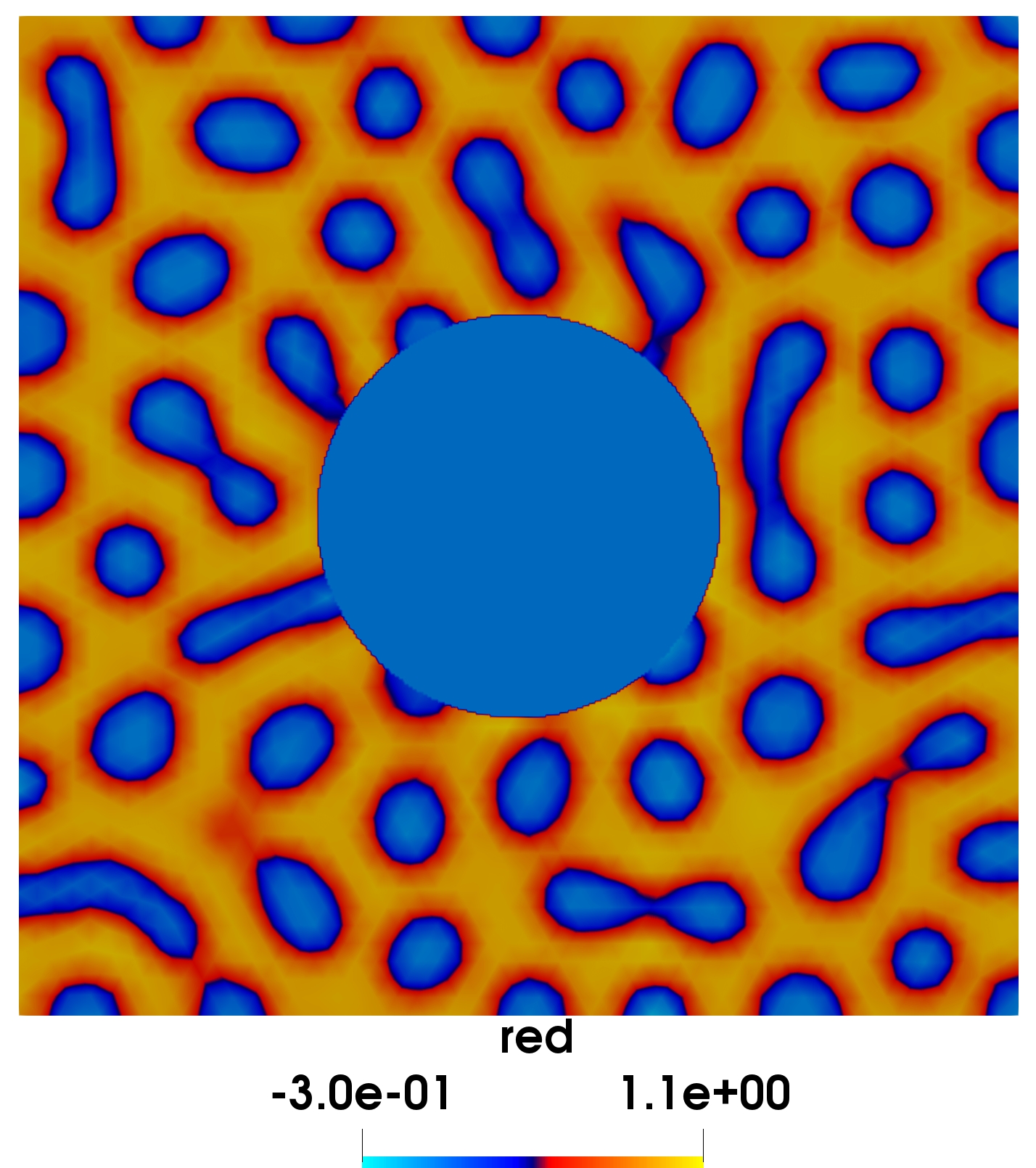}
\end{minipage}
\\
\begin{minipage}{0.24\textwidth}
  \includegraphics[width=\textwidth]{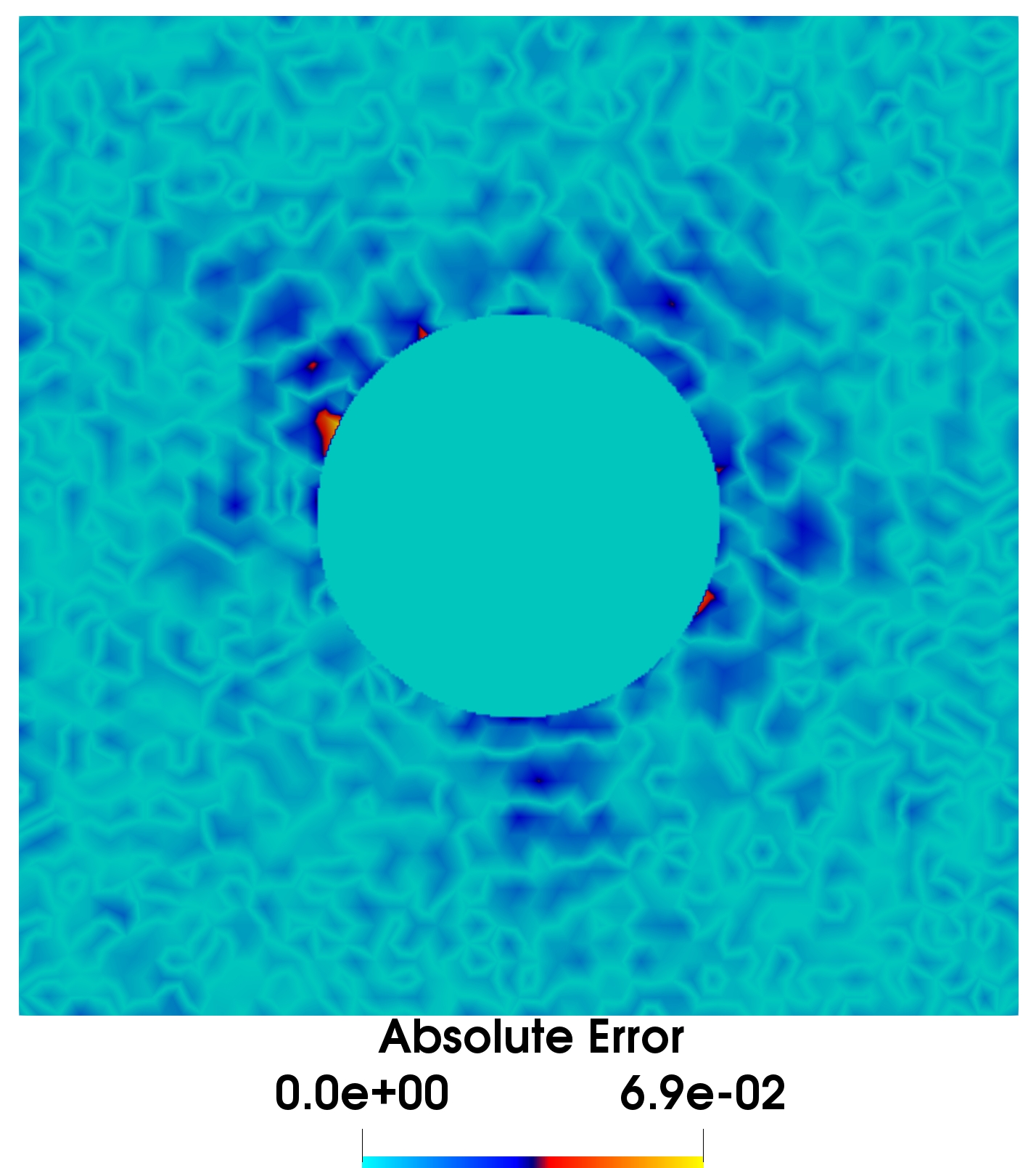}
\end{minipage}
\begin{minipage}{0.24\textwidth}
  \includegraphics[width=\textwidth]{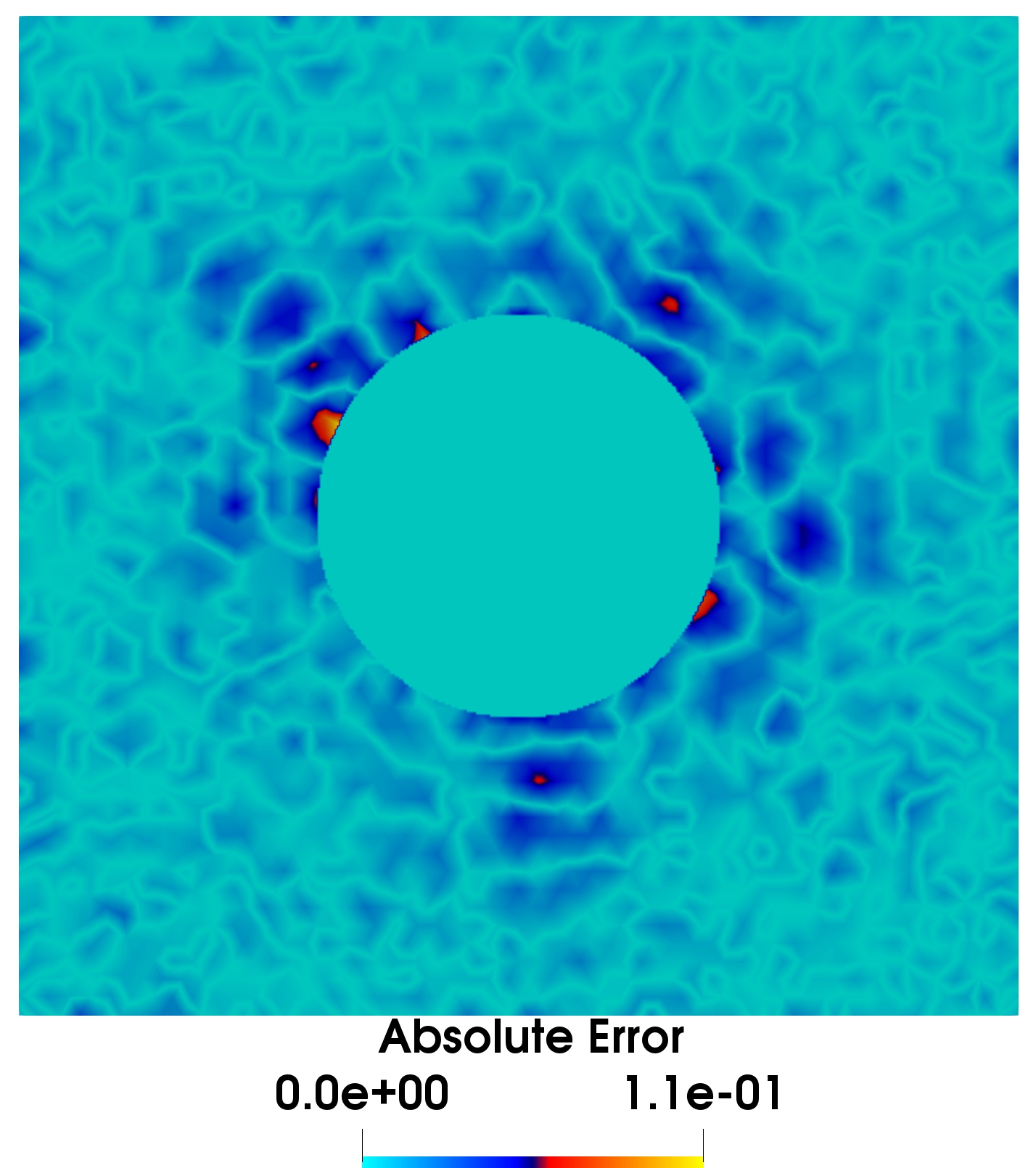}
\end{minipage}
\begin{minipage}{0.24\textwidth}
  \includegraphics[width=\textwidth]{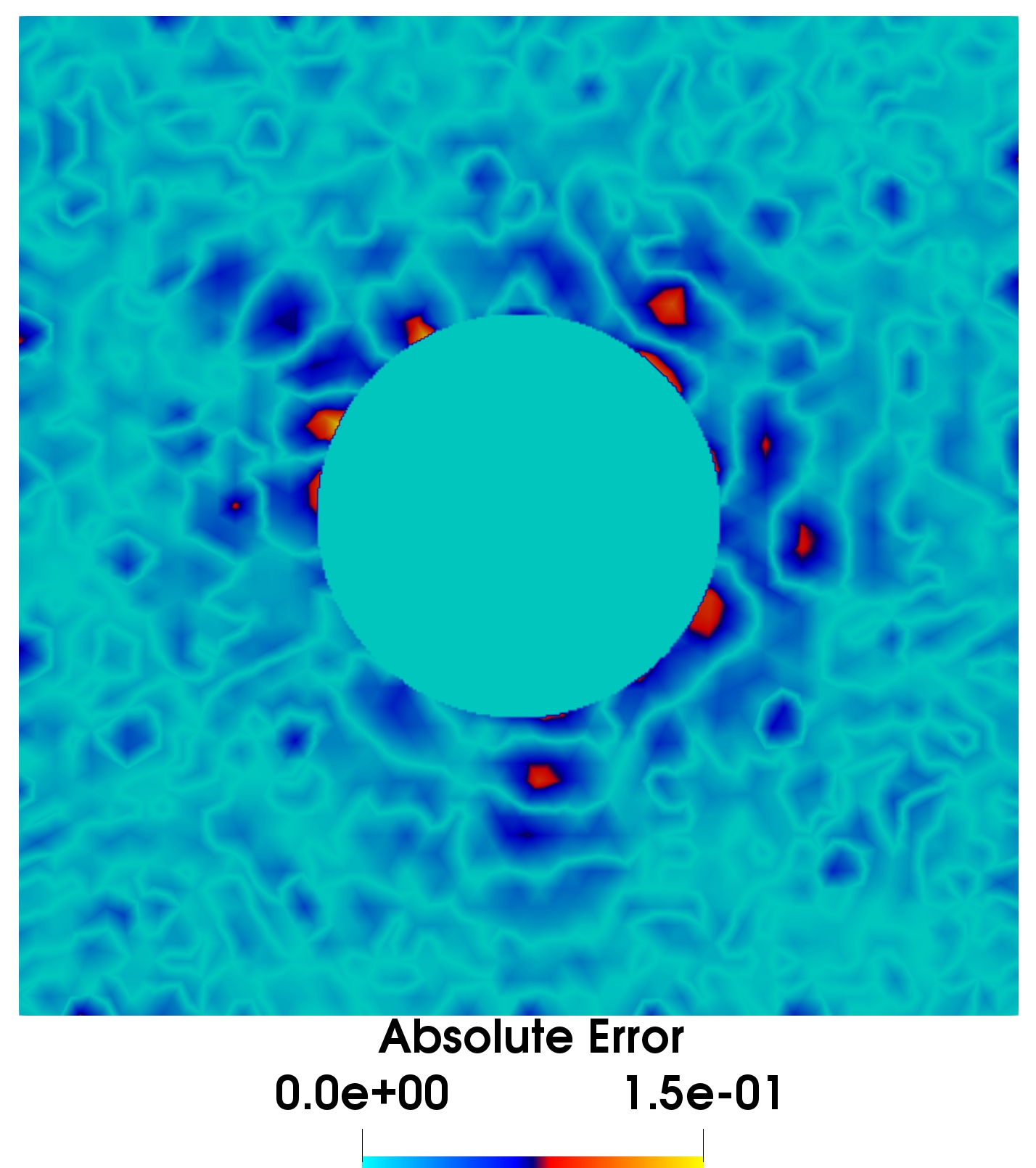}
\end{minipage}
\begin{minipage}{0.24\textwidth}
  \includegraphics[width=\textwidth]{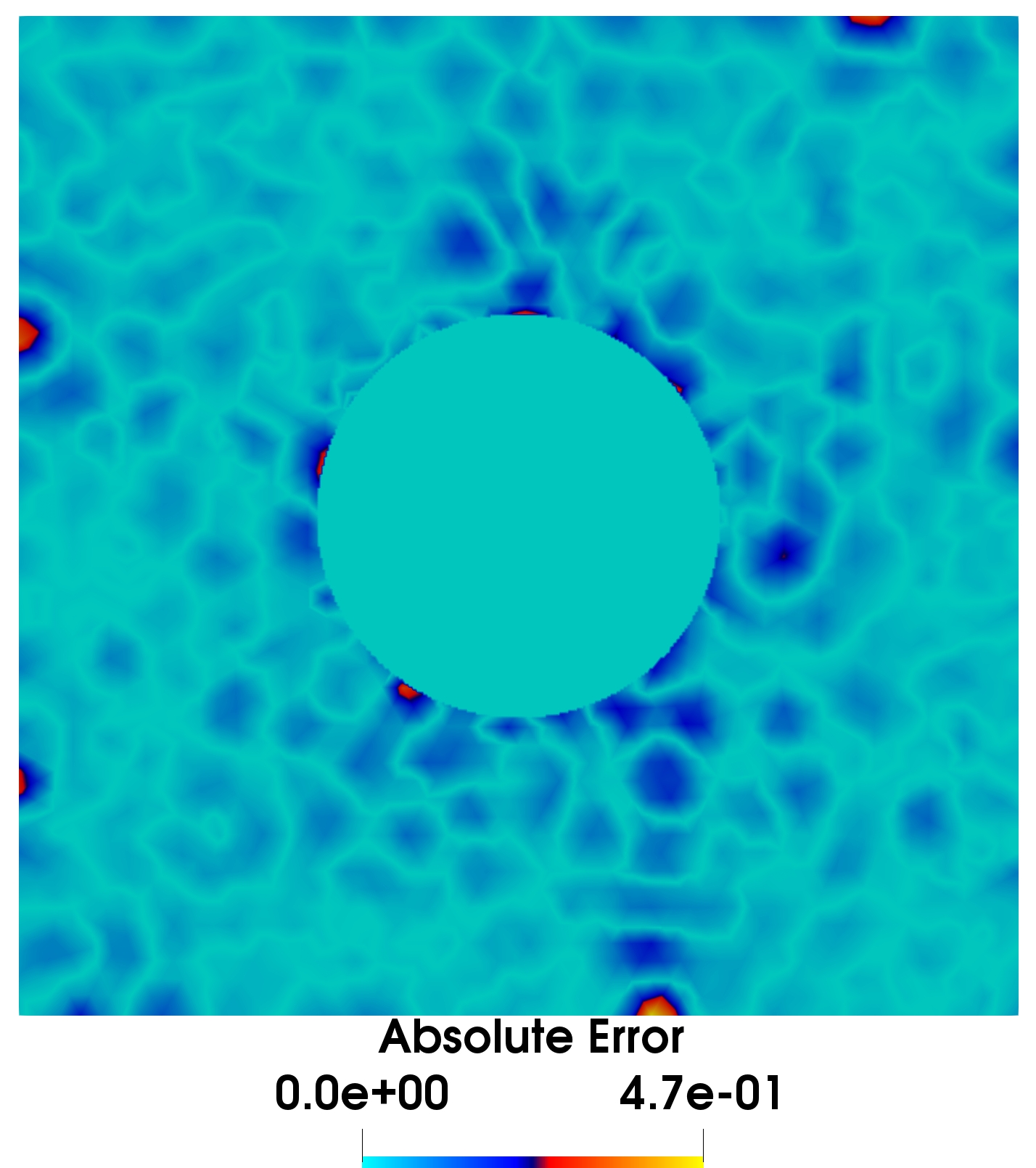}
\end{minipage}
\end{minipage}
\caption{Results for the embedded circle geometrical parametrization $\mu= 0.4261$, and $t=[28,36,46,100]dt$.
} 
\label{fig:CH_FULL_RED_ERROR_1P}
\end{figure} 
The full-order and the reduced mass evolution
 with respect to time verifies the conservation of mass 
for proper number of modes and for the  {truth} solver
and parameter $\mu=0.4261$. The experiment took place for time instances $t=ndt$, $n = 1,...100$, Figure \ref{Energy_Mass}. 
 \begin{figure} \centering
    \includegraphics[width=0.8\textwidth]{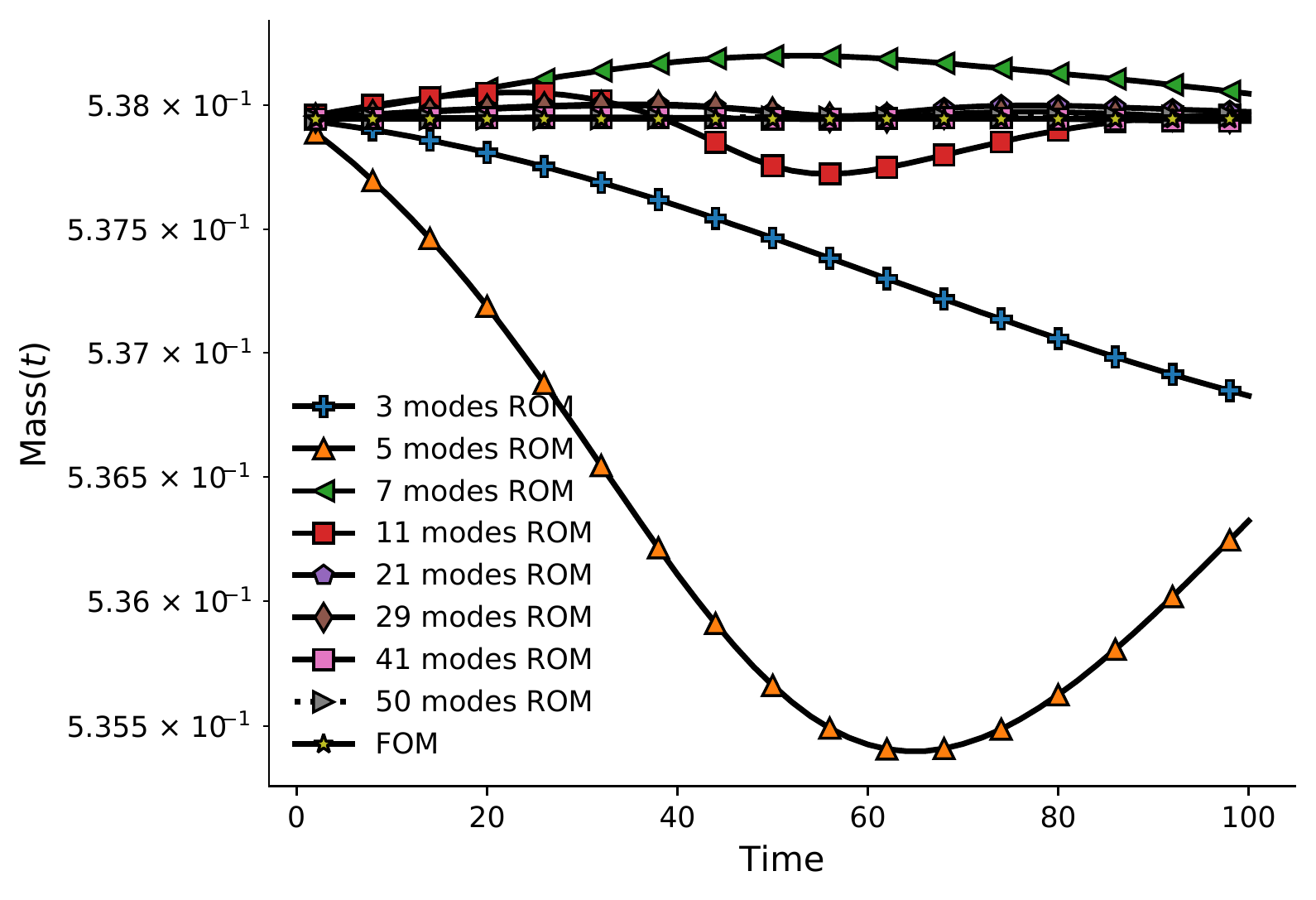} 
\caption{Cahn Hilliard (CutFEM): The FOM mass evolution with respect to time together with its RB approximation and the conservation of mass.}
  \label{Energy_Mass}
\end{figure}
We mention that we trained our basis with $900$ snapshots and in the parameter range for the  
diameter $\mu_{\text{test}} \in [0.36, 0.48]
$.
\section{Conclusions and future developments}\label{sec:conclusions}
In this chapter we introduced 
a POD-Galerkin ROM based on SBM/CutFEM high fidelity simulations, 
for  linear and nonlinear, steady and unsteady PDE problems, 
characterized by a 
geometrical parametrization with possibly large deformations.
The embedded boundary method discretization naturally allows to use a level set description of the parametrized geometry. In our opinion, 
 these results have been derived in a simpler and more versatile high fidelity and ROM method when compared to a FE formulation with pull back to a reference domain.
The transportation approach and the developed ROM is able to reproduce the high fidelity solution in an accurate manner, with relative errors of the order of $10^{-4}$, and 
for the Stokes case of the order of $10^{-3}$ employing only four modes for both velocity and pressure. 
Also, efficiency in nonlinear time depended systems, namely Cahn-Hilliard has been verified.

As perspectives we mention, as 
the proposed ROM 
is not offline-online separable in the usual sense, 
the empirical interpolation method,
and application of greedy algorithms during the generation of the reduced basis space.
As a further future development, we mention the 
extension of the proposed ROM into a more general framework of {nonlinear problems} in fluid dynamics, 
as well as fluid-structure interaction problems, 
multiphase flow and Navier-Stokes coupled systems,
 multiphysics, coupled processes or systems like heat transfer, stress and strain, optimal control in hydrodynamics, chemical reactions systems,
as well as shallow water flow systems.

\subsection{Acknowledgements}
We acknowledge the support by European Union Funding for Research and Innovation -- Horizon 2020 Program -- in the framework of European Research Council Executive Agency: Consolidator Grant H2020 ERC CoG 2015 AROMA-CFD project 681447 ``Advanced Reduced Order Methods with Applications in Computational Fluid Dynamics'' (PI Prof. Gianluigi Rozza).
We also acknowledge the INDAM-GNCS project ``Tecniche Numeriche Avanzate per Applicazioni Industriali''.
The first author has received funding from the Hellenic Foundation for Research and Innovation (HFRI) and  the  General  Secretariat  for  Research  and  Technology (GSRT), under  grant agreement No[1115], the ”First Call for H.F.R.I. Research Projects to support Faculty members and Researchers and the procurement of high-cost research equipment” grant 3270 and the National Infrastructures for Research and Technology S.A. (GRNET S.A.) in the National HPC facility - ARIS - under project ID pa190902.
Numerical simulations have been obtained, for the high fidelity solver with Nalu and Athena C++ Duke University in-home software, the extension \emph{ngsxfem} of \emph{ngsolve} software package, \cite{ngsolve,LHPW}, and for the reduced order part RBniCS, \cite{rbnics}.

\bibliographystyle{plain}
\bibliography{main_clean.bib}

\end{document}